\documentclass[11pt,reqno]{amsart}
\usepackage[T1]{fontenc}
\usepackage{lmodern}
\usepackage{amsmath,amssymb,amsthm,mathrsfs,mathtools}
\usepackage{booktabs,array,longtable}
\usepackage{xcolor}
\usepackage[expansion=false]{microtype}
\usepackage[a4paper,margin=30mm,headheight=14pt]{geometry}
\usepackage{hyperref}
\usepackage{orcidlink}
\definecolor{navy}{RGB}{18,48,96}
\colorlet{orcidlogocol}{black}

\definecolor{citationpurple}{RGB}{144,0,128}
\newif\ifshowcitationrevisions
\showcitationrevisionsfalse 
\newenvironment{citationrevision}{\par\begingroup\ifshowcitationrevisions\color{citationpurple}\fi}{\par\endgroup}
\newcommand{\citationadded}[1]{{\ifshowcitationrevisions\color{citationpurple}\fi#1}}
\hypersetup{colorlinks=true,
    linkcolor=black,
    citecolor=black,
    urlcolor=black,
bookmarksnumbered=true,pdftitle={Positive solutions of power equations in planar corner domains},pdfauthor={Hongtao Hu; Zhouyu Long; Hongzhe Tao; Qiaoran Wu; Wenming Zou},pdfsubject={Semilinear elliptic equations and Mazya Problem 50},pdfkeywords={corner domains, Martin trace, Hessian estimates, singular elliptic equations},bookmarksdepth=2}
\numberwithin{equation}{section}
\newtheorem{theorem}{Theorem}[section]
\newtheorem{proposition}[theorem]{Proposition}
\newtheorem{lemma}[theorem]{Lemma}
\newtheorem{corollary}[theorem]{Corollary}
\theoremstyle{definition}
\newtheorem{definition}[theorem]{Definition}
\theoremstyle{remark}
\newtheorem{remark}[theorem]{Remark}
\newcommand{\R}{\mathbb R}

\newcommand{\Om}{\Omega}
\newcommand{\dd}{\,\mathrm d}
\newcommand{\G}{\mathcal G}
\newcommand{\HH}{\mathcal H}
\newcommand{\cA}{\mathcal A}
\newcommand{\eps}{\varepsilon}
\newcommand{\tr}{\operatorname{tr}}
\newcommand{\dist}{\operatorname{dist}}
\newcommand{\supp}{\operatorname{supp}}

\newcommand{\one}{\mathbf 1}

\newcommand{\norm}[1]{\lVert #1\rVert}
\newcommand{\abs}[1]{\lvert #1\rvert}
\newcommand{\weak}{\rightharpoonup}

\newcommand{\Hmat}{\mathsf H}
\newcommand{\Dmat}{\mathsf D}
\newcommand{\diag}{\operatorname{diag}}
\newcommand{\Span}{\operatorname{span}}

\newcommand{\orcid}[1]{\textsuperscript{\orcidlink{#1}}}
\allowdisplaybreaks[2]
\title[Power equations in planar corner domains]{Positive solutions of power equations\\ in planar corner domains}

\author[Hongtao Hu]{Hongtao Hu\orcid{0009-0000-5966-3130}}
\address{Department of Mathematical Sciences, Tsinghua University, Beijing 100084, China}
\email{\href{mailto:huht23@mails.tsinghua.edu.cn}{\texttt{huht23@mails.tsinghua.edu.cn}}}

\author[Zhouyu Long]{Zhouyu Long}
\address{Department of Mathematical Sciences, Tsinghua University, Beijing 100084, China}
\email{\href{mailto:longzy25@mails.tsinghua.edu.cn}{\texttt{longzy25@mails.tsinghua.edu.cn}}}

\author[Hongzhe Tao]{Hongzhe Tao}
\address{Department of Mathematical Sciences, Tsinghua University, Beijing 100084, China}
\email{\href{mailto:taohz25@mails.tsinghua.edu.cn}{\texttt{taohz25@mails.tsinghua.edu.cn}}}

\author[Qiaoran Wu]{Qiaoran Wu}
\address{Department of Mathematical Sciences, Tsinghua University, Beijing 100084, China}
\email{\href{mailto:wuqr24@mails.tsinghua.edu.cn}{\texttt{wuqr24@mails.tsinghua.edu.cn}}}

\author[Wenming Zou]{Wenming Zou \orcid{0000-0002-4178-9568}}
\address{Department of Mathematical Sciences, Tsinghua University, Beijing 100084, China}
\email{\href{mailto:zou-wm@mail.tsinghua.edu.cn}{\texttt{zou-wm@mail.tsinghua.edu.cn}}}

\date{}
\begin{document}
\begin{abstract}
We study positive solutions of $\Delta u=\lambda u^k$ on bounded planar
corner domains with zero boundary values away from the vertex. For
negative-power absorption we identify the admissible interval of Martin
masses and prove uniqueness at the minimum mass. Along the ordinary boundary,
we identify the weak Hessian endpoint and its sharp tail when the set of
zero inward derivative has positive length. For superlinear sources we prove
attainment of the largest positive Martin mass in the subcritical range.
At criticality we construct a family of singular solutions with zero
Martin mass and distinguish its members by their second logarithmic
coefficient and the endpoint integrability of their Hessian differences.
Above criticality an angular eigenvalue governs the Morse index away from
its zero threshold. We also give sharp physical-Hessian asymptotics for
finite-mass remainders and joint corner and side estimates for negative-power
sources.
\end{abstract}
\maketitle

\tableofcontents

\section{Introduction}\label{sec:intro}
Maz'ya's Problem~50~\cite[Section~6.4]{Mazya} asks for existence, uniqueness, and regularity of positive solutions of
\begin{equation}\label{eq:original}
 \Delta u-\lambda u^k=0\quad\text{in }\Om,\qquad
 u=0\quad\text{on }\partial\Om\setminus\{O\},\qquad\lambda,k\in\R,
\end{equation}
where $\Om$ is a bounded planar domain, smooth except at an angular point $O$. We assume that the incident arcs are $C^2$ up to the vertex and make the angle $\theta\in(0,2\pi)$. We use the natural class of positive classical interior solutions, continuous with the stated boundary values away from $O$. Neither a value nor a singular coefficient is prescribed at $O$.

The source of this question is Section~6.4 of Maz'ya's 2018 collection~\cite{Mazya}, pp.~26--27, where the harmonic problem on a truncated sector is given as the elementary model. Related earlier work includes the study of quasilinear boundary problems near a conical point by Maz'ya and Plamenevskii~\cite{MP73}. The free singular amplitude in the harmonic model already explains why zero values on the punctured boundary do not determine a unique solution. The nonlinear question concerns how the reaction term changes the admissible singularities, the solution families, and their regularity. We correct the sector formula and separate its local geometry from the one-corner domain assumption in Section~\ref{subsec:sector-model}. Section~6.4 itself specifies no physical derivation; the connections below describe interpretations of the equation and its boundary data.

The missing corner datum is essential for uniqueness. Even the harmonic problem has the whole ray $cK$, where $K$ is the positive Martin kernel at the vertex. For a fixed nonlinear equation with $\lambda\ne0$, a nontrivial constant rescaling changes the equation, so the construction of a family requires a separate argument. Conversely, fixing a Martin coefficient can give uniqueness for a negative-power source or a sublinear source, but does not give a general uniqueness principle. For a superlinear source even a positive prescribed coefficient can have more than one solution. At the first linear eigenvalue the coefficient is zero for an entire eigenfunction ray.

The purpose of this paper is to give quantitative results for these distinctions and for second derivatives at curved corners. The exact parameter range for existence is summarized in Theorem~\ref{thm:existencetable}. In that statement a solution is permitted to be bounded at $O$. Write $\beta=\pi/\theta$ and define the harmonic second-order threshold
\begin{equation}\label{eq:harmonicthreshold}
 P_\beta=
 \begin{cases}2/(2-\beta),&1/2<\beta<2,\quad\beta\ne1,\\
 \infty,&\beta=1\text{ or }\beta\ge2.
 \end{cases}
\end{equation}

The admissible positive Martin masses form a closed upper interval for
negative-power absorption and a bounded interval with an attained upper
endpoint for subcritical superlinear sources. For the critical source
equation we construct a global family of singular solutions of zero mass
and separate its members using the solution-dependent constant $c_u$
in the second logarithmic expansion. The same constant characterizes
the endpoint integrability of the Hessian of the difference of two
critical singular solutions. The ordinary-boundary endpoint for negative-power absorption and the sharp
physical-Hessian asymptotics at the vertex give the corresponding
second-order regularity results.

\subsection{Notation}
Throughout the paper write
\begin{equation}\label{eq:geometry}
 \beta=\frac\pi\theta,\qquad I=(0,\theta),\qquad
 \psi(\varphi)=\sin(\beta\varphi),\qquad N=\int_I\psi^2\dd\varphi=\frac\theta2\,,
\end{equation}
and $d(x)=\dist(x,\partial\Om)$. A natural solution is a function
$u\in C^2(\Om)\cap C(\overline\Om\setminus\{O\})$, positive in $\Om$ and zero on $\partial\Om\setminus\{O\}$, satisfying the equation pointwise. In particular, no boundedness, energy, stability, or trace condition at $O$ is included.

Let $G(x,y)$ be the positive Dirichlet Green kernel of $-\Delta$, and define
\[
 \G[f](x)=\int_\Om G(x,y)f(y)\dd y
\]
whenever the integral is finite. Fix $x_*\in\Om$ and normalize the positive Martin kernel at $O$ by $K(x_*)=1$. After translating $O$ to the origin and rotating, fix a conformal map $F$ from a truncated sector into a neighborhood of $O$ in $\Om$, taking its straight sides to the incident arcs, with $F(0)=0$ and $F'(0)=1$. Its construction and expansion are given in Section~\ref{sec:linear}. The positive number $\kappa$ is the coefficient in the leading relation
\[
 K(F(re^{i\varphi}))\sim\kappa r^{-\beta}\psi(\varphi).
\]
We write $f=O_\psi^1(h)$ for a positive weight $h$ if $f/\psi$, $\partial_{\log r}f$, and $\partial_\varphi f$ are all $O(h)$. Matrix norms are Frobenius norms.

For a natural solution of~\eqref{eq:original}, let
\begin{equation*}
 Q_u[\xi]=\int_\Om\bigl(|\nabla\xi|^2+\lambda k u^{k-1}\xi^2\bigr)\dd x
 \qquad(\xi\in C_c^\infty(\Om))\,.
\end{equation*}
Stability on an open subset means that $Q_u$ is nonnegative on tests supported there. Strict stability means a positive lower bound by a multiple of the Dirichlet energy. The Morse index is the supremum of the dimensions of subspaces of $C_c^\infty(\Om)$ on which $Q_u$ is negative definite.

\subsection{Existence}
\begin{theorem}\label{thm:existencetable}
For the natural positive problem $\Delta u=\lambda u^k$ with zero ordinary boundary data, a positive solution exists if and only if one of the following four conditions holds:
\begin{enumerate}
    \item \(\lambda = 0\) and \(k \in \mathbb{R}\);
    \item\(\lambda < 0\) and \(k \in \mathbb{R} \setminus \{1\}\);
    \item \(\lambda < 0\), \(k = 1\), and \(-\lambda \le \Lambda_1\);
    \item\(\lambda > 0\) and \(-1 < k < 1 + \frac{2}{\beta}\).
\end{enumerate}
Here $\Lambda_1$ denotes the first Dirichlet eigenvalue of $-\Delta$ on $\Om$ and $\beta$ is defined by \eqref{eq:geometry}. 
\end{theorem}

\subsection{Trace, uniqueness, and classification}
\subsubsection{Negative-power absorption}
For $0<p<1$, the negative-power absorption equation and its boundary exponents are
\begin{equation}\label{eq:negabs}
 \Delta u=a u^{-p},\qquad a>0,\quad 0<p<1\,,
\end{equation}
\begin{equation}\label{eq:gammaabs}
 \gamma=\frac2{1+p},\qquad \alpha=\gamma-1,\qquad
 m=2-\gamma=p\gamma,\qquad s_b=\frac1m=\frac{1+p}{2p}\,.
\end{equation}

\begin{theorem}\label{thm:finiteabs}
Every natural solution of~\eqref{eq:negabs} satisfies $u^{-p}\in L^1(\Om)$. There is a unique finite number $c>0$ such that
\begin{equation}\label{eq:repabs}
 u=cK-a\G[u^{-p}],\qquad
 \frac{u(x)}{K(x)}\longrightarrow c\quad(x\to O)\,.
\end{equation}
The limit holds uniformly along all interior paths.
\end{theorem}

For the quantitative construction put $H=\G[K^{-p}]$, $M=\sup_\Om(H/K)$, and
\begin{equation}\label{eq:Jgain}
 J=\G[K^{-p-1}H],\qquad \vartheta=\frac1M\sup_\Om\frac JH\,.
\end{equation}
Section~\ref{sec:negabs} proves that these quantities are finite and $0<\vartheta<1$. Define
\begin{align}
 s_p&=\frac p{p+1},& F_p&=s_p^{-p},& B_p&=1+\vartheta(F_p-1),\label{eq:twoconst}\\
 C_G&=\big[aM(p+1)B_p\big]^{1/(p+1)},&
 q_G&=\frac{\vartheta F_p}{B_p}<1,&
 C_0&=(p+1)p^{-p/(p+1)}(aM)^{1/(p+1)}.\nonumber
\end{align}

\begin{theorem}\label{thm:twoGreen}
For each $c\ge C_G$ there is a greatest natural solution $U_c$ among solutions with trace $c$, and
\begin{equation}\label{eq:shaped}
 cK-\frac{cH}{(p+1)M}\le U_c<cK,\qquad U_c\ge s_pcK\,.
\end{equation}
It is the only solution with that trace satisfying $u\ge s_pcK$. If $\tau_c=ac^{-p-1}$ and $q_c=p\tau_c\vartheta M s_p^{-p-1}$, the iteration
\begin{equation}\label{eq:iterate}
 z_0=0,\qquad z_{n+1}=\tau_c\G[(K-z_n)^{-p}]
\end{equation}
increases to $z=K-U_c/c$ and satisfies
\begin{equation}\label{eq:error}
 0\le c(K-z_n)-U_c\le ac^{-p}\frac{q_c^n}{1-q_c}H,\qquad q_c\le q_G<1\,.
\end{equation}
If $c_2>c_1\ge C_G$, then $U_{c_2}-U_{c_1}>(c_2-c_1)K$. In particular the original data admit uncountably many solutions.
\end{theorem}

The preceding construction gives an explicit upper bound for the minimum
admissible trace. The following theorem describes the entire maximal
branch and singles out the solution at its endpoint.

\begin{theorem}\label{thm:inverse-mass-branch}
For the natural solutions of~\eqref{eq:negabs}, the set of Martin
traces is exactly $[c_*,\infty)$, where
\begin{equation}\label{eq:inverse-mass-bounds}
 aM\le c_*^{1+p}\le C_G^{1+p}.
\end{equation}
For every $c\ge c_*$ there is a greatest solution $U_c$ of trace $c$.
A solution $u$ of trace $c$ equals $U_c$ if and only if it satisfies
\begin{equation}\label{eq:inverse-mass-stability}
 ap\int_\Om u^{-p-1}\xi^2\dd x
 \le\int_\Om|\nabla\xi|^2\dd x
 \qquad (\xi\in C_c^\infty(\Om))\,.
\end{equation}
The branch satisfies $U_{c'}\ge U_c+(c'-c)K$ for $c'>c\ge c_*$,
is pointwise strictly concave on $[c_*,\infty)$, and is smooth for
$c>c_*$ in the norm $\norm{v}_K=\sup_\Om|v|/K$.
The minimum-trace solution $U_{c_*}$ is unique among all natural
solutions of trace $c_*$. For every $E\Subset\Om$ and $L>c_*$,
every natural solution $u$ of trace $c\in[c_*,L]$ satisfies
\begin{equation}\label{eq:inverse-mass-local-rate}
 \int_E |u-U_{c_*}|^2\dd x\le C_{E,L}(c-c_*).
\end{equation}
On each interior compact subset, the $C^{2,\eta}$ norm of the
difference is $O((c-c_*)^{1/2})$ as $c\downarrow c_*$, for
every $0<\eta<1$.
\end{theorem}

\subsubsection{Negative-power sources and sublinear equations}
For negative-power sources write
\begin{equation}\label{eq:negsource}
 -\Delta u=au^{-q},\qquad a>0,\quad q>0\,.
\end{equation}

\begin{proposition}\label{prop:negsourceclass}
For every $c\ge0$ there is exactly one natural solution $u_c$ of~\eqref{eq:negsource} with
\begin{equation}\label{eq:negsourceclass}
 u_c=cK+a\G[u_c^{-q}],\qquad
 \max\{u_0,cK\}\le u_c\le u_0+cK\,.
\end{equation}
These are all the natural solutions. The solution $u_0$ is bounded and continuous with zero values on the whole boundary, and $u_c-cK\to0$ absolutely on the whole boundary. One has $u_0\in H^1_0(\Om)$ if and only if $q<3$, whereas $u_c\notin H^1(\Om)$ for $c>0$. Thus the original data, which do not prescribe $c$, are nonunique for every $q>0$.
\end{proposition}

\begin{theorem}\label{thm:subclass}
Let $a>0$ and $0<p<1$. For $-\Delta u=au^p$ there is exactly one positive solution $u_c$ of
\begin{equation}\label{eq:subsource}
 u_c=cK+a\G[u_c^p]
\end{equation}
for every $c\ge0$. These exhaust the natural class, increase strictly with $c$, and satisfy $u_c/K\to c$ uniformly at the corner. The solution $u_0$ is bounded, belongs to $H^1_0(\Om)$, and is zero on the whole boundary.

For $\Delta u=au^p$ and every $c>0$ there is exactly one nonnegative solution $v_c$ with
\begin{equation}\label{eq:subabs}
 v_c=cK-a\G[v_c^p],\qquad v_c/K\to c
\end{equation}
uniformly at the corner. Every positive natural solution is one of these. All sufficiently large $c$ give a strictly positive solution; sufficiently small $c$ give a solution with an interior zero set. For either sign there are infinitely many solutions of the original strictly positive problem.
\end{theorem}

\subsubsection{Superlinear equations}
For $p>1$ define the critical corner exponent
\begin{equation}\label{eq:criticalpower}
 p_c=1+\frac2\beta=1+\frac{2\theta}{\pi}\,.
\end{equation}

\begin{proposition}\label{prop:superabs}
For $\Delta u=au^p$, $a>0$, natural positive solutions exist exactly when $1<p<p_c$. In this range they are the increasing finite-trace family $u_c$, $c>0$, and its limit $u_\infty$. A finite-trace solution is unique for each prescribed $c$, and there is exactly one strongly singular solution $u_\infty$. Their corner behaviors are
\begin{equation}\label{eq:superabsprofiles}
 u_c/(cK)\to1,\qquad
 u_\infty\circ F\sim r^{-d_p}\Psi_p(\varphi),\quad d_p=\frac2{p-1}>\beta,
\end{equation}
where $\Psi_p$ is the unique positive zero-endpoint solution of
\begin{equation}\label{eq:superabsangular}
 \Psi_p''+d_p^2\Psi_p=a\Psi_p^p\,.
\end{equation}
Both equivalences are relative and uniform up to the sides. Thus the original data are nonunique throughout the existence range.
\end{proposition}

\begin{theorem}\label{thm:positive-source-mass}
Let $a>0$ and $p>1$, and let $c$ be the coefficient of $K$ in the greatest harmonic minorant of a natural solution of $-\Delta u=au^p$. A solution with $c>0$ exists if and only if $p<p_c=1+2/\beta$. In that range, for some $0<c^*<\infty$, the set of realized positive masses is $(0,c^*]$. For each such $c$ there is a pointwise least solution $U_c$ of mass $c$. These solutions are increasing in $c$ and $U_c\uparrow U_{c^*}$ as $c\uparrow c^*$. They satisfy
\begin{equation}\label{eq:positive-source-stability}
 ap\int_\Om U_c^{p-1}\zeta^2\,\dd x
 \le \int_\Om |\nabla\zeta|^2\,\dd x
 \qquad (\zeta\in C_c^\infty(\Om)).
\end{equation}
For sufficiently small $c>0$, one has $cK\le U_c\le 2cK$ and $U_c$ is the only solution in this order interval.
\end{theorem}

\begin{proposition}\label{prop:supersourceexist}
For $-\Delta u=au^p$ with $a>0$ and every finite $p>1$, there is a bounded positive solution in $H^1_0(\Om)$. If $1<p<p_c$, then for every sufficiently small $c>0$ there are at least two natural solutions with the same positive trace $c$. There is a least such solution, and it is strictly stable for sufficiently small $c$.
\end{proposition}

\begin{theorem}\label{thm:sourcefirst}
Let $u$ be any natural solution of $-\Delta u=au^p$, with $a>0$ and $p>1$, and put $d_p=2/(p-1)$. The following are all the possible leading behaviors. In each formula the relative equivalence is uniform in the open angle, including approaches to the sides.
\begin{enumerate}
\item If $1<p<p_c$, then either $u\circ F\sim B r^\beta\psi$, $B>0$, or $u\circ F\sim c\kappa r^{-\beta}\psi$, $c>0$.
\item If $p=p_c$, then either the same regular branch occurs, or
\begin{equation}\label{eq:criticalfirst}
 u\circ F\sim C_*r^{-\beta}(\log(1/r))^{-\mu}\psi,\qquad
 \mu=\beta/2,\quad
 C_*=\left(\frac{\beta^2N}{a\int_I\psi^{p+1}}\right)^\mu\,.
\end{equation}
\item If $p>p_c$, then either the regular branch occurs, or $u\circ F\sim r^{-d_p}\Phi_p$, where $\Phi_p$ is the unique positive solution of
\begin{equation}\label{eq:sourceangular}
 -\Phi_p''-d_p^2\Phi_p=a\Phi_p^p,\qquad\Phi_p(0)=\Phi_p(\theta)=0\,.
\end{equation}
\end{enumerate}
The regular branch has finite Dirichlet energy and vanishes continuously at $O$. Every singular branch has infinite energy. The Riesz--Martin coefficient is positive only in the subcritical weak branch; critical and supercritical singular solutions have coefficient zero and satisfy $u=a\G[u^p]$.
\end{theorem}

\begin{proposition}\label{prop:source-angular-stability}
Let $p>p_c$ and let $u$ have the strong singular asymptotic in
Theorem~\ref{thm:sourcefirst}. Let $\nu_1$ be the first Dirichlet eigenvalue
on $I=(0,\theta)$ of
\begin{equation}\label{eq:source-angular-linearized}
 \mathcal A_p=-\partial_\varphi^2-ap\Phi_p^{p-1}.
\end{equation}
If $\nu_1>0$, then $u$ is stable near $O$ and has finite Morse index. If
$\nu_1<0$, then $u$ has infinite Morse index. At fixed angle, $\nu_1>0$
for $p$ sufficiently close to $p_c$ from above, while
$p d_p^2<(p-1)\beta^2$ implies $\nu_1<0$.
\end{proposition}

\subsection{Regularity and sharp asymptotics}
\subsubsection{Ordinary-boundary regularity for negative-power absorption}
\begin{theorem}\label{thm:C1abs}
Every natural solution of~\eqref{eq:negabs} belongs to $C^{1,\alpha}_{\mathrm{loc}}(\overline\Om\setminus\{O\})$. In particular $b(P)=\partial_\nu u(P)\ge0$ exists at every ordinary boundary point and is locally $C^{0,\alpha}$. On a uniform ordinary boundary piece write $D_R(P)=\Om\cap B_R(P)$ and $M=\sup_{D_{8R}(P)}u$. Then
\begin{equation}\label{eq:C1abs}
 R\norm{\nabla u}_{L^\infty(D_R)}+R^\gamma[\nabla u]_{C^{0,\alpha}(\overline D_R)}
 \le C\big(M+a^{1/(1+p)}R^\gamma\big)\,.
\end{equation}
If $b(P)=0$, in coordinates with $\nu(P)=e_2$,
\begin{equation}\label{eq:absprofile}
 r^{-\gamma}u(P+r\cdot)\longrightarrow A_p(y_2)_+^\gamma,
 \qquad A_p=\left(\frac{a}{\gamma(\gamma-1)}\right)^{1/(1+p)},
\end{equation}
locally in $C^1$ after identifying the true domains by boundary graphs, and in $C^2$ on compact subsets of $\{y_2>0\}$. Thus the exponent $\alpha$ is optimal at every such point.
\end{theorem}

Let $\Gamma$ be a compact ordinary boundary arc. Choose $\delta>0$ within the tubular radius of a slightly larger ordinary arc, with its closed collar disjoint from $O$ and the other boundary pieces. Then the inward normal coordinates are injective and $d(P+t\nu(P))=t$ for $P\in\Gamma$ and $0<t\le\delta$. Put
\[
 V_\delta(\Gamma)=\{P+t\nu(P):P\in\Gamma,\ 0<t<\delta\},\qquad
 Z=\{P\in\partial\Om\setminus\{O\}:b(P)=0\}.
\]
We denote one-dimensional Hausdorff measure by $\HH^1$.

\begin{theorem}\label{thm:mixed}
For every natural solution of~\eqref{eq:negabs},
\begin{equation}\label{eq:mixedspace}
 u^{-p},\ |D^2u|\in L^{s_b,\infty}_{\mathrm{loc}}(\overline\Om\setminus\{O\}),
 \qquad u\in W^{2,s}_{\mathrm{loc}}(\overline\Om\setminus\{O\}),\quad1\le s<s_b\,.
\end{equation}
More precisely,
\begin{align}
 \lim_{t\downarrow0}t^{-1}\abs{\{u<t^\gamma\}\cap V_\delta(\Gamma)}
 &=A_p^{-1/\gamma}\HH^1(Z\cap\Gamma),\label{eq:mixedlow}\\
 \lim_{L\to\infty}L^{s_b}\abs{\{|D^2u|>L\}\cap V_\delta(\Gamma)}
 &=(aA_p^{-p})^{s_b}\HH^1(Z\cap\Gamma).\label{eq:mixedtail}
\end{align}
The same limit as in~\eqref{eq:mixedtail} holds for $au^{-p}$. No stability or lower bound on the positive values of $b$ is assumed. If $\HH^1(Z\cap\Gamma)>0$, the weak endpoint is sharp and $D^2u\notin L^{s_b,\ell}$ for every finite $\ell$.
\end{theorem}

\subsubsection{Finite-trace Green remainders}
For complete conformal annuli put $\cA=\{z:1<|z|<2,\ 0<\arg z<\theta\}$. For an angular function $Y$, the polar-frame Hessian and logarithmic-derivative matrices are
\begin{equation}\label{eq:Hmat}
 \Hmat_\lambda[Y]=
 \begin{pmatrix}\lambda(\lambda-1)Y&(\lambda-1)Y'\\(\lambda-1)Y'&Y''+\lambda Y\end{pmatrix},\qquad
 \Dmat_\lambda[Y]=
 \begin{pmatrix}(2\lambda-1)Y&Y'\\Y'&Y\end{pmatrix}.
\end{equation}
All matrix norms below are Frobenius norms.

\begin{theorem}\label{thm:remainder}
Suppose $a>0$, $p>0$, $p\ne1$, and $\beta(p-1)<2$. Let $u$ be a solution, positive near $O$, with
\begin{equation}\label{eq:remainderhyp}
 u=cK+\epsilon w,\quad w=a\G[u^p]>0,\quad c>0,\quad
 \epsilon\in\{-1,1\},\qquad u/(cK)\longrightarrow1
\end{equation}
uniformly along all paths to the corner. Put
\begin{gather}
 \sigma=2-\beta p,\qquad Q=a(c\kappa)^p,\nonumber\\
 P_\beta=
 \begin{cases}2/(2-\beta),&1/2<\beta<2,\ \beta\ne1,\\
 \infty,&\beta=1\text{ or }\beta\ge2,
 \end{cases}
 \qquad S=\min\{P_\beta,2/(\beta p)\}\,.\label{eq:remainderindices}
\end{gather}
For every finite $s\ge1$ and sufficiently small $\delta$,
\begin{equation}\label{eq:remainderW}
 w\in W^{2,s}(\Om\cap B_\delta(O))\quad\Longleftrightarrow\quad s<S\,.
\end{equation}
In particular $S\le1$ implies $w\notin W^{2,1}$.

Except when $\beta=1$ and $0<p<1$, the function and complete annular physical Hessian have the leading profiles
\begin{equation}\label{eq:remainderprofiles}
\begin{array}{c|c|c}
 \text{condition}&w\circ F&\text{normalized Hessian limit}\\ \hline
 \beta<\sigma&B r^\beta\psi&r^{2-\beta}D_x^2w(F(r\cdot))\to B D_z^2(|z|^\beta\psi)\\
 \beta>\sigma&Q r^\sigma Y_\sigma&r^{2-\sigma}D_x^2w(F(r\cdot))\to QD_z^2(|z|^\sigma Y_\sigma)\\
 \beta=\sigma&D r^\beta\log(1/r)\psi&\frac{r^{2-\beta}}{\log(1/r)}D_x^2w(F(r\cdot))\to DD_z^2(|z|^\beta\psi)
\end{array}
\end{equation}
The limits hold in $L^\ell(\cA)$ for every finite $1<\ell<\infty$. Here $B>0$, $Y_\sigma$ is the positive zero-endpoint solution of $-Y''-\sigma^2Y=\psi^p$ in the second row, and
\begin{equation}\label{eq:Dres}
 D=\frac{Q}{2\beta N}\int_I\psi^{p+1}\dd\varphi>0\,.
\end{equation}
When $\beta=1$ and $0<p<1$, remove the physical harmonic coordinate $BS_1$, $S_1=\operatorname{Im}F^{-1}$. The remainder has the second-row limit with $\sigma=2-p$ and the zero-endpoint solution $Y_{2-p}$, which need not be positive. Its tail determines that of $w$. When $\beta=1$ and $1<p<3$, the ordinary second row applies without this subtraction.

In every nonresonant case $D^2w\in L^{S,\infty}\setminus L^S$ and
\begin{equation}\label{eq:remaindertail}
 \lim_{L\to\infty}L^S\abs{\{|D^2w|>L\}\cap B_\delta(O)}
 =\frac12\int_I|T(\varphi)|^S\dd\varphi>0\,,
\end{equation}
where $T=B\Hmat_\beta[\psi]$ in the first row and $T=Q\Hmat_\sigma[Y_\sigma]$ in the second row and in the exceptional flat sublinear case. In resonance $S=2/(2-\beta)$, and the limit holds instead with $L^S/(\log L)^S$ and constant $\frac12(2-\beta)^{-S}\int_I|D\Hmat_\beta[\psi]|^S$. In this case $D^2w\notin L^{S,\infty}$.
\end{theorem}

\subsubsection{The critical second logarithm and Hessian endpoint}
For the critical alternative in Theorem~\ref{thm:sourcefirst} take $p=p_c$ and $\mu=\beta/2$. Let $P$ be the $L^2(I)$ projection onto $\Span\{\psi\}$ and $Q_0=\operatorname{Id}-P$. On the range of $Q_0$ the operator $H_Q=-\partial_\varphi^2-\beta^2$ with zero endpoints is positive and invertible. Define
\begin{align}
 \Gamma_*&=N^{-1}\int_I\psi^{p+1},&
 A^{p-1}\Gamma_*&=2\beta\mu,& C_*&=a^{-\mu}A,\label{eq:critconstants1}\\
 \chi&=A^pH_Q^{-1}Q_0(\psi^p),&
 D_*&=N^{-1}\int_I\psi^p\chi,& Z_*&=\chi/A,\label{eq:critconstants2}\\
 R_*&=A\mu(\mu+1)+pA^{p-1}D_*,&
 \kappa_*&=-R_*/(2\beta A).\label{eq:critconstants3}
\end{align}
These constants depend only on the fixed equation and geometry; Section~\ref{sec:supersource} proves $D_*>0$ and $\kappa_*<0$.

\begin{theorem}\label{thm:critsecond}
For every critical singular natural solution in Theorem~\ref{thm:sourcefirst} there is a real number $c_u$ such that, with $t=\log(1/r)$,
\begin{equation}\label{eq:critsecond}
 u(F(re^{i\varphi}))=C_*r^{-\beta}t^{-\mu}
 \left[\psi+\frac{(\kappa_*\log t+c_u)\psi+Z_*}{t}
       +O_\psi^1\!\left(\frac{(\log t)^2}{t^2}\right)\right]\,.
\end{equation}
Here $O_\psi^1$ has the meaning fixed in the standing notation. The coefficient $c_u$ depends on the solution and on the fixed coordinate normalization. The coefficient $\kappa_*$ and the angular correction $Z_*$ depend only on the angle and the critical exponent.
\end{theorem}

Put $H_0=\Hmat_{-\beta}[\psi]$, so $|H_0|=\sqrt2\beta(\beta+1)$ and $\tr H_0=0$.

\begin{theorem}\label{thm:critH}
Let $R_\varphi$ be the matrix with columns $e_r(\varphi)$ and $e_\varphi(\varphi)$ in the fixed Cartesian coordinates of the corner chart. For a solution in Theorem~\ref{thm:critsecond}, write $\widehat H_u(r,\varphi)=R_\varphi^{\mathsf T}D_x^2u(F(re^{i\varphi}))R_\varphi$. Then
\begin{equation}\label{eq:critH}
 \widehat H_u(r,\varphi)=C_*r^{-\beta-2}t^{-\mu}
 \left[H_0+\frac{(\kappa_*\log t+c_u)H_0+\Hmat_{-\beta}[Z_*]+\mu\Dmat_{-\beta}[\psi]}t
       +\mathcal R\right]\,.
\end{equation}
Here $\mathcal R=\mathcal R(\rho,\varphi)$ is defined by the identity at a sector radius $\rho$. For every finite $1<\ell<\infty$, with $t_r=\log(1/r)$, its annular norm satisfies
\begin{equation}\label{eq:critHerror}
 \norm{\mathcal R(r|z|,\arg z)}_{L^\ell(\cA)}\le C_\ell t_r^{-2}(\log t_r)^2\,.
\end{equation}
In~\eqref{eq:critH} evaluated at $\rho=r|z|$, the local time is $\log(1/\rho)=t_r-\log|z|$. The estimate is on a complete annulus; it does not assert a pointwise boundary Hessian bound.

Put $S_0=2/(\beta+2)$. Then
\begin{equation}\label{eq:criticaltail}
 \lim_{L\to\infty}L^{S_0}(\log L)^{\mu S_0}
 \abs{\{|D^2u|>L\}\cap B_\delta(O)}
 =\frac\theta2\big(C_*\sqrt2\beta(\beta+1)\big)^{S_0}(\beta+2)^{\mu S_0}\,.
\end{equation}
For every $0<q_0<\infty$, the classical Hessian belongs to $L^{S_0,q_0}$ exactly when $q_0>1/\mu=2/\beta$; it also belongs to $L^{S_0,\infty}$. The Lorentz spaces use the rearrangement convention above, including indices below one. The Hessian belongs to $L^s$ exactly when $0<s<S_0$. For two critical singular solutions on the same domain and with the same coordinate normalization,
\begin{equation}\label{eq:criticaldifference}
 D^2(u-v)\in L^{S_0}(\Om\cap B_\delta(O))\quad\Longleftrightarrow\quad c_u=c_v\,.
\end{equation}
These are statements about integrability of the classical Hessian; $S_0<1$.
\end{theorem}

\begin{theorem}\label{thm:criticalglobal}
Let $\Om$ be a domain as in Section~\ref{sec:linear}, let $a>0$, and let $p=p_c=1+2/\beta$. There are $C_0>0$ and a family $(u_C)_{C\ge C_0}$ of natural solutions of $-\Delta u=au^p$ with zero Martin atom. Every member realizes the singular alternative~\eqref{eq:criticalfirst}, is stable near $O$, and has finite Morse index. In the fixed coordinates of Theorem~\ref{thm:critsecond}, their asymptotic constants satisfy
\begin{equation}\label{eq:criticalglobalphase}
 c_{u_C}-c_{u_D}\asymp D-C\qquad(C_0\le C<D),
\end{equation}
with constants independent of $C,D$. In particular the family is injective, and for distinct $C,D$,
\[
 D^2(u_C-u_D)\notin L^{2/(\beta+2)}(\Om\cap B_\delta(O))
\]
for all sufficiently small $\delta>0$.
\end{theorem}

\subsubsection{Joint corner and side scales for negative-power sources}
For~\eqref{eq:negsource} put $\alpha_q=2/(q+1)$. The side scales used in the joint estimates are, for $0<y<y_0$,
\begin{align}
 S_\eps(y)&=
 \begin{cases}
 y,&0<q<1,\\
 y\sqrt{1+\eps\log(e/y)},&q=1,\\
 y+\eps^{1/(q+1)}y^{\alpha_q},&q>1,
 \end{cases}\label{eq:Ssource}\\
 T_\eps(y)&=y^{-2}\int_0^y t\eps S_\eps(t)^{-q}\dd t\,.
 \label{eq:Tsource}
\end{align}
Let $c\ge0$ be the Martin coefficient of the solution. Its corner amplitude $M(r)$ is defined by
\begin{equation}\label{eq:sourceM}
\begin{array}{c|c|c}
 \text{solution}&M(r)&\text{interior annular profile}\\ \hline
 c>0&c\kappa r^{-\beta}&|z|^{-\beta}\psi\\
 c=0,\ \beta<\alpha_q&B r^\beta&|z|^\beta\psi\\
 c=0,\ \beta=\alpha_q&r^\beta b(\log(1/r))&|z|^\beta\psi\\
 c=0,\ \beta>\alpha_q&r^{\alpha_q}&|z|^{\alpha_q}\Phi
\end{array}
\end{equation}
Here $B>0$ is the regular harmonic coefficient, and $\Phi$ is the unique positive solution of
\begin{equation}\label{eq:negativeangular}
 -\Phi''-\alpha_q^2\Phi=a\Phi^{-q},\qquad \Phi(0)=\Phi(\theta)=0\,.
\end{equation}
In the resonant row, $b$ is the positive first sine coefficient after dividing by $r^\beta$. The existence of these profiles and the meaning of the interior annular limits are established in Section~\ref{sec:negsource}.

\begin{theorem}\label{thm:jointsource}
For every natural solution of~\eqref{eq:negsource}, let $M$ be as in~\eqref{eq:sourceM} and put $\eps(r)=ar^2/M(r)^{q+1}$. Near a side, let $y=d(F(re^{i\varphi}))/r$, truncated above by a fixed small constant. Then
\begin{align}
 u(F(re^{i\varphi}))&\asymp M(r)S_{\eps(r)}(y),\label{eq:jointvalue}\\
 |\nabla_xu(F(re^{i\varphi}))|&\le C\frac{M(r)}r\frac{S_{\eps(r)}(y)}y,\label{eq:jointgrad}\\
 |D_x^2u(F(re^{i\varphi}))|&\le C\frac{M(r)}{r^2}\big(1+T_{\eps(r)}(y)\big).\label{eq:jointhess}
\end{align}
In the middle angular region $y$ may be taken constant. The constants may depend on the solution and the fixed domain, but not on $r$ or $y$. For $q<2$,
\begin{equation}\label{eq:Texplicit}
 T_\eps(y)\asymp
 \begin{cases}
 \eps y^{-q},&q<1,\\
 \eps\big(y\sqrt{1+\eps\log(e/y)}\big)^{-1},&q=1,\\
 \eps\big(y+\eps^{1/(q+1)}y^{\alpha_q}\big)^{-q},&1<q<2.
 \end{cases}
\end{equation}
For $q\ge2$, put $\delta_\eps=\eps^{1/(q-1)}$. Then
\begin{equation}\label{eq:Tlarge}
 T_\eps(y)\asymp
 \begin{cases}
 \eps^{1/(q+1)}y^{\alpha_q-2},&y\le\delta_\eps,\\
 \eps y^{-2}(1+\log(y/\delta_\eps)),&y\ge\delta_\eps,\ q=2,\\
 \delta_\eps y^{-2},&y\ge\delta_\eps,\ q>2.
 \end{cases}
\end{equation}
These are upper estimates for the full Hessian; no optimality assertion for the outer region $q\ge2$ is made.
\end{theorem}

\subsection{Related results}
\begin{citationrevision}
Classical treatments of singular Dirichlet problems and ordinary boundary behavior include~\cite{CRT,LM}. Guo, Li, and Zhang~\cite{GLZ} study well-posedness for continuous Dirichlet data on bounded Lipschitz domains and boundary growth near zero-data portions. For the critical cone $\beta=\alpha_q$ and $q>1$, the growth orders and side transitions in~\eqref{eq:negsourceresonance} and~\eqref{eq:jointvalue} agree in the zero-mass case with~\cite[Theorems~1.1--1.2]{WuZhang}. The planar conformal chart puts this case in their setting with bounded positive coefficient. Our estimates track the Martin-mass branches and control the physical gradient and Hessian through the side layer. The quotient comparison for bounded sublinear solutions is due to Brezis--Oswald~\cite[Section~2]{BO}. Lemma~\ref{lem:Bernstein} proves the interior gradient estimate for negative-power absorption by a weighted logarithmic Bernstein argument related to Ma--Wei~\cite{MW}. See also~\cite{AT} for related geometric methods. The boundary estimates yield the finite trace and the mixed-boundary Hessian endpoint for the natural class. The elementary one-potential construction is an application of Grigor'yan--Verbitsky~\cite[Theorem~3.5(ii)]{GV}. Our two-potential argument improves that elementary constant and gives the stated quantitative contraction and conditional uniqueness.
\end{citationrevision}

Linear asymptotics on angular domains originate in the work of Kondrat'ev~\cite{Kondratiev}; see also Dauge~\cite{Dauge} and Costabel--Dauge~\cite{CD}. Wigley~\cite{Wigley85} studies semilinear corner behavior. These results explain why a harmonic corner term can restrict regularity even when the forcing is less singular. Our use of physical annular estimates distinguishes finite-exponent Hessian convergence from pointwise second-derivative convergence at a merely $C^2$ side. In particular a formally twice-differentiated conformal expansion would not suffice.

For superlinear absorption we use the precise cone classification and removability results of Marcus--V\'eron~\cite{MV}, and give the curved-coordinate comparison needed here. For superlinear sources, Bidaut--V\'eron, Ponce, and V\'eron~\cite{BPV} classify isolated singularities at smooth boundary points. Del Pino, Musso, and Pacard~\cite{DMP} construct continua of singular solutions on smooth domains at and slightly above the boundary-critical exponent. Bidaut--V\'eron and Vivier~\cite[Theorem~1.2]{BVV} prove attainment of the upper mass threshold for subcritical boundary measures on smooth domains. Conical and wedge-like constructions are considered in~\cite{HMR,Gkikas}. The latter results have specific exponent and geometry hypotheses. We keep the local classification of an existing solution separate from global realization of a singular branch. The second logarithmic refinement is derived from the uniform first profile by a quantitative inverse on higher angular modes, a scalar resonant equation, and an integrable remainder. No assumed second-order decay is used in that argument.

\subsection{Reaction, diffusion, and singular amplitudes}\label{subsec:context}
A stationary reaction--diffusion model gives the scalar equation in~\eqref{eq:original}. For $D,\rho>0$ and positive $U$, the evolution equation $\partial_tU=D\Delta U-\rho U^k$ has equilibria satisfying $\Delta u=(\rho/D)u^k$; changing the sign of the reaction gives the source sign. The exponent determines whether the reaction vanishes or diverges as $u$ approaches zero. Prescribed concentrations or temperatures lead to Dirichlet data. Concentrating these data near a corner gives the Martin boundary loading constructed in Section~\ref{subsec:corner-mass}. Its normalization uses harmonic measure and is distinct from prescribing a total normal flux.

There is also a direct probabilistic connection. Harmonic measure describes the exit location of killed Brownian motion, and Martin kernels enter conditional Brownian motion~\cite{Doob}. For the quadratic absorption equation $\Delta u=u^2$ on smooth domains, Le Gall~\cite{LeGall} uses the Brownian snake to represent nonnegative solutions and establishes a correspondence between harmonically bounded solutions and finite boundary measures that do not charge the relevant polar sets. This is a probabilistic realization of a nonlinear equation occurring in our parameter family. It does not assign the same interpretation to every real power. In Lipschitz domains, the measure-trace theory of Marcus--V\'eron~\cite{MV} supplies the corresponding analytic setting for admissible absorption traces, including atoms at conical boundary points.

Singular amplitudes also occur in elasticity. Williams~\cite{Williams} derives angular stress fields at plate corners and notches whose exponents depend on the angle and the mechanical boundary conditions. The coefficient of such a mode records its intensity, as the coefficient of a Martin kernel records a scalar boundary singularity here. This is a comparison of asymptotic structures; an elastic stress coefficient is a different quantity from the Martin mass, and its theory is not used in our uniqueness proofs.
\subsection{Organization and conventions}
Section~\ref{sec:linear} fixes the geometry, potentials, and physical annular estimates. Section~\ref{sec:elementary} treats the harmonic, constant, and linear equations. Sections~\ref{sec:negabs} and~\ref{sec:negsource} treat negative powers, including the minimum mass and the maximal branch for absorption. Section~\ref{sec:sublinear} proves the sublinear classifications and the finite-trace second-order theorem. Sections~\ref{sec:superabs} and~\ref{sec:supersource} treat superlinear classification, source mass thresholds, and the critical singular family. Section~\ref{sec:range} proves Theorem~\ref{thm:existencetable}.

Constants denoted by $c,C$ may change from one occurrence to the next. We write $f\asymp g$ for two positive functions whose quotient stays between positive constants, and $f\sim g$ for relative convergence to one in the specified region. An unqualified boundary derivative is an inward derivative. A finite Martin coefficient means the coefficient in the Green--harmonic representation; its identification with an all-path ratio is stated only when proved. For a measurable $f$ on $D$, write $f^*$ for the decreasing rearrangement of $|f|$. For $S>0$ and $0<\ell<\infty$, membership in $L^{S,\ell}(D)$ means that $\int_0^{|D|}(s^{1/S}f^*(s))^\ell\dd s/s$ is finite. At $\ell=\infty$ use the supremum of $s^{1/S}f^*(s)$. We use these spaces also for measurable classical derivatives when an index is below one, and reserve $W^{2,s}$ for $s\ge1$.

\section{Geometry and linear estimates}\label{sec:linear}
\subsection{The domain and the solution class}
Throughout, $\Om\subset\R^2$ is a bounded connected Lipschitz domain with finitely many boundary components. There is a distinguished point $O\in\partial\Om$. The boundary is $C^2$ away from $O$, and its two incident arcs are $C^2$ up to $O$, with opening $0<\theta<2\pi$. We use $\beta,I,\psi,N,d$ from the standing notation. The case $\theta=\pi$ is included when the two arcs meet with a common tangent. No differentiability of the combined boundary beyond the stated arc assumptions is implicit in this convention.

We call a boundary point or a compact boundary arc \emph{ordinary} if it stays at positive distance from the vertex $O$.

\begin{definition}\label{def:natural}
A natural solution of $\Delta u=\lambda u^k$ is a function
\begin{equation}\label{eq:natural}
 u\in C^2(\Om)\cap C(\overline\Om\setminus\{O\}),\qquad
 u>0\ \text{in }\Om,\qquad
 u=0\ \text{on }\partial\Om\setminus\{O\},
\end{equation}
satisfying the equation pointwise. No boundedness, energy condition, stability, or trace at $O$ is part of this definition.
\end{definition}

We use the Green kernel $G$, the potential operator $\G$, and the normalized Martin kernel $K$ from the standing notation. \citationadded{The normalized kernel and its uniqueness are the classical results recalled in~\cite[Section~2.1, A.1--A.2, p.~923]{MV}. Consequently positive harmonic functions which vanish continuously on $\partial\Om\setminus\{O\}$ are precisely $cK$, $c>0$.}

\begin{lemma}\label{lem:conformal}
After translating $O$ to the origin and rotating the coordinates, there is a conformal map $F$ from a truncated sector into a neighborhood of $O$ in $\Om$, taking its straight sides to the two incident arcs, such that $F(0)=0$ and $F'(0)=1$. We use these coordinates for all corner formulas. For some $\eta_0>0$,
\begin{align}
 F(z)&=z+O(|z|^{1+\eta_0}),& F'(z)&=1+O(|z|^{\eta_0}),\label{eq:conformal}\\
 d(F(re^{i\varphi}))&\asymp r\psi(\varphi),&
 K(F(re^{i\varphi}))&=\kappa r^{-\beta}\psi(\varphi)+O_\psi^1(r^\beta),\label{eq:Martin}
\end{align}
where $\kappa>0$. Here and below $O_\psi^1(h)$ means that the function divided by $\psi$, and its derivatives with respect to $\log r$ and $\varphi$ without division by $\psi$, are $O(h)$. A regular positive harmonic function vanishing on the sides has leading term $B r^\beta\psi$, $B>0$.
\end{lemma}
\begin{proof}
Apply the power map of exponent $\beta$ to the curved sector. The two resulting boundary arcs have a common tangent and $C^{1,\eta}$ regularity for some $\eta>0$. The boundary differentiability theorem for conformal maps gives a nonzero derivative at the origin; composing with the inverse power map gives~\eqref{eq:conformal}. \citationadded{For this boundary differentiability input see~\cite[Chapter~3]{Pommerenke}. The higher corner expansions in~\cite[Section~2]{Wigley65} assume $C^{2,\eta}$ or stronger arc regularity; only the first-order power-map reduction just given is used under the present $C^2$ assumption.} Boundary distance is preserved up to constants because $F$ and $F^{-1}$ are locally Lipschitz. In the straight sector the zero-side harmonic expansion is a sine series in $r^{n\beta}$ and $r^{-n\beta}$. Positivity excludes all negative modes except the first. The local singular harmonic function is therefore $\kappa r^{-\beta}\psi$ plus a regular zero-side harmonic function. Its sine expansion and boundary estimates give~\eqref{eq:Martin}. The same expansion and positivity give the last assertion.
\end{proof}

We will not differentiate $F$ twice at the vertex. Under the present assumptions $F''$ need not be bounded there. This matters for physical Hessians.

\subsection{The finite corner mass}\label{subsec:corner-mass}
With the normalization $K(x_*)=1$, a natural solution is said to have finite corner mass $c$ when it admits the Green--harmonic representation
\begin{equation*}
 u+\G[\lambda u^k]=cK,\qquad 0\le c<\infty\,,
\end{equation*}
with a finite Green potential. Equivalently, the harmonic component has Martin boundary measure $\mu=c\delta_O$. We use \emph{trace~$c$} as shorthand for this prescribed measure and call $c$ its Martin mass. Every finite positive measure supported at the single point $O$ has this form. Prescribing $c$ therefore specifies the one remaining scalar boundary datum in the finite-trace class. A result for every admissible $c$ concerns the whole family of such traces; specifying the mass alone imposes no symmetry, stability, or global lower bound on the solution. Changing the normalization point to $x^\dagger$ replaces $K$ by $K/K(x^\dagger)$ and $c$ by $cK(x^\dagger)$, so equality of the masses of two solutions is independent of a common choice of normalization.

The boundary datum can be obtained by concentrating ordinary Dirichlet data. Let $\omega^x$ denote harmonic measure in~$\Om$ with pole~$x$, and put $E_\varepsilon=\partial\Om\cap B_\varepsilon(O)$. For sufficiently small $\varepsilon>0$, its harmonic measure is positive. Define
\begin{equation*}
 g_\varepsilon=\frac{c\one_{E_\varepsilon}}{\omega^{x_*}(E_\varepsilon)},\qquad
 \mu_\varepsilon=g_\varepsilon\omega^{x_*}\weak c\delta_O\,.
\end{equation*}
The harmonic extension $h_\varepsilon$ of $g_\varepsilon$ satisfies
\begin{equation*}
 h_\varepsilon(x)=c\frac{\omega^x(E_\varepsilon)}{\omega^{x_*}(E_\varepsilon)}
 =\frac{c}{\omega^{x_*}(E_\varepsilon)}
   \int_{E_\varepsilon}K(x,y)\dd\omega^{x_*}(y)
 \longrightarrow cK(x,O)\,.
\end{equation*}
Here $K(x,y)$ is the Martin kernel normalized at $x_*$ for a boundary point $y$, so $K(x,O)=K(x)$. The identity $K(x,y)=\dd\omega^x/\dd\omega^{x_*}(y)$ and its continuity in the boundary variable are recalled in~\cite[Section~2.1, A.1--A.2]{MV}. The displayed limit follows by averaging this continuous kernel over sets shrinking to $O$; it is locally uniform in~$\Om$. The measure convergence follows by testing against continuous functions on~$\partial\Om$. This construction concerns the boundary datum for the harmonic problem. Existence of a nonlinear solution with that datum is a separate assertion in the equation-specific sections.

The word mass refers to this boundary measure. It is neither the volume integral of $u$ nor an ordinary total normal flux. Nor is fixing a mass a harmless rescaling for a fixed nonlinear equation: if $c>0$ and $u=cv$, then $\Delta v=\lambda c^{k-1}v^k$. An identification of $c$ with $\lim u/K$ requires the relative smallness of the Green term and is made only in the statements proving that property. For negative-power absorption this is Theorem~\ref{thm:finiteabs}; for negative-power sources tangential side layers require the distinction stated in Proposition~\ref{prop:negsourceclass}.

Prescribed mass is not by itself a universal uniqueness criterion. Theorem~\ref{thm:twoGreen} requires a global lower bound for negative-power absorption, Proposition~\ref{prop:supersourceexist} gives distinct superlinear source solutions with the same positive mass, and Proposition~\ref{prop:elementary} gives an eigenfunction ray of zero mass at resonance. Strong absorption singularities lie outside the finite-mass representation, whereas critical and supercritical source singularities can have zero Martin mass, as in Theorem~\ref{thm:sourcefirst}. These distinctions separate the missing boundary datum from multiplicity within a prescribed-trace class.

\subsection{Linear potentials for a negative power}
\begin{lemma}\label{lem:negativepotential}
If $0<p<1$, then $H=\G[K^{-p}]$ is finite, belongs to $H^1_0(\Om)$, and satisfies
\begin{equation}\label{eq:HM}
 0<M:=\sup_\Om\frac HK<\infty,\qquad
 \lim_{x\to O}\frac{H(x)}{K(x)}=0\,.
\end{equation}
The last limit is uniform along all interior paths. Moreover $H\le Cd$ on compact ordinary boundary pieces.
\end{lemma}
\begin{proof}
First, $K^{-p}$ defines a bounded functional on $H^1_0(\Om)$ by Hardy's inequality: $dK^{-p}\in L^2$, using~\eqref{eq:Martin} at the corner and $K\asymp d$ on ordinary pieces. Its variational Dirichlet solution equals the Green potential by monotone approximation, so $H$ is finite and in $H^1_0$. On an ordinary boundary piece $K^{-p}\le Cd^{-p}$. The local Green gradient estimate is integrable because $p<1$, and gives $H\le Cd$ after separating a remote source. Near $O$, choose
\begin{equation*}
 0<s<\min\{\beta,2+\beta p\}\,.
\end{equation*}
The positive solution of $-f''-s^2f=\psi^{-p}$, $f(0)=f(\theta)=0$, satisfies $f\asymp\psi$. A large multiple of $r^sf$ is a supersolution for $H\circ F$, since $r^{s-2}\ge Cr^{\beta p}$ at small radii. For a precise comparison take bounded compactly supported sources $f_n\uparrow K^{-p}$. Their potentials are regular harmonic near $O$, hence $O(r^\beta\psi)$, and are dominated by the barrier on sufficiently small inner arcs because $s<\beta$. The constant on the outer arc can be chosen independent of $n$, since these potentials are bounded above by $H$. Compare in the punctured sector, let the inner radius tend to zero, and then let $n\to\infty$. Thus $H\circ F\le Cr^s\psi$ and $H/K\le Cr^{s+\beta}\to0$. Interior compactness and the ordinary boundary estimates prove finiteness of $M$.
\end{proof}

\begin{lemma}\label{lem:positivepotential}
If $p>0$ and $\beta(p-1)<2$, then $H_p=\G[K^p]$ is finite, $M_p=\sup_\Om H_p/K<\infty$, and $H_p/K\to0$ uniformly at $O$.
\end{lemma}
\begin{proof}
Near $O$, $G(x_*,\cdot)$ is a regular positive harmonic function that vanishes on the sides. By Lemma~\ref{lem:conformal}, it is comparable in the chart to $r^\beta\psi$. The moment $\G[K^p](x_*)$ thus has radial factor $r^{\beta+1-\beta p}$ and angular factor $\psi^{p+1}$. This moment is finite precisely when $\beta(p-1)<2$. The same comparison at any interior observation point gives finiteness there. On compact ordinary boundary pieces, a local bounded source and boundary comparison for the remote part give $H_p\le Cd$.

Choose $-\beta<s<\min\{\beta,2-\beta p\}$ and solve $-f''-s^2f=\psi^p$ with zero endpoints. The positive inverse gives $f\asymp\psi$. The bound on ordinary boundary pieces controls the outer-arc values of all compact source truncations. A large multiple of $r^sf$ dominates these values and is a supersolution near $O$. Comparing the truncated potentials as in Lemma~\ref{lem:negativepotential} gives $H_p\circ F\le Cr^s\psi$, and therefore $H_p/K=O(r^{s+\beta})$ at the vertex. This estimate, the ordinary boundary bound, and interior compactness give $M_p<\infty$.
\end{proof}

\subsection{Estimates on annuli with curved sides}
Put $\cA=\{z:1<|z|<2,\ 0<\arg z<\theta\}$ and let $\cA'$ denote a fixed enlargement with radii $1/2$ and $4$. Estimates on these sets are local at the artificial circular boundaries.

\begin{lemma}\label{lem:physical}
Fix $1<\ell<\infty$. For sufficiently small $r>0$, let $v$ be continuous up to the boundary of the enlarged physical annulus, with $v=0$ on the true sides and distributional $\Delta v\in L^\ell$ there. There is a constant $C_\ell$, uniform in $r$, such that
\begin{equation}\label{eq:annestimate}
 \norm{D_x^2v}_{L^\ell(\Om\cap\{r<|x-O|<2r\})}
 \le C_\ell\left(r^{-2}\norm v_{L^\ell(\Om\cap\{r/2<|x-O|<4r\})}
       +\norm{\Delta v}_{L^\ell(\Om\cap\{r/2<|x-O|<4r\})}\right)\,.
\end{equation}
Let $v_r$ satisfy these assumptions on the corresponding enlarged annuli and let $M_r>0$. If
\begin{equation*}
 \frac{v_r(F(rz))}{M_r}\longrightarrow V(z)\quad\hbox{uniformly on }\cA',\qquad
 \frac{r^2\Delta_xv_r(F(rz))}{M_r}\longrightarrow\Delta_zV(z)\quad\text{in }L^\ell(\cA'),
\end{equation*}
then $r^2D_x^2v_r(F(r\cdot))/M_r\to D_z^2V$ in $L^\ell(\cA)$. For a bounded family of explicit $C^{2,\alpha}$ profiles, positive-power errors in the functions and Laplacians give a positive-power error in this conclusion.
\end{lemma}
\begin{proof}
\begin{citationrevision}
On a rescaled enlarged annulus, the true sides have uniformly bounded $C^{1,1}$ graph norms: their difference from the tangent rays is $O(r)$ in $C^2$. The interior and local boundary estimates of~\cite[Theorems~9.11 and~9.13, pp.~235 and~239]{GT}, followed by a finite covering and rescaling, give~\eqref{eq:annestimate}. They are applied only along the true sides; the artificial circular ends are avoided by the enlargement. Where boundary $W^{2,\ell}$ membership is not initially known, take a smooth-sided auxiliary patch and approximate both the Laplacian in $L^\ell$ and the continuous artificial-boundary values uniformly, keeping zero values on the true-side portion. Solve these smooth Dirichlet problems with the approximated data. Convergence of their harmonic liftings and Green potentials identifies the limit with $v$; here the Green potentials converge uniformly since the dimension is two and $\ell>1$. The cited local estimates then pass to the limit. Thus no zero condition on the artificial boundary, or global energy assumption on a corner-singular solution, is introduced.
\end{citationrevision}

For the convergence assertion, use slightly nested enlarged annuli and choose graph maps $\Xi_r$ from the straight annuli to the rescaled physical ones, with $\norm{\Xi_r-\mathrm{Id}}_{C^2}=O(r)$. The normalized pullbacks
$U_r(z)=v_r(O+r\Xi_r(z))/M_r$ converge uniformly to $V$. The pulled-back equations are
$\mathcal L_rU_r=f_r$, with $\mathcal L_r=a_r:D^2+b_r\cdot\nabla$,
$a_r-I_2$ and $b_r$ tending to zero uniformly, and $f_r\to\Delta V$ in $L^\ell$. Here $I_2$ is the two-dimensional identity matrix.
The two pullback parametrizations approach one another in $C^1$ and have uniformly controlled Jacobians. This gives the asserted convergence of the pulled-back sources by continuity of composition in $L^\ell$, first for continuous functions and then by density.

The uniform annular estimate and weak compactness first give $V\in W^{2,\ell}$ on an intermediate annulus. The uniform local estimate applied to $U_r-V$, which has zero true side trace, then gives strong $W^{2,\ell}$ convergence on the smaller annulus, because
\[
  \mathcal L_r(U_r-V)=f_r-\Delta V+(I_2-a_r):D^2V-b_r\cdot\nabla V\longrightarrow0
  \quad\text{in }L^\ell.
\]
The chain rule for $\Xi_r$ transfers this convergence to the physical Hessians, since $D\Xi_r\to I$ and $D^2\Xi_r\to0$. Composition with the near-identity change to the conformal parametrization transfers the $L^\ell$ limit without using $F''$. For a bounded family of $C^{2,\alpha}$ profiles, the composition error in the limiting Hessians is bounded by a positive power of $r$. The same estimates then preserve a positive-power error in the conclusion.
\end{proof}

The matrices in~\eqref{eq:Hmat} are written in the orthonormal polar frame. Then $D^2(r^\lambda Y)=r^{\lambda-2}\Hmat_\lambda[Y]$. More generally, for $t=\log(1/r)$,
\begin{equation}\label{eq:matrixidentity}
 D^2\big(r^\lambda W(t,\varphi)\big)
 =r^{\lambda-2}\left(\Hmat_\lambda[W]-\Dmat_\lambda[W_t]+\diag(W_{tt},0)\right)\,.
\end{equation}

\begin{lemma}\label{lem:tail}
Suppose a physical Hessian has, on complete shrinking annuli, strong $L^\ell$ asymptotics $r^{-b}T(\varphi)$ for some $b>0$, with uniform normalized annular bounds for an exponent $\ell>2/b$. If $T\not\equiv0$ and $S=2/b$, then
\begin{equation}\label{eq:tailhom}
 \lim_{L\to\infty}L^S\abs{\{x\in\Om\cap B_\delta(O):|D^2v(x)|>L\}}
 =\frac12\int_I|T(\varphi)|^S\dd\varphi\,.
\end{equation}
Consequently $D^2v\in L^{S,\infty}$ and $D^2v\notin L^S$. If the leading factor is $r^{-b}\log(1/r)$, the left side must instead be multiplied by $(\log L)^{-S}$ and the right side by $b^{-S}$; in this case $D^2v\notin L^{S,\infty}$.
\end{lemma}
\begin{proof}
For the exact model integrate $r\dd r\dd\varphi$ over $r<(\abs T/L)^{1/b}$. For a remainder, set $r=L^{-1/b}s$, restrict first to $\eps<s<\eps^{-1}$, and use strong annular convergence. The limiting level set $|z|^{-b}|T(\arg z)|=1$ has area zero, since it has at most one radius at each angle. The inner set has area $O(\eps^2)$ in the rescaled variables. The outer annuli have a summable Chebyshev bound because $b\ell>2$. Letting $\eps\downarrow0$ proves~\eqref{eq:tailhom}. The Jacobian of $F$ tends to one. The logarithmic version uses $L=R^{-b}\log(1/R)$ and the same argument. It follows that the Hessian belongs to $L^s$ precisely for $0<s<S$ in the homogeneous case. When $S<1$, these assertions concern the classical Hessian as a measurable function; they are not statements about the Banach Sobolev space $W^{2,S}$.
\end{proof}

\section{Harmonic, constant, and linear equations}\label{sec:elementary}
\subsection{The sector model}\label{subsec:sector-model}
Write $S_\theta=\{re^{i\varphi}:0<r<1,\ 0<\varphi<\theta\}$ and $\beta=\pi/\theta$. The positive harmonic functions on this sector with continuous zero boundary values away from $O=0$ are
\begin{equation*}
 u_A(r,\varphi)=A(r^{-\beta}-r^\beta)\sin(\beta\varphi),\qquad A>0\,.
\end{equation*}
The radial factor in the remark following Problem~50~\cite[Section~6.4, p.~27]{Mazya} is printed as $r^{-\beta}-1$. It must be replaced by $r^{-\beta}-r^\beta$: the identity
\begin{equation*}
 \Delta\bigl(r^s\sin(\beta\varphi)\bigr)
 =(s^2-\beta^2)r^{s-2}\sin(\beta\varphi)
\end{equation*}
shows that the printed expression has Laplacian $A\beta^2r^{-2}\sin(\beta\varphi)$, whereas the corrected function is harmonic and vanishes also on $r=1$. Completeness of the corrected family follows from the single-point Martin representation. The amplitude $A$ is proportional to, and need not equal, the mass defined using $K(x_*)=1$.

The points $P_1=1$ and $P_2=e^{i\theta}$ are also geometric corners, both of angle $\pi/2$. Thus $S_\theta$ is a truncated corner model rather than literally a domain smooth outside $O$. The boundary condition has a different role at these points: continuity and zero values are required there, whereas $O$ is excluded. Consequently the two outer corners supply no additional Martin masses. They cannot, however, be treated as smooth boundary points in a global regularity estimate.

One may remove these outer corners while preserving the entire sector near $O$. Fix $0<\varepsilon<\min\{1/4,\theta/4\}$, choose an even nonnegative smooth function $\rho$ with integral one, positive on $(-1,1)$ and zero outside that interval, and set
\begin{align*}
 \rho_\varepsilon(s)&=\varepsilon^{-1}\rho(s/\varepsilon),&
 h_\varepsilon(t)&=\int_\R |t-s|\rho_\varepsilon(s)\dd s,\\
 X_\varepsilon(t)&=\frac{h_\varepsilon(t)+t}{2},&
 Y_\varepsilon(t)&=\frac{h_\varepsilon(t)-t}{2}\,.
\end{align*}
Convolution and convexity give $h_\varepsilon\in C^\infty(\R)$, $h_\varepsilon\ge|t|$, $h_\varepsilon(t)=|t|$ for $|t|\ge\varepsilon$, and $|h_\varepsilon'|\le1$, with strict inequality when $|t|<\varepsilon$. On $[-\varepsilon,\varepsilon]$ one also has $h_\varepsilon\le\varepsilon$. Define the two regular arcs
\begin{align*}
 \gamma_1(t)&=(1-X_\varepsilon(t))e^{iY_\varepsilon(t)},\\
 \gamma_2(t)&=(1-X_\varepsilon(t))e^{i(\theta-Y_\varepsilon(t))},\qquad |t|\le\varepsilon\,.
\end{align*}
Their coordinate curves have derivative $X_\varepsilon'-Y_\varepsilon'=1$, hence are regular; the polar maps used above have nonzero Jacobian since $1-X_\varepsilon>0$. Extending these formulas in $t$, the part $t\ge\varepsilon$ is exactly a radial side and the part $t\le-\varepsilon$ is exactly the unit circle. Thus all derivatives match at the junctions.

Let $\Gamma_\varepsilon$ be the union of the radial segments from $0$ to $(1-\varepsilon)$ and $(1-\varepsilon)e^{i\theta}$, the arcs $\gamma_1,\gamma_2$, and the unit-circle arc $\varepsilon\le\varphi\le\theta-\varepsilon$. The strict monotonicity of $Y_\varepsilon$ on the interiors of the replacement arcs orders these pieces by polar angle and excludes self-intersections. They form a Jordan curve; let $\Om_\varepsilon$ be its bounded interior. The construction gives
\begin{align*}
 &\Om_\varepsilon\subset S_\theta,\qquad
 S_\theta\setminus\Om_\varepsilon\subset B_{2\varepsilon}(P_1)\cup B_{2\varepsilon}(P_2),\\
 &\Om_\varepsilon\cap B_{1-\varepsilon}(0)
   =S_\theta\cap B_{1-\varepsilon}(0),\qquad
 \partial\Om_\varepsilon\setminus\{0\}\text{ is }C^\infty\,.
\end{align*}
Indeed each modified part lies in $1-\varepsilon\le r\le1$ and within angular distance $\varepsilon$ of its endpoint; there $|re^{i\varphi}-1|\le|r-1|+|e^{i\varphi}-1|\le2\varepsilon$, with the rotated estimate at $P_2$. The sides near $O$ and their angle are unchanged. For $\theta\ne\pi$, this is a domain with exactly one geometric corner; for $\theta=\pi$ the marked point is smooth as well. On the new domain the harmonic family is $cK_{\Om_\varepsilon}(\cdot,O)$, with the same local corner exponent $\beta$. Its global kernel and the constants in nonlinear estimates must be computed on the new domain; the explicit sector formula is not transported by this boundary modification.

\subsection{Existence and uniqueness after normalization}
\begin{proposition}\label{prop:elementary}
For $\lambda=0$, all natural solutions are $cK$, $c>0$, irrespective of $k$. For $k=0$ they are exactly
\begin{equation}\label{eq:constant}
 u=cK-\lambda\G[1],\qquad c\ge0,
\end{equation}
subject to strict positivity in $\Om$; when $\lambda>0$ this requires $c>0$ and is satisfied for every sufficiently large $c$.

For $k=1$, let $\Lambda_1$ be the first Dirichlet eigenvalue of $-\Delta$. If $\lambda>-\Lambda_1$, all natural solutions are positive multiples of the Martin kernel at $O$ of $-\Delta+\lambda$. If $\lambda=-\Lambda_1$, they are positive multiples of the first eigenfunction. If $\lambda<-\Lambda_1$, there is no natural solution. Each existing ray is nonunique before normalization. At resonance all its members have zero Laplace Martin coefficient, so prescribing that coefficient does not select a solution.
\end{proposition}
\begin{proof}
\citationadded{The harmonic assertion is the one-pole Martin representation~\cite[Section~2.1, A.1--A.2, p.~923]{MV}.}
For $k=0$, the function $h=u+\lambda\G[1]$ is harmonic, bounded below, and zero at the ordinary boundary. Choose $R_0>\sup_{x\in\Om}|x-O|$. For any $\delta>0$ the harmonic function $h+\delta\log(R_0/|x-O|)$ is nonnegative on the ordinary boundary and on a sufficiently small inner circle. The maximum principle and then $\delta\downarrow0$ give $h\ge0$. Thus $h=cK$ with $c\ge0$, including the identically zero case. Near the corner, the angular barrier used in Lemma~\ref{lem:positivepotential}, with $-Y''-s^2Y=1$ and $0<s<\min\{\beta,2\}$, gives $\G[1]\circ F\le Cr^s\psi$. Therefore $\G[1]/K\to0$ there. On ordinary boundary pieces $\G[1]\le Cd$ and $K\asymp d$, so the quotient is bounded globally. This proves the large-$c$ assertion.

Suppose $\lambda>-\Lambda_1$. The Dirichlet form of $L=-\Delta+\lambda$ is coercive. Its variational inverse preserves nonnegativity by testing with the negative part. For $Lz=1$ with $z\in H^1_0(\Om)$ this gives $z\ge0$. The function $z=\G[1-\lambda z]$ is bounded and continuous with zero boundary values by the planar Sobolev embedding into finite $L^q$ and the Laplace Green estimate for $L^q$ sources, $q>1$, on a bounded regular domain. Thus $\sigma=1+Cz\ge1$ satisfies $L\sigma=\lambda+C>0$ for a sufficiently large $C$. An interior negative minimum of an $L$-harmonic function divided by $\sigma$ contradicts its equation. This proves the maximum principle for continuous boundary data on $\Om$ and its truncations.

For each interior pole, cut off the logarithmic fundamental solution in a compact interior neighborhood and correct its remaining $L^q$ forcing by the variational inverse of $L$. The same Laplace Green estimate makes the correction bounded and continuous with zero boundary values. This constructs the Dirichlet Green kernel $G_\lambda$, which is positive by the maximum principle off the pole and the positive logarithmic divergence there.

We record the local harmonic comparison used for existence and the corner profile. Let $D$ be a bounded Lipschitz planar subdomain and let $v>0$ be a bounded solution of $-\Delta v+Bv=0$, continuous on $\overline D$, where $B$ is bounded. Write $h$ for the harmonic lifting of its boundary values and $\G_D$ for the Laplace Green operator. The Cranston--McConnell bound~\cite[Theorem~A]{AikawaMurata} is $\G_D[h]\le C_0|D|h$, with an absolute constant $C_0$. If $C_0|D|\norm{B}_\infty\le1/4$, the Neumann series for $(I+\G_DB)^{-1}h$ converges in the norm $\sup_D|\cdot|/h$ and differs from $h$ by at most $h/3$. It equals $v$: both solve $v+\G_D[Bv]=h$, and this equation is unique among bounded functions by the same contraction bound with the positive harmonic function $1$. In particular,
\[
 \tfrac23 h\le v\le\tfrac43 h.
\]

To construct the pole, take $y_j\to O$ in $\Om$ and normalize $G_\lambda(\cdot,y_j)$ to be one at $x_*$. Interior Harnack and local elliptic estimates give a positive locally convergent subsequence, whose limit solves $Lu=0$. On a small ordinary boundary patch, use the harmonic comparison on $D=\Om\cap B_R(P)$ with $B=\lambda$, keeping the poles outside $D$. The harmonic liftings are zero on the true boundary and bounded at a fixed interior point by interior Harnack for the normalized kernels. The harmonic boundary Harnack principle on a smaller patch, with a fixed positive zero-side harmonic function, bounds these liftings by $Cd$. The limit therefore vanishes continuously at every ordinary boundary point and is a natural solution.

For any natural solution of $Lu=0$, put $v=u\circ F$ and $B=\lambda|F'|^2$ on the small straight sector $D_R$. Apply the comparison to each annular truncation $D_R\cap\{r>\eps\}$. Harmonic compactness gives a positive harmonic limit $h$ comparable to $v$, which vanishes on the open sides. The sector expansion yields $v=O(r^{-\beta}\psi)$. If $h$ were regular at the vertex, $v$ would extend continuously by zero there and the coercive maximum principle on $\Om$ would force $u=0$. Thus $h$ has a positive first singular coefficient and $v\ge c_0r^{-\beta}\psi$ near the vertex.

To identify the coefficient of $v$ itself, choose $-\beta<s<\min\{\beta,2-\beta\}$. Then $Cr^s\psi$ is a supersolution for the Laplace potential of $|Bv|$ in $D_R$, since $-\Delta(r^s\psi)=(\beta^2-s^2)r^{s-2}\psi$. Compare first with compact source truncations, whose potentials are regular harmonic at the vertex, and pass to the limit. This gives $\abs{\G_{D_R}[Bv]}\le Cr^s\psi=o(r^{-\beta}\psi)$. The harmonic function $v+\G_{D_R}[Bv]$ vanishes on the sides and is $O(r^{-\beta}\psi)$, so its sine expansion has at most the first singular mode. The lower bound on $v$ shows that its coefficient $A$ is positive. Consequently
\[
  \frac{u(F(re^{i\varphi}))}{r^{-\beta}\psi(\varphi)}\longrightarrow A>0
  \quad\text{uniformly for }0<\varphi<\theta .
\]
For any two natural solutions normalized to have the same $A$, this limit orders each below $1+\delta$ times the other on sufficiently small inner caps. The maximum principle on the truncated domain and then $\delta\downarrow0$ prove equality. All subsequential normalized Green limits are therefore the same. This identifies the normalized solution as the Martin kernel at $O$ and proves the ray assertion.

To handle the other values, write $a=-\lambda$ and use the ground-state identity for any positive solution of $-\Delta u=au$:
\begin{equation*}
 \int_\Om(|\nabla\xi|^2-a\xi^2)=\int_\Om u^2|\nabla(\xi/u)|^2\ge0,
 \qquad\xi\in C_c^\infty(\Om)\,.
\end{equation*}
Approximation of the first eigenfunction in $H^1_0$ gives $a\le\Lambda_1$. At equality the right side tends to zero on every interior compact set, so $u$ is proportional to that eigenfunction. This argument does not assume that $u$ itself has finite energy.
\end{proof}

\subsection{Second-order regularity}
Every solution in Proposition~\ref{prop:elementary} with a nonzero first singular coefficient $A$ has the raw Hessian tail~\eqref{eq:rawtail}, with $c\kappa$ replaced by $A$. In the harmonic and constant-source cases the normalized function profile follows from~\eqref{eq:Martin} and $\G[1]/K\to0$. In the linear case the profile was proved above, and its normalized Laplacian is $r^{\beta+2}\lambda u(F(r\cdot))=O(r^2)$ on complete annuli. Lemmas~\ref{lem:physical} and~\ref{lem:tail} give the stated tail. The first eigenfunction and the positive zero-mass constant-source solution are bounded and continuous with zero data at the corner, by the same Laplace $L^q$ Dirichlet estimate. They satisfy Lemma~\ref{lem:boundedforcing}, which gives their exact finite Sobolev range $1\le s<P_\beta$.

\section{Negative-power absorption}\label{sec:negabs}
We use the equation~\eqref{eq:negabs} and the exponents in~\eqref{eq:gammaabs}. The restriction $p<1$ is forced by the ordinary zero boundary condition.

\subsection{Every natural solution has a finite trace}

We first prove the pointwise gradient estimate used below. The
logarithmic Bernstein approach is related to that of
Ma--Wei~\cite[Theorem~1]{MW}.

\begin{lemma}\label{lem:Bernstein}
If $\Delta u=a u^{-p}$ in $B_R$ and $0<u\le M$, then
\begin{equation}\label{eq:Bernstein}
 \sup_{B_{R/2}}\abs{\nabla u^{1/\gamma}}^2
 \le C_p\left(a+\frac{M^{1+p}}{R^2}\right),\qquad
 |\nabla u|\le C_pu^{(1-p)/2}\left(\sqrt a+\frac{M^{(1+p)/2}}R\right)\,.
\end{equation}
\end{lemma}
\begin{proof}
Put $z=\log u$, $S=|\nabla z|^2$, $A=1+p$, $X=au^{-A}$,
and $D=\Delta+2\nabla z\cdot\nabla$. Choose a parameter $b>0$
such that $2-A-bA^2/4>0$ and set
\[
 \Psi=\frac1b\log(X^b+R^{-2b}),\qquad
 \chi =\frac{X^b}{X^b+R^{-2b}}.
\]
Using $\Delta z=X-S$, we obtain
\[
 DS=2|D^2z|^2-2AXS,\qquad
 D\Psi=-A\chi (X+S)+bA^2\chi (1-\chi )S.
\]
Take a smooth cutoff $0\le\eta\le1$ supported in $B_R$, equal to
one on $B_{R/2}$, with $|\nabla\eta|\le C/R$ and
$|\Delta\eta|\le C/R^2$. At a positive maximum of
$\eta^2S e^{-\Psi}$, put $J=(R\eta)^{-1}$. The first and second
derivative tests give
\[
 \nabla\log S=-A\chi \nabla z-2\nabla\log\eta,\qquad
 D(\log S-\Psi)\le C(J^2+J\sqrt S).
\]
Take the first coordinate direction along $\nabla z$. The two diagonal
entries of $D^2z$ are $-A\chi S/2+O(J\sqrt S)$ and
$(-1+A\chi /2)S+O(X+J\sqrt S)$, respectively. The first identity follows
by differentiating $S$ in the first direction, and the second by
taking the trace. Bounding $|D^2z|^2$ below by the sum of the
squares of these entries gives
\[
 D(\log S-\Psi)
 \ge \bigl(2-A\chi -bA^2\chi (1-\chi )\bigr)S-C(X+J\sqrt S+J^2).
\]
Here we used $D\log S=(DS)/S-|\nabla\log S|^2$ and bounded
the cross terms in the squares by $C(X+J\sqrt S)$ after division
by $S$. Since $0\le \chi \le1$, the coefficient of $S$ is bounded
below by $2-A-bA^2/4>0$. Young's inequality at the maximum gives
$S\le C_p(X+J^2)$. Because $e^\Psi\ge\max(X,R^{-2})$ and
$\eta^2J^2=R^{-2}$, the maximum of $\eta^2Se^{-\Psi}$ is
bounded by $C_p$. On $B_{R/2}$ we have $\eta=1$ and
$e^\Psi\le C_p(X+R^{-2})$. Thus throughout $B_{R/2}$,
\[
 |\nabla u|^2\le C_p\left(R^{-2}u^2+au^{1-p}\right).
\]
Multiplying this estimate by
$\gamma^{-2}u^{p-1}$ and using $u\le M$ proves the first estimate
of~\eqref{eq:Bernstein}. The chain rule gives the second.
\end{proof}

\begin{lemma}\label{lem:localabs}
On a $C^2$ boundary piece avoiding $O$, let $D_R=\Om\cap B_R(P)$ and $M=\sup_{D_{4R}}u$. In a smaller piece,
\begin{align}
 u&\le CMd/R,& |\nabla u|&\le C_p\big(M/R+\sqrt a\,M^{(1-p)/2}\big),\label{eq:localabs1}\\
 u|D^2u|+|\nabla u|^2&\le C_p\big(M^2/R^2+aM^{1-p}\big),&
 \int_{D_R}au^{-p}&\le C_p\big(M+\sqrt a\,R M^{(1-p)/2}\big).\label{eq:localabs2}
\end{align}
The zero extension of $u^2$ is locally $C^{1,1}$. The Hessian of $u$ is locally weak $L^1$ up to this boundary piece.
\end{lemma}
\begin{proof}
A harmonic function with value zero on the true boundary and value $M$ on the artificial boundary dominates the subharmonic function $u$. Its boundary Lipschitz estimate gives the first bound. Apply~\eqref{eq:Bernstein} in balls with radius comparable to $d(x)$ to obtain the gradient bound. If $L$ is that bound and $u(x)=u_x$, the ball of radius $u_x/(8L)$ lies in $\Om$ and $u\asymp u_x$ there. Interior Schauder estimates give
\begin{equation*}
 |D^2u(x)|\le C_p\big(L^2/u_x+a u_x^{-p}\big)\,.
\end{equation*}
This proves the weighted Hessian bound. Now $D^2(u^2)=2\nabla u\otimes\nabla u+2uD^2u$ is bounded and $\nabla(u^2)\to0$ on the true boundary. Its zero extension has no boundary measure in its second derivatives.

To estimate the source, test with $\zeta\chi_\eps(d)$, where $\zeta$ is a local cutoff and $\chi_\eps$ vanishes for $d<\eps$ and equals one for $d>2\eps$. The layer has area $O(\eps R)$ and $|\nabla u|\le L$, so $\int au^{-p}\zeta\chi_\eps\le CLR$. Fatou's lemma gives~\eqref{eq:localabs2}. Finally localize $u$ in a larger $C^{1,1}$ auxiliary domain. Its Laplacian is in $L^1$, with the same bound. \citationadded{The Dirichlet weak type $(1,1)$ Hessian estimate~\cite[Theorem~1.20]{DEK} gives the last assertion: the auxiliary domain is $C^{1,1}$ and its coefficient matrix is the identity, so the required Dini mean oscillation is zero.} For $L^1$ data use truncation, convergence in $W^{1,q}$ for $1<q<2$, and interior convergence of the Hessians. Thus no boundary lower bound for $u$ has been assumed.
\end{proof}

\begin{lemma}\label{lem:additive}
For positive solutions in $B_{2R}$,
\begin{align}
 R|\nabla u|&\le C_p\big(u+a^{1/(1+p)}R^\gamma\big)&&\text{in }B_R,\label{eq:addgrad}\\
 \sup_{B_R}u+a^{1/(1+p)}R^\gamma
 &\le C_p\big(\inf_{B_R}u+a^{1/(1+p)}R^\gamma\big).\label{eq:addH}
\end{align}
On a uniform $C^{1,1}$ boundary piece, with a fixed interior point $A(P,R)$ at depth comparable to $R$,
\begin{equation}\label{eq:Carleson}
 \sup_{\Om\cap B_R(P)}u
 \le C\big(u(A(P,R))+a^{1/(1+p)}R^\gamma\big)\,.
\end{equation}
\end{lemma}
\begin{proof}
Set $A_R=a^{1/(1+p)}R^\gamma$. Taking square roots in the pointwise estimate from the proof of
Lemma~\ref{lem:Bernstein} on $B_R$ gives
\[
 R|\nabla u|\le C_p\bigl(u+A_R^{(1+p)/2}u^{(1-p)/2}\bigr)
 \le C_p(u+A_R).
\]
The last inequality is the weighted arithmetic--geometric mean
inequality. Integrate the gradient of $\log(u+A_R)$ along line
segments in $B_R$ to obtain~\eqref{eq:addH}.

For~\eqref{eq:Carleson}, scale to one and put $H=u(A(P,1))+a^{1/(1+p)}$. Interior Harnack chains give $u(x)\le Cd(x)^{-N}H$ on a fixed larger piece. Let $M(s)=\sup_{\Om\cap B_s(P)}u$. Splitting points according to $d(x)\ge\delta(t-s)$, and using a local harmonic upper bound for the other points, gives
\begin{equation*}
 M(s)\le C[\delta(t-s)]^{-N}H+C\delta M(t),\qquad 1\le s<t\le2\,.
\end{equation*}
Choose $C\delta<2^{-N-1}$ and iterate over $s_j=2-2^{-j}$. The resulting geometric series converges; the terminal term tends to zero because the larger piece avoids $O$ and its supremum is finite. This proves~\eqref{eq:Carleson} without an assumed bound on artificial boundary data.
\end{proof}

\begin{proof}[Proof of Theorem~\ref{thm:finiteabs}]
We will use the following vanishing-coefficient compactness. If
$\Delta z_j=\eps_j z_j^{-p}$ with $\eps_j\to0$ and the $z_j$
are locally bounded, then~\eqref{eq:Bernstein} gives a locally
uniformly convergent subsequence. Its limit is nonnegative and
harmonic on its positive set. If that set has an interior boundary
point, a ball in the positive set tangent to a zero point gives
linear lower growth by Hopf's lemma. The local Lipschitz bound on
$z_j^{1/\gamma}$ gives growth $O(d^\gamma)$ on the same radius,
a contradiction since $\gamma>1$. Thus on a connected limit domain
the limit is either zero or positive harmonic.

Fix $1<\zeta<\gamma$, so $p\zeta<1$. A positive harmonic function in a sufficiently large half-disc, zero on its diameter and equal to one at $e_2$, has value at $e_2/2$ strictly larger than $2^{-\zeta}$. To see the uniform assertion, use the odd Poisson kernel on a half-disc of radius $L$. Its value ratio is bounded below by $(L^2-1)/(2(L^2-1/4))$, which tends to $1/2$.

It follows by compactness, using Lemmas~\ref{lem:localabs} and~\ref{lem:additive}, that there are $T,L$ and a flatness tolerance such that
\begin{equation}\label{eq:normaldescent}
 u(P+d\nu)\ge T a^{1/(1+p)}d^\gamma
 \quad\Longrightarrow\quad
 u(P+\tfrac d2\nu)\ge2^{-\zeta}u(P+d\nu)\,.
\end{equation}
Here $\nu$ is the inward normal, and a fixed multiple $Ld$ of the boundary piece must avoid $O$ and other components. Indeed, a failure as the scaled curvature and $ad^2/u(P+d\nu)^{1+p}$ tend to zero would give the harmonic counterexample just excluded. The necessary uniform bound follows from~\eqref{eq:Carleson}. Boundary Lipschitz compactness follows from~\eqref{eq:localabs1}. The compactness argument is valid up to the true zero boundary.

Once~\eqref{eq:normaldescent} holds, it can be iterated: the ratio to $a^{1/(1+p)}d^\gamma$ increases by $2^{\gamma-\zeta}$ under halving. Interpolation by~\eqref{eq:addH} gives
\begin{equation}\label{eq:descentpower}
 u(P+d\nu)\ge C u(P+d_0\nu)(d/d_0)^\zeta,\qquad 0<d\le d_0\,.
\end{equation}

Write $U=u\circ F$, $t=\log(1/r)$, and
\begin{equation*}
 A(r)=N^{-1}\int_I U(r,\varphi)\psi\dd\varphi,\qquad B(t)=r^\beta A(r)\,.
\end{equation*}
On every fixed annulus Lemma~\ref{lem:localabs} makes the source integrable. Projection of the equation is therefore legitimate, first distributionally, and gives
\begin{equation}\label{eq:momentabs}
 B''+2\beta B'=\mathcal F(t)>0,\qquad
 \mathcal F(t)=e^{-d_0t}J(t),\quad d_0=2+\beta(1+p),
\end{equation}
where $J=N^{-1}\int_Ia|F'|^2(r^\beta U)^{-p}\psi$. The quantity $e^{2\beta t}B'$ is increasing, so $B$ is eventually monotone and has a limit in $[0,\infty]$. If the limit were zero, $B'<0$ eventually and integration of $B'(t)\ge B'(T)e^{-2\beta(t-T)}$ would give $B(t)\le Ce^{-2\beta t}$. Harnack chains and~\eqref{eq:Carleson} compare the annular supremum with $A(r)+a^{1/(1+p)}r^\gamma$. Thus $u\to0$ at $O$, contradicting the maximum principle for the positive subharmonic function $u$. Consequently $B(t)\ge b_0>0$ for large $t$.

Interior chains now give $U\ge cr^{-\beta}B(t)$ away from the angular endpoints. At a side choose an initial normal depth $d_0'=\delta r$. Its value is at least $cr^{-\beta}B(t)$, much larger than the natural scale $(d_0')^\gamma$. Curvature tends uniformly to zero after scaling. Apply~\eqref{eq:descentpower} on all these normals to obtain
\begin{equation}\label{eq:angularlower}
 U(r,\varphi)\ge cr^{-\beta}B(t)\psi(\varphi)^\zeta,\qquad
 U^{-p}\le Cr^{\beta p}\psi^{-p\zeta}\,.
\end{equation}
The latter is integrable in the whole corner since $p\zeta<1$. Together with~\eqref{eq:localabs2} this proves $u^{-p}\in L^1(\Om)$. It also gives $J(t)\le CB(t)^{-p}$, hence $\mathcal F\in L^1(T,\infty)$. Variation of constants in~\eqref{eq:momentabs} shows that $B$ has a finite positive limit.

The Green potential is finite at every interior point. The exhaustion by domains avoiding $O$ gives a positive harmonic function $h=u+a\G[u^{-p}]$. We check its ordinary boundary trace explicitly. On a smaller ordinary piece the near part of the source has an $L^1$ Dirichlet potential with zero $W^{1,q}$ trace, $1<q<2$. The remote potential is bounded by $Cd(x)$ times its value at a fixed interior point, by Green comparison. Thus $h$ is locally a zero-trace $W^{1,q}$ harmonic function, and boundary regularity makes it continuous and zero there. Its Martin measure is supported at $O$; hence $h=cK$ and $c>0$.

It remains to prove the relative limit along paths very close to a side. Choose $0<s<\min\{\beta,2+\beta p\}$ and solve $-f''-s^2f=\psi^{-p\zeta}$ with zero endpoints. Since $p\zeta<1$, $f\asymp\psi$. By~\eqref{eq:angularlower}, a large multiple of $r^sf$ is a supersolution for $w=a\G[u^{-p}]\circ F$ near $O$ and dominates its outer circular values. To justify comparison at the puncture, approximate $au^{-p}$ by bounded compactly supported increasing sources. Their potentials $w_n\le w\le cK\circ F$ have outer constants independent of $n$, while each is regular harmonic of order $r^\beta\psi$ sufficiently near $O$. Since $s<\beta$, compare on punctured sectors and pass first to zero inner radius and then to $n\to\infty$. Thus $0\le w\le Cr^s\psi$, and division by~\eqref{eq:Martin} proves the all-path relative limit. It also identifies $\lim B=c\kappa$. The uniqueness of $c$ follows at once.
\end{proof}

\subsection{Two Green potentials: existence and conditional uniqueness}
\citationadded{The elementary construction is~\cite[Theorem~3.5(ii)]{GV}, with kernel $G$, $h=cK$, $V=a$, and exponent $-p$. Its hypothesis is $aMc^{-p-1}\le p^p/(p+1)^{p+1}$, equivalently $c\ge C_0$, and its conclusion includes $u\ge s_pcK$. Here $C_0$ and $s_p$ are defined below. Lemma~\ref{lem:negativepotential} supplies the boundary behavior needed for its integral solution to lie in the natural class. The following second-potential argument improves this bound.} We use $H,M,J,\vartheta$ and the constants from~\eqref{eq:Jgain}--\eqref{eq:twoconst}.

\begin{lemma}\label{lem:gain}
One has $0<\vartheta<1$.
\end{lemma}
\begin{proof}
Write $f=K^{-p}$, $h=H/K$, and choose a cutoff $\chi$ supported near $O$, equal to one in a smaller neighborhood, where $h\le M/2$. Put $H_0=\G[\chi f]$ and $H_1=H-H_0$. Then
\begin{equation*}
 MH-J=\G[f(M-h)]\ge(M/2)H_0\,.
\end{equation*}
There is a global constant $C$ such that $H_1\le CH_0$. For completeness choose an interior ball $B$ in the region where $\chi=1$ and put $P=\G[\one_B]$. Then $H_0\ge bP$ for some $b>0$. Near $O$, both $P$ and $H_1$ are regular positive harmonic functions and are comparable to $r^\beta\psi$. On other boundary pieces $H_1\le Cd$ by $p<1$, whereas $P\ge cd$ by Hopf's lemma. The remaining compact part is immediate from $P>0$. Therefore $H_0\ge H/(1+C)$, and $J/(MH)\le1-1/(2(1+C))<1$. Positivity gives the other inequality.
\end{proof}

Then $C_G/C_0=(B_p/F_p)^{1/(p+1)}<1$.

\begin{proof}[Proof of Theorem~\ref{thm:twoGreen}]
Use the Banach space $Y_H$ with norm $\sup|z|/H$ and the interval $0\le z\le\overline z=H/((p+1)M)$. On this interval $K-z\ge s_pK$. The chord inequality
\begin{equation*}
 (1-t/(p+1))^{-p}\le1+(F_p-1)t,\qquad 0\le t\le1,
\end{equation*}
and~\eqref{eq:Jgain} imply $\tau_c\G[(K-\overline z)^{-p}]\le\tau_cB_pH\le\overline z$. The map thus preserves the interval. For two of its elements the mean value theorem gives
\begin{equation*}
 \abs{\tau_c\G[(K-z)^{-p}-(K-\widetilde z)^{-p}]}
 \le p\tau_cs_p^{-p-1}\G[K^{-p-1}|z-\widetilde z|]
 \le q_c\norm{z-\widetilde z}_{Y_H}H\,.
\end{equation*}
At $c=C_G$ the constant is $q_G<1$. Contraction gives the fixed point, and order preservation gives monotonicity of~\eqref{eq:iterate}. The first increment has norm $\tau_c$; summing subsequent increments proves~\eqref{eq:error}. Local positive lower bounds give classical regularity. Lemma~\ref{lem:negativepotential} gives the boundary data and trace.

If another solution with trace $c$ satisfies $u\ge s_pcK$, Theorem~\ref{thm:finiteabs} gives $\widetilde z=K-u/c\le\tau_cs_p^{-p}H$. It belongs to $Y_H$ even if it lies outside the smaller invariant interval. The same contraction estimate applies to it, proving $u=U_c$. Starting the original iteration at $cK$ dominates every represented solution at every step, so $U_c$ is maximal. Iterating simultaneously for two traces gives a difference at least $(c_2-c_1)K$. Subtract the limiting representations and use strict positivity of the Green kernel to obtain the strict inequality.
\end{proof}

\begin{corollary}\label{cor:energy}
Let $V_c=apU_c^{-p-1}$. For $c\ge C_G$ and $\xi\in H^1_0(\Om)$,
\begin{equation}\label{eq:energybound}
 \int_\Om V_c\xi^2\le q_c\int_\Om|\nabla\xi|^2\,.
\end{equation}
Any two solutions with the same trace and with $u,v\ge\ell K$, where $\ell>(ap\vartheta M)^{1/(p+1)}$, coincide.
\end{corollary}
\begin{proof}
The bound $\G[V_cH]\le q_cH$ gives a Schur estimate. With weight
$\sqrt{V_c}H$, the operator $\sqrt{V_c}\,\G\,\sqrt{V_c}$ has norm at most $q_c$. Apply it to $\sqrt{V_c}\xi$, or equivalently test the Dirichlet equation for $\G[V_c\xi]$, to obtain~\eqref{eq:energybound}. Truncate the potential and exhaust the domain before passing to the limit; no boundedness of $V_c$ at the boundary is required. For the last assertion use the same $Y_H$ estimate on the difference, with contraction constant $ap\ell^{-p-1}\vartheta M<1$.
\end{proof}

\begin{remark}\label{rem:uniquenesslimit}
Theorems~\ref{thm:finiteabs} and~\ref{thm:twoGreen} do not prove uniqueness in every fixed-trace class. In particular finiteness and a positive relative limit at $O$ do not imply the global lower bound required in Theorem~\ref{thm:twoGreen}. The minimum trace and uniqueness there are addressed in Theorem~\ref{thm:inverse-mass-branch}.
\end{remark}

\subsection{The minimum Martin trace and the maximal branch}

\begin{lemma}\label{lem:inverse-mass-integrability}
For fixed $a>0$ and $0<p<1$, there is an $r>1$ such that
$\norm{u^{-p}}_{L^r(\Om)}$ is bounded uniformly over all positive
solutions of~\eqref{eq:negabs}. A locally bounded sequence of stable
positive solutions has a subsequence converging locally uniformly to
a strictly positive function.
\end{lemma}
\begin{proof}
Use the pointwise estimate in the proof of Lemma~\ref{lem:Bernstein}.
First scale the ball and the solution so
that the radius and the equation coefficient are fixed. The estimate
bounds $|\nabla\log(1+u^{1/\gamma})|$ on a smaller concentric
ball. Thus either $u\ge1$ there, or $u$ has a uniform upper bound.
In the latter case test $\Delta u=u^{-p}$ with
$u^{-\eps}\phi^2$, where $\phi$ is a cutoff in that ball. Integration
by parts gives
\[
 \int u^{-p-\eps}\phi^2
 =\eps\int u^{-1-\eps}|\nabla u|^2\phi^2
   -2\int u^{-\eps}\phi\,\nabla u\cdot\nabla\phi.
\]
The first term is bounded by $C_p\eps\int u^{-p-\eps}\phi^2+C_p$
for $\eps<1$. The second is bounded independently of $u$ for
$\eps<(1-p)/2$, by the pointwise gradient estimate and the uniform
upper bound on $u$. Absorption for small $\eps>0$ and rescaling
give an integral bound $C R^{2-\gamma(p+\eps)}$ on the interior
quarter of the original ball. This bound also holds in the normalized
case $u\ge1$. Choose $\eps$ with $\gamma(p+\eps)<1$.
A Whitney decomposition of the planar Lipschitz domain has
$O(R^{-1})$ squares of side $R$ at each dyadic scale, since the
boundary layer of thickness $CR$ has area $O(R)$. Cover each
square by a bounded number of the interior balls just estimated.
Summing gives the assertion with $r=1+\eps/p$.

For a locally bounded sequence, Lemma~\ref{lem:Bernstein} bounds
the local Lipschitz norms of the powers $u^{1/\gamma}$. A diagonal
subsequence therefore converges locally uniformly to a nonnegative
function $v$. If $v(Q)=0$ at an interior point, the same estimate
gives $v(x)\le C|x-Q|^\gamma$ near $Q$. On each fixed annulus
about $Q$ contained in that neighborhood, local convergence gives
$ap u^{-p-1}\ge c_0|x-Q|^{-2}$ for the approximating solutions
sufficiently far along the sequence, with $c_0>0$ independent of the
annulus. A smooth radial test on an annulus of sufficiently large
logarithmic width has
$\int (|\nabla\xi|^2-c_0|x-Q|^{-2}\xi^2)<0$, by the
change of variable $t=\log|x-Q|$. This contradicts stability.
\end{proof}

\begin{proof}[Proof of Theorem~\ref{thm:inverse-mass-branch}]
Write $T_c(v)=cK-a\G[v^{-p}]$. By Theorem~\ref{thm:finiteabs},
every natural solution of trace $c$ is a positive fixed point of
$T_c$ and lies below $cK$. Conversely, a fixed point below $cK$
and bounded below by a positive continuous function is classical
by local elliptic regularity and vanishes continuously off $O$.
It has trace $c$ by Theorem~\ref{thm:finiteabs} and uniqueness of
the harmonic part of the representation. In particular, if
$0<v\le cK$ and $T_c(v)\ge v$ with $v$ positive and continuous,
iteration from $cK$ decreases to a fixed point by monotone
convergence of the potentials. Whenever $c$ is admissible,
iteration from $cK$ dominates every fixed point, so its limit
is the greatest solution $U_c$. From $u\le cK$ the fixed point
equation gives $cK\ge ac^{-p}\G[K^{-p}]$. This proves the
lower bound in~\eqref{eq:inverse-mass-bounds}, and
Theorem~\ref{thm:twoGreen} gives the upper bound.

If $c'>c$ and $U_c$ exists, then
$T_{c'}(U_c+(c'-c)K)\ge U_c+(c'-c)K$. This gives the asserted
order and shows that the admissible traces form an upper interval.
For $c>c_*$ choose an admissible $c_0<c$ and set $w=U_c-U_{c_0}>0$.
Convexity gives $-\Delta w\ge apU_c^{-p-1}w$. Testing this
inequality with $\xi^2/w$ and completing the gradient square
proves~\eqref{eq:inverse-mass-stability}. Along $c_j\downarrow c_*$
the maximal solutions are decreasing and locally bounded, and
Lemma~\ref{lem:inverse-mass-integrability} excludes zeros of the
limit. Monotone convergence gives its Green equation, and the
ordering with each possible solution of trace $c_*$ makes the
limit maximal. Stability passes to the limit on compact tests.

If $u$ is stable at trace $c$ and $w=U_c-u$, then
$w=a\G[u^{-p}-U_c^{-p}]\ge0$. The source belongs to $L^r$ by
Lemma~\ref{lem:inverse-mass-integrability}. Since $r>1$ in
dimension two, it belongs to $H^{-1}(\Om)$ and $w\in H^1_0(\Om)$.
If $w\ne0$, strict convexity and stability imply
\[
 \int_\Om|\nabla w|^2
 =a\int_\Om(u^{-p}-U_c^{-p})w
 <ap\int_\Om u^{-p-1}w^2
 \le\int_\Om|\nabla w|^2.
\]
The last inequality follows from compactly supported tests by
energy approximation and Fatou's lemma. This proves uniqueness in
the stable class.

At an interior trace $c>c_*$, use $c_0<c$ as above and set
$A_c z=ap\G[U_c^{-p-1}z]$. The Green equations give
$A_cw\le w-(c-c_0)K\le(c_0/c)w$, with
$(c-c_0)K\le w\le cK$. Thus $I-A_c$ has a positive bounded
inverse in the equivalent norm $\sup|z|/w$. The fixed point map is
smooth on the open cone of functions comparable to $K$ in the
$K$-norm, because its power derivatives applied to increments
bounded by multiples of $K$ are dominated by multiples of $K^{-p}$.
The implicit function theorem gives a smooth solution curve.
For $c'>c$, convexity gives
$0\le U_{c'}-U_c\le(c'-c)K+A_c(U_{c'}-U_c)$, proving right
continuity in the $K$-norm. For nearby $c'<c$, the implicit curve is
a lower bound for $U_{c'}$, while $U_{c'}\le U_c$, which gives
left continuity. The curve is therefore the maximal branch.
Strict concavity follows because a strict convex combination of
maximal solutions at distinct traces is a strict lower bound
for its image under $T_c$ at the interpolated trace $c$.
Indeed the two solutions are strictly ordered and $s\mapsto s^{-p}$
is strictly convex. Iteration from $cK$ proves the claim, including
when one of the traces equals $c_*$.

It remains to prove uniqueness at the minimum trace without
assuming stability of its competitors. Set $U_*=U_{c_*}$.
Fix $E\Subset\Om$ and $L>c_*$. Choose
$0\le\chi\le1$ in $C_c^\infty(\Om)$ with $\chi=1$ on $E$.
For $u$ of trace $c\in[c_*,L]$, put
\[
 W=\tfrac12(u+U_*),\qquad
 g=\chi(u-U_*)^2,\qquad D_g=\int_\Om g\dd x.
\]
On $\supp\chi$, both solutions are bounded above by a fixed
multiple of $K$. The second derivative of $s^{-p}$ has a positive
lower bound on the corresponding range. Hence
$\Delta W-aW^{-p}\ge b g$ for some $b>0$ independent of $u$.
Theorem~\ref{thm:finiteabs} gives $W\ge U_*/2\ge b_0K$ on a
neighborhood of $O$, for some $b_0>0$.

Assume $D_g>0$. Choose $B$ so large that replacing $W$ by
$W-\eps Z$, wherever $W\ge b_0K$ and $0\le Z\le2K$, increases
$aW^{-p}$ by at most $\eps BK^{-p}$ for
$0<\eps\le b_0/4$. Lemma~\ref{lem:negativepotential} lets us
choose $D=\Om\cap B_\delta(O)$ disjoint from $\supp\chi$,
with $W\ge b_0K$ and $BH\le K$ in $D$. Set
\[
 F_D=\one_D K^{-p},\qquad
 Z=K+B\G[F_D]-\frac{A}{D_g}\G[g].
\]
For $x\notin D$ and $y\in\supp\chi$, the quotient
$G(x,y)/K(x)$ has a positive uniform lower bound. Near the
ordinary boundary this follows from boundary Harnack and interior
Harnack in the source variable, and on interior compact subsets it
follows from positivity and the diagonal singularity. Since
$\G[F_D]\le H\le MK$, a fixed $A>0$ makes $Z\le0$ off $D$.
In $D$ we have $Z\le2K$.

Choose $\eps=\kappa_0D_g$ with $\kappa_0>0$ small enough that
$\eps\le b_0/4$, $\eps<c_*/2$ and $\kappa_0A\le b$.
The bound on $c$ makes $D_g$ uniformly bounded. Then
$W_\eps=W-\eps Z>0$ and
$\Delta W_\eps\ge aW_\eps^{-p}$ in distributions. In $D$,
the added Laplacian $\eps BF_D$ pays for the possible
increase of the power. Off $D$, $W_\eps\ge W$, and the
convexity surplus pays for the subtracted term $\kappa_0Ag$.
The Green representations give exactly
\[
 W_\eps=\bar cK-\G[\Delta W_\eps],\qquad
 \bar c=\tfrac12(c+c_*)-\eps.
\]
Consequently $W_\eps\le\bar cK$ and
$T_{\bar c}(W_\eps)\ge W_\eps$, so $\bar c$ is admissible.
Thus $\kappa_0D_g\le(c-c_*)/2$. If $D_g=0$ the estimate is
immediate. This proves~\eqref{eq:inverse-mass-local-rate} and
uniqueness at $c_*$.

Finally, Lemma~\ref{lem:Bernstein} and the upper bound by $LK$
give local equicontinuity. The integral estimate thus gives local
uniform convergence to $U_*>0$. Interior regularity then bounds
the local H\"older norms of the coefficient of the linear equation
for $u-U_*$. The interior $L^2$-to-supremum and Schauder estimates,
on nested interior neighborhoods, give the stated $C^{2,\eta}$
rate.
\end{proof}

\subsection{Ordinary boundary derivatives and degenerate points}
We first prove the $C^{1,\alpha}$ estimate in Theorem~\ref{thm:C1abs}
and distinguish positive and zero boundary flux. The mixed endpoint in
Theorem~\ref{thm:mixed} then follows from a Weiss monotonicity formula
and shallow-layer asymptotics at degenerate points.

\begin{proof}[Proof of Theorem~\ref{thm:C1abs}]
For the estimate~\eqref{eq:C1abs}, first rescale $D_{8R}(P)$ by $R$ and divide $u$ by $M+a^{1/(1+p)}R^\gamma$. In these units the solution is bounded by one and the equation coefficient is at most one. We use the same notation in the normalized patch.

We first record an interior consequence of~\eqref{eq:Bernstein}. For bounded positive solutions with bounded coefficient, $u^{1/\gamma}$ is uniformly Lipschitz. If $u(x)=h$, a ball of radius $ch^{1/\gamma}$ has $u\asymp h$; rescaling and interior Schauder estimates give $|D^2u(x)|\le C(1+h^{-p})$. For two points at distance $s$, either one value is $O(s^\gamma)$, in which case both gradients are $O(s^\alpha)$, or the values along the segment are bounded below by $Cs^\gamma$. Integration of the Hessian bound then gives $Cs^{1-\gamma p}=Cs^\alpha$. This proves a uniform interior $C^{1,\alpha}$ estimate without a positive infimum.

At a boundary point put $H_r=u(P+r\nu)$. If
$H_r\ge T a^{1/(1+p)}r^\gamma$ for $T$ large, additive Harnack
at a fixed scaled depth gives a comparable initial value on
neighboring normals. Apply~\eqref{eq:descentpower} on those
normals. After rescaling by $r$ and $H_r$, this gives
$v\ge c d^\zeta$ on a smaller boundary patch, where $d$ now
denotes the scaled distance. The scaled forcing is bounded by
$C\eps d^{-p\zeta}$ with $\eps=ar^2/H_r^{1+p}$. On a smaller $C^{1,1}$ auxiliary cap within this patch let $h$ be the harmonic lift of $v$ and $w$ the
zero-Dirichlet potential of its forcing, so $v=h-w$. The
normalized interior lower bound and the local upper bound from
\eqref{eq:Carleson} give $h\asymp d$ near the central boundary
piece by boundary Harnack.

For clarity, a bound $C\eps d^{-q}$ on the forcing, $0<q<1$,
gives $|\nabla w|\le C\eps$ there. The Green gradient kernel is
bounded by $C|x-y|^{-1}$ and the source weight has integral at most
$C\rho^{2-q}$ in any ball of radius $\rho$ in the scaled patch.
Thus $0\le w\le C\eps d$ near the zero boundary. Taking $T$
large gives $v\asymp d$ and improves the forcing bound to
$C\eps d^{-p}$. If $s=|x-x'|$, the part of the gradient
difference within distance $Cs$ of either target is
$O(\eps s^{1-p})$. On farther annuli use the kernel difference
bound $Cs^\mu\rho^{-1-\mu}$ for a fixed $1-p<\mu<1$ and
annular radius $\rho$. This is the local interior and zero-boundary
$C^{1,\mu}$ estimate on balls of radius comparable to $\rho$
that stay away from the source point. The gradient bound gives bounded oscillation on such balls. At the true side one uses the zero boundary value in place of subtracting an interior constant. The resulting dyadic sum is
bounded by $C\eps s^\mu\sum_{\rho\ge s}\rho^{1-p-\mu}
\le C\eps s^{1-p}$. The harmonic lift is also $C^{1,\mu}$
on the smaller piece. Therefore the normalized $v$ is
$C^{1,1-p}$ there with a uniform norm, and
\begin{equation}\label{eq:highflux}
 b(P)\asymp H_r/r,\qquad
 \sup_{D_t(P)}|\nabla u-b(P)\nu(P)|\le C(H_r/r)(t/r)^{1-p},\quad t<c r\,.
\end{equation}
The reverse upper bound uses the harmonic comparison in Lemma~\ref{lem:localabs}.

Fix a sufficiently small $r_0$ for the preceding comparisons. At a low scale $r\le r_0$, additive Harnack and~\eqref{eq:Carleson} bound the solution on the patch of radius $4r$ by $Cr^\gamma$. Lemma~\ref{lem:localabs} then gives $|\nabla u|\le Cr^\alpha$ on the patch of radius $r$. If all these scales are low, the gradient tends to zero at $P$ and $b(P)=0$. If any scale is high,~\eqref{eq:highflux} gives a positive boundary derivative. After this derivative is defined, the low-scale bound also gives $b(P)\le C_{\mathrm{low}}r^\alpha$ at each low scale.

For $b(P)>0$ take
\begin{equation*}
 r_b=\min\{r_0,(b(P)/(2C_{\mathrm{low}}))^{1/\alpha}\}\,.
\end{equation*}
This is a high scale. We have $b(P)\le Cr_b^\alpha$, using the fixed-patch gradient bound if $r_b=r_0$. For $t<c_h r_b$, where $c_h$ is the radius factor in~\eqref{eq:highflux}, that estimate bounds the gradient difference by $Cr_b^\alpha(t/r_b)^{1-p}\le Ct^\alpha$. If $t\ge c_h r_b$, apply the high or low estimate at $r=L_h t$ with fixed $L_h>1/c_h$, provided $r\le r_0$. In the high case use $H_r/r\asymp b(P)\le Ct^\alpha$. In the low case both the gradient and $b(P)$ are $O(r^\alpha)$. The fixed-patch bound handles the remaining $t$, which are bounded away from zero in normalized units. The same low-scale estimate applies when $b(P)=0$. Integrating the gradient estimate along a path of length $O(t)$ inside the boundary chart gives the affine approximation. In particular,
\begin{equation}\label{eq:Taylorabs}
 |u(x)-b(P)(x-P)\cdot\nu(P)|\le C|x-P|^\gamma,\qquad
 |\nabla u(x)-b(P)\nu(P)|\le C|x-P|^\alpha\,.
\end{equation}
Comparison through an interior point at the distance between two boundary points gives the $C^{0,\alpha}$ bound on $b\nu$. Interior pairs at separation much smaller than their boundary distance are handled by the interior estimate above if $b(P)\le Cd^\alpha$. In the other case subtract the affine tangent function, use~\eqref{eq:Taylorabs}, and apply interior Schauder estimates to get $|D^2u|\le Cd^{\alpha-1}$. Scaling back gives~\eqref{eq:C1abs} and the local estimates~\eqref{eq:Taylorabs} in the original variables.

Suppose now $b(P)=0$. The upper bounds $u=O(r^\gamma)$ and $\nabla u=O(r^\alpha)$ give compactness of the rescalings. They cannot converge to zero: in an interior ball of radius $cr$, comparison with a quadratic of Laplacian $a(\sup_{D_r}u)^{-p}$ gives $(\sup_{D_r}u)^{1+p}\ge ca r^2$.

For completeness the required homogeneity follows from an identity valid for these solutions, without minimization. Put $F_0(u)=au^{1-p}/(1-p)$, $\Sigma_r=\Om\cap\partial B_r(P)$, and
\begin{equation}\label{eq:Weiss}
 W(P,r)=r^{-2\gamma}\int_{D_r(P)}(|\nabla u|^2+2F_0(u))
             -\gamma r^{-2\gamma-1}\int_{\Sigma_r}u^2\,.
\end{equation}
Integration of the equation and of the stress identity
\begin{equation*}
 \operatorname{div}\big((x-P)(|\nabla u|^2+2F_0)-2((x-P)\cdot\nabla u)\nabla u\big)=4F_0
\end{equation*}
gives, with $n$ the outward normal,
\begin{equation}\label{eq:Weissder}
 W'(P,r)=2r^{-2\gamma}\int_{\Sigma_r}(\partial_ru-\gamma u/r)^2
    +r^{-2\gamma-1}\int_{\partial\Om\cap B_r(P)}((x-P)\cdot n)b^2\,.
\end{equation}
All boundary integrations follow by first removing a boundary strip: Lemma~\ref{lem:localabs} gives the integrable source, and the already proved $C^1$ estimate gives convergence of the boundary terms. This is the Weiss identity~\cite{Weiss} with its explicit curvature term. At a zero-flux point that term is bounded, since $|(x-P)\cdot n|\le C|x-P|^2$ and $b=O(r^\alpha)$. Thus $W+Cr$ is increasing and has a finite limit at zero. Every rescaling limit $U$ has $\partial_rU=\gamma U/r$ and is $r^\gamma\Phi(\varphi)$.

The nonnegative $C^1$ angular function on $[0,\pi]$ satisfies $\Phi''+\gamma^2\Phi=a\Phi^{-p}$ on its positive intervals. Its first integral is
\begin{equation*}
 (\Phi')^2+\gamma^2\Phi^2-\frac{2a}{1-p}\Phi^{1-p}=E\,.
\end{equation*}
At a zero endpoint $E\ge0$; at an interior zero endpoint $E=0$. If $M_\Phi$ is the maximum on a positive interval and $\tau=2a/((1-p)\gamma^2M_\Phi^{1+p})$, its length is
\begin{equation*}
 \frac2\gamma\int_0^1\big((1-s^2)-\tau(1-s^{1-p})\big)^{-1/2}\dd s\,.
\end{equation*}
It increases strictly from $\pi/\gamma$ to $\pi$ as $\tau$ increases from zero to one. The endpoint integral at $\tau=1$ follows by the substitution $z=s^{(1+p)/2}$. Hence the only nonzero profile on a half-plane is $\Phi=A_p\sin^\gamma\varphi$. A positive interval with an interior zero endpoint would itself require length $\pi$, so no other interval is possible. Every subsequence has this limit; positivity inside the half-plane gives interior $C^2$ convergence. This proves~\eqref{eq:absprofile}. Along the normal $|\nabla u(P+t\nu)|\sim A_p\gamma t^\alpha$, which proves optimality when $b(P)=0$.
\end{proof}

\subsection{The second-order endpoint on a mixed boundary piece}
Let $\Gamma$ be a compact ordinary boundary arc, and use $V_\delta(\Gamma)$ and $Z$ from the introduction.

\begin{proof}[Proof of Theorem~\ref{thm:mixed}]
The proof counts low-value scales on each normal before summing the boundary-layer estimates. We first derive a local level-set bound. For $v=u^{1/\gamma}$,
\begin{equation}\label{eq:vvariation}
 \int_{B_1}|\Delta v|\le C(a,p,M),\qquad 0<u\le M\quad\text{in }B_2\,.
\end{equation}
To prove this, set $A=a/\gamma$, $B=A/\alpha$, and $Q=(|\nabla v|^2-B)/v$, so $\Delta v=-\alpha Q$. Differentiating $Qv=|\nabla v|^2-B$ gives the exact identity
\begin{equation*}
 \Delta Q+\frac{2\gamma}{v}\nabla v\cdot\nabla Q
 =\frac2v|D^2v|^2+\frac\alpha v Q^2\,.
\end{equation*}
Take $\eta=\chi_0^2$, where $\chi_0$ is a smooth cutoff
supported in the larger ball and equal to one on a smaller ball.
Then $|\nabla\eta|^2/\eta$ is bounded wherever $\eta>0$.
At a positive maximum of $\eta Q$, the identity gives
\begin{equation*}
 \alpha\eta Q\le Cv+2\gamma\sqrt{B+vQ}\,|\nabla\eta|
 \le\frac\alpha2\eta Q+C_p(v+\sqrt B)\,.
\end{equation*}
Thus $Q_+$ is uniformly bounded in a smaller ball, and $\Delta v$ has a uniform lower bound. The gradient of $v$ is bounded by~\eqref{eq:Bernstein}. For a further cutoff $\chi$, integration of $\Delta v$ gives
\begin{equation*}
 \int\chi(\Delta v)_+
 =-\int\nabla\chi\cdot\nabla v+\int\chi(\Delta v)_-\le C\,.
\end{equation*}
This proves~\eqref{eq:vvariation}.

Testing $\Delta v$ with $\chi\min(v,t)$ now gives $\int_{\{v<t\}}\chi|\nabla v|^2\le Ct$. Integration of $v\Delta v+\alpha|\nabla v|^2=A$ over the same set gives $A\int_{\{v<t\}}\chi\le Ct$. Therefore
\begin{equation}\label{eq:interiorlevel}
 \abs{\{u<t^\gamma\}\cap B_{r/2}(x)}\le Crt
 \quad\hbox{if }\sup_{B_{2r}(x)}u\le M_0r^\gamma\,.
\end{equation}

Fix $0<\tau<A_p$. Call $(P,r)$ a low scale if $u(P+r\nu)\le\tau r^\gamma$. The Taylor estimates~\eqref{eq:Taylorabs} imply $b(P)\le Cr^\alpha$ and uniform $C^{1,\alpha}$ bounds for the rescaling at a low scale. There are $r_0,\eps_0>0$ such that every low scale below $r_0$ satisfies
\begin{equation}\label{eq:lowenergy}
 2\int_{r/2}^{2r}s^{-2\gamma}\int_{\Sigma_s(P)}(\partial_\rho u-\gamma u/s)^2\dd S\dd s\ge\eps_0\,.
\end{equation}
Otherwise rescaling a violating sequence would give a nonzero homogeneous half-annulus limit with value at $e_2$ at most $\tau$. The angular first-integral classification used in Theorem~\ref{thm:C1abs} gives only $A_p y_2^\gamma$, a contradiction. The limit is nonzero because on any interior ball $\Delta u_j\ge aC^{-p}$, which excludes uniform convergence to zero.

For a geometric sequence $r_j=r_0\vartheta_0^j$, each $P$ has at most a fixed number of low scales. Indeed take any finite collection with smallest radius $r_-$ and largest radius $r_+$. The bound $b(P)\le Cr_-^\alpha$ implies $b(P)\le Cs^\alpha$ for $r_-/2\le s\le2r_+$. For $Q\in\partial\Om\cap B_s(P)$, the H\"older bound on $b$ gives $b(Q)\le b(P)+C|Q-P|^\alpha\le Cs^\alpha$. Thus the curvature error in~\eqref{eq:Weissder} is bounded on the whole interval. The estimates~\eqref{eq:Taylorabs} also bound both endpoint energies. Integrate~\eqref{eq:Weissder} and use the bounded overlap of $[r_j/2,2r_j]$. Each selected scale contributes at least $\eps_0$, so their number is uniformly bounded. This does not assume a finite limit of $W(P,r)$ at a positive-flux point.

Choose a small $\eta>0$ and take $\vartheta_0=(1+\eta)^{-1}$.
Divide a slightly larger boundary arc into intervals $I_Q$ of
length between $\eta r_j/2$ and $\eta r_j$. Form the normal boxes
$Q$ above these intervals with depths in $[r_j,(1+\eta)r_j)$.
Mark $Q$ if it contains a point $P+s\nu(P)$ with
$u(P+s\nu(P))\le(A_p/2)r_j^\gamma$. At that point
\eqref{eq:Taylorabs} gives $b(P)\le Cr_j^\alpha$, since
$s\asymp r_j$. The gradient estimate in~\eqref{eq:Taylorabs}
therefore gives $|\nabla u|\le Cr_j^\alpha$ throughout the box.
Any point at depth $r_j$ above $I_Q$ is joined to the marked
point by a normal and tangential path of length $O(\eta r_j)$
in the box. For $\eta$ small the value changes by at most
$(A_p/4)r_j^\gamma$. Thus every point of $I_Q$ has a low
scale at $r_j$ with $\tau=3A_p/4$. The scale count integrated
over the boundary gives
\begin{equation}\label{eq:packing}
 \sum_{Q\ \mathrm{marked}}|I_Q|<\infty.
\end{equation}
For sufficiently large $L_0$, the condition $d>L_0t$ and
$u<t^\gamma$ ensures that a point in these boxes belongs to a
marked box, since $d<(1+\eta)r_j$. The part of the collar above
the top of these boxes lies in a fixed interior compact set and
contributes nothing for small $t$. The same Taylor estimate
anchored at the marked normal bounds $u$ by $Cr_j^\gamma$ on a
slightly larger interior ball of radius comparable to $r_j$
whose smaller concentric ball contains the box. Hence
\eqref{eq:interiorlevel} gives, with a constant allowed to
depend on the fixed $\eta$,
\begin{equation*}
 \abs{Q\cap\{u<t^\gamma\}}\le Ct|I_Q|\,.
\end{equation*}
Each fixed box is compactly contained in $\Om$ and its positive minimum excludes it when $t$ is sufficiently small. Dominated convergence in the summable series~\eqref{eq:packing} therefore gives
\begin{equation}\label{eq:deep}
 \abs{\{u<t^\gamma,\ d>L_0t\}\cap V_\delta(\Gamma)}=o(t)\,.
\end{equation}
The remaining layer has area $O(t)$, proving the weak estimate for $u^{-p}$.

The Hessian satisfies the useful pointwise bound
\begin{equation}\label{eq:Hmixed}
 |D^2u(x)|\le C\big(u(x)^{-p}+d(x)^{-m}\big)\,.
\end{equation}
If $b(P)\le Cd^\alpha$, Taylor's estimate bounds $u$ by $Cd^\gamma$ on a larger interior ball. The power gradient estimate on this ball gives comparability on a ball of radius $cu(x)^{1/\gamma}$ and Schauder gives $|D^2u(x)|\le Cu(x)^{-p}$. If $b(P)>Cd^\alpha$, subtract the tangent affine function and use its $O(d^\gamma)$ remainder on a ball of radius comparable to $d$. Schauder gives $|D^2u|\le C(d^{\gamma-2}+(b(P)d)^{-p})\le Cd^{-m}$. Since $d^{-m}\in L^{s_b,\infty}$, this proves the Hessian assertion and all the strong Sobolev statements below the endpoint.

Finally, let $d=tz$ in the shallow layer. For fixed $z>0$
and $P\in Z$, Theorem~\ref{thm:C1abs} gives
\[
 t^{-\gamma}u(P+tz\nu)\longrightarrow A_pz^\gamma,\qquad
 t^mD^2u(P+tz\nu)\longrightarrow
 aA_p^{-p}z^{-m}\nu\otimes\nu .
\]
For $P\notin Z$, the first ratio tends to infinity and the
second tends to zero. Indeed $b(P)>0$ and~\eqref{eq:highflux},
together with interior Schauder estimates, gives
$|D^2u|=O(d^{-p})$ on that normal, and $m>p$. The normal
coordinate Jacobian tends uniformly to one. Dominated convergence
for the value indicators on $\Gamma\times(0,L_0)$, with
$L_0>A_p^{-1/\gamma}$, together with~\eqref{eq:deep},
proves~\eqref{eq:mixedlow}.

For the Hessian take $t=L^{-s_b}$, so that $L=t^{-m}$.
By~\eqref{eq:Hmixed}, a point with $|D^2u|>L$ at depth larger
than a sufficiently large multiple of $t$ satisfies
$u<(C t)^\gamma$. By increasing the depth multiple,
\eqref{eq:deep} with $Ct$ in place of $t$ makes the measure of
this deep part $o(t)$. Choose the shallow depth multiple also
larger than $(aA_p^{-p})^{1/m}$. Dominated convergence for the
indicators $t^m|D^2u|>1$ now gives~\eqref{eq:mixedtail}.
Here the norm of $\nu\otimes\nu$ is one in the Frobenius
norm. Substituting $t=(a/L)^{s_b}$ in~\eqref{eq:mixedlow}
gives the source limit. A positive tail coefficient implies
logarithmic divergence of the endpoint integral and of all
finite-index Lorentz integrals.
\end{proof}

\begin{remark}
If $\HH^1(Z\cap\Gamma)=0$,~\eqref{eq:mixedtail} yields $o(L^{-s_b})$ but does not decide strong $L^{s_b}$ membership. The normal asymptotic at a fixed zero-flux point also does not imply a uniform pointwise bound $u\ge cd^\gamma$ throughout a mixed boundary neighborhood. These stronger assertions are not used in the proof.
\end{remark}

\section{Negative-power sources}\label{sec:negsource}
We now study~\eqref{eq:negsource}. On ordinary boundary pieces, the inward derivative is positive and finite for $0<q<1$. For $q\ge1$ the inward difference quotient tends to $+\infty$.

\subsection{Existence and uniqueness at a prescribed trace}
We first construct the zero-trace solution by regularization and then use the Martin kernel to prescribe $c$. Comparison proves uniqueness and exhaustiveness, while separate boundary and energy arguments identify the properties of the resulting family.

\begin{proof}[Proof of Proposition~\ref{prop:negsourceclass}]
Write $\tau=\G[1]$ for the torsion function and $\gamma=1/(q+1)$. It is bounded, positive in $\Om$, and continuous with zero values on the whole boundary by the exterior cone condition. For $A$ with $\gamma A^{q+1}\ge a$, the chain rule gives
\[
 -\Delta(A\tau^\gamma)
 =A\gamma\tau^{\gamma-1}
   +A\gamma(1-\gamma)\tau^{\gamma-2}|\nabla\tau|^2
 \ge a(A\tau^\gamma)^{-q}.
\]
Solve the regularized Dirichlet problem with bounded nonincreasing reaction $a(v_++\eps)^{-q}$ and zero data. For each $\eps>0$ a weak solution $u_\eps\in H^1_0(\Om)$ follows by minimizing the Dirichlet energy minus the primitive of this reaction. It is nonnegative and continuous up to the boundary because its bounded forcing gives $0\le u_\eps\le a\eps^{-q}\tau$. Interior regularity and the strong maximum principle give positivity. Pointwise comparison with $A\tau^\gamma$ bounds $u_\eps$ independently of $\eps$ and makes the boundary estimate valid at the corner. The solutions increase as $\eps\downarrow0$. Their positive lower bound on each interior compact subset follows by comparison with one fixed regularization. Local elliptic estimates give a classical positive limit $u_0$ with continuous zero values on the whole boundary.

On a compact ordinary boundary patch the classical singular Dirichlet estimates~\cite{CRT,LM} give
\begin{equation}\label{eq:ordinarysource}
 u_0\asymp
 \begin{cases}d,&q<1,\\d\sqrt{\log(1/d)},&q=1,\\d^{2/(q+1)},&q>1.\end{cases}
\end{equation}
These estimates also follow locally. On a sufficiently thin smooth patch $d$ is $C^2$ with bounded Laplacian. A positive harmonic lifting of a nonzero part of the artificial boundary data gives the lower bound $cd$. For $q<1$ the resulting source bound $Cd^{-q}$ has a potential bounded by $Cd$ on a smaller patch, by comparison with a positive multiple of $d-Bd^{2-q}$ and an outer harmonic lifting. For $q=1$ and $q>1$ use, respectively, $g(d)=d\sqrt{\log(C/d)}$ and $g(d)=d^{2/(q+1)}$. In a sufficiently thin collar $-\Delta g(d)\asymp g(d)^{-q}$, since the term $g'(d)\Delta d$ is lower order. A large multiple of $g$ plus the harmonic lifting of the artificial boundary data is a supersolution. A small multiple of $g$ minus the harmonic lifting of its own artificial boundary data is a subsolution on its positive set. On the smaller patch these liftings are $O(d)$ and $g(d)/d\to\infty$, giving the other two bounds.

For prescribed $c$ use the subsolution $\max\{u_0,cK\}$ and the supersolution $u_0+cK$ on smooth interior exhaustion domains, with the subsolution as boundary data. Comparison is valid since $t\mapsto at^{-q}$ is decreasing. The exhaustion solutions increase and remain between these barriers, giving a classical limit. Its nonnegative harmonic minorant is $cK$ because the difference from $cK$ is bounded by $u_0$ and the harmonic minorant is supported at $O$. This gives~\eqref{eq:negsourceclass}. The maximum principle for the difference gives uniqueness at fixed $c$, since the difference tends to zero on the whole boundary.

For an arbitrary natural solution $u$, the positive part of $u_0-u$ is subharmonic and has zero boundary upper limit, including at $O$ since $u_0\to0$ and $u\ge0$. Thus $u\ge u_0$. The Riesz decomposition has nonnegative harmonic part supported at $O$. Since $\G[u^{-q}]\le\G[u_0^{-q}]=u_0/a$, its remainder tends to zero on the whole boundary. This proves exhaustion of the natural class.

For $0<q<1$, the bound $u\ge u_0\gtrsim d$ gives a local source bound $a u^{-q}\le Cd^{-q}$. Write $u$ as a harmonic lifting plus its Green potential on a smaller ordinary cap. The local Green-gradient estimates used in~\eqref{eq:highflux}, with exponent $q$ in the source bound, give continuity of the gradient up to the side. The lower bound $u\gtrsim d$ makes the inward derivative positive. When $q\ge1$, the lower bound $u\ge u_0$ and~\eqref{eq:ordinarysource} give $u(P+t\nu(P))/t\to+\infty$.

For the energy criterion, testing the regularized equation with $u_\eps$ gives
\[
 \int_\Om|\nabla u_\eps|^2\dd x
 =a\int_\Om u_\eps(u_\eps+\eps)^{-q}\dd x.
\]
If $q\le1$, the integrand is uniformly bounded. If $1<q<3$, comparison for the unregularized equation gives $u_\eps+\eps\ge u_0$ because $u_\eps+\eps$ solves that equation and is positive on the boundary. By~\eqref{eq:sourceball} below,
\[
 u_\eps(u_\eps+\eps)^{-q}\le u_0^{1-q}
 \le C d^{-2(q-1)/(q+1)}\in L^1(\Om).
\]
The last inclusion follows from the linear volume bound for boundary layers of a bounded Lipschitz domain. Thus $u_0\in H^1_0(\Om)$ by weak compactness. If $q\ge3$, the lower bound in~\eqref{eq:ordinarysource} makes $u_0^2/d^2$ nonintegrable on any smooth boundary patch, contrary to Hardy's inequality for an $H^1_0$ function. For $c>0$, the inequality $u_c\ge cK$ on an interior angular subsector gives $u_c\notin L^s$ for finite $s\ge 2/\beta$. The planar Sobolev embedding on a bounded Lipschitz domain excludes $u_c\in H^1(\Om)$.
\end{proof}

The absolute remainder in Proposition~\ref{prop:negsourceclass} allows a positive-mass solution to approach its Martin profile on interior subsectors. The all-path ratio also depends on the side scale.

\subsection{Corner scales and joint derivative estimates}
Set $\alpha_q=2/(q+1)$ and write $U=u\circ F$. We first record two local comparisons that do not assume a corner profile. For every natural solution of~\eqref{eq:negsource},
\begin{equation}\label{eq:sourceball}
 u(x)\ge c_q a^{1/(q+1)}d(x)^{\alpha_q}\,.
\end{equation}
Indeed, in $B_\rho(x)$ with $\rho=d(x)/2$, the function $w(y)=k a^{1/(q+1)}\rho^{\alpha_q}(1-|y-x|^2/\rho^2)$ satisfies $-\Delta w\le aw^{-q}$ if $4k\le k^{-q}$. It is zero on the sphere, whereas $u>0$ there. Comparison for the decreasing reaction gives $u\ge w$, proving~\eqref{eq:sourceball}. In conformal coordinates this implies $U(r,\varphi)\ge cr^{\alpha_q}\psi(\varphi)^{\alpha_q}$.

The side scales $S_\eps$ and $T_\eps$ are those in~\eqref{eq:Ssource}--\eqref{eq:Tsource}.

\begin{lemma}\label{lem:sourceside}
Consider a family $V_j>0$ solving $-\Delta V_j=\eps_j A_jV_j^{-q}$ in a fixed half-disc, with continuous zero values on its straight diameter. Suppose $0<\eps_j\le C_0$, $0<c_0\le A_j\le C_0$, and on each interior compact set the functions $V_j$ have uniform positive lower bounds and finite upper bounds. On a smaller half-disc adjacent to an interior part of the diameter,
\begin{equation}\label{eq:sourceside}
 V_j\asymp S_{\eps_j}(y),
\end{equation}
where $y$ is the distance to the diameter. The constants use only fixed local interior bounds and the stated coefficient bounds. They are independent of $\eps_j\downarrow0$.
\end{lemma}
\begin{proof}
Use three nested half-discs, so every artificial circular boundary stays outside the final estimate. Let $h_j$ be the positive harmonic lifting of the outer half-disc data. At a fixed interior point $h_j\le V_j\le C$. Boundary Harnack gives $h_j\le Cy$ on the middle piece. A fixed segment of the artificial circular boundary is contained in an interior compact set and has $V_j\ge c$; its harmonic measure gives $h_j\ge cy$ there. On the middle half-disc use a new lifting $\widetilde h_j$. It is at least $h_j$, and on the final piece its upper bound follows in the same way. Hence $\widetilde h_j\asymp y$ where needed, and $\widetilde h_j\ge cy$ on the middle domain.

Write $V_j=\widetilde h_j+P_j$, with $P_j$ its positive zero-Dirichlet potential on that domain. If $q<1$, then $\eps_jA_jV_j^{-q}\le C\eps_j y^{-q}$. The Green potential of this density is bounded by $C\eps_j y$ on the final piece. Thus $V_j\asymp y$, using $\eps_j\le C_0$.

If $q>1$, the function $\widetilde h_j+C\eps_j^{1/(q+1)}y^{\alpha_q}$ is a supersolution: its negative Laplacian is $C\eps_j^{1/(q+1)}\alpha_q(1-\alpha_q)y^{\alpha_q-2}$, while the reaction decreases when the harmonic term is added. It dominates the artificial data. The ball comparison~\eqref{eq:sourceball}, with the lower coefficient bound, gives the matching lower fractional-power term. Together with $V_j\ge\widetilde h_j\ge cy$, this proves~\eqref{eq:sourceside}.

If $q=1$, put $g(y)=y\sqrt{\log(C_1/y)}$, with fixed $C_1$ large. Direct differentiation gives $-g''\asymp g^{-1}$. Thus $\widetilde h_j+C\sqrt{\eps_j}g$ is a supersolution. For the lower bound compare with the zero-boundary singular solution in a fixed ball tangent to the diameter. Its ordinary boundary estimate~\eqref{eq:ordinarysource} is $cy\sqrt{\log(C_1/y)}$; multiplication by $\sqrt{\eps_j c_0}$ gives a subsolution for the current coefficient. The radius of that ball is fixed, so the constant is independent of $\eps_j$. Combining this lower bound with $cy$ and the upper bound gives $V_j\asymp y\sqrt{1+\eps_j\log(e/y)}$.
\end{proof}

The amplitude $M(r)$ and the angular state $\Phi$ are defined in~\eqref{eq:sourceM}--\eqref{eq:negativeangular}. In the resonant row, $b$ tends to infinity and
\begin{equation}\label{eq:negsourceresonance}
 b(t)\asymp
 \begin{cases}
 t^{1/(q+1)},&0<q<2,\\
 (t\log t)^{1/3},&q=2,\\
 t^{(q-1)/(q+1)},&2<q<3.
 \end{cases}
\end{equation}
Resonance is possible only for $q<3$, since $\beta>1/2$. The convergence in~\eqref{eq:sourceM} is uniform on interior angular subsectors and locally $C^2$ there.

For $c>0$ the bound $0\le u_c-cK\le u_0$ gives the profile in~\eqref{eq:sourceM} after division by $M(r)=c\kappa r^{-\beta}$. On each interior compact subset of a normalized annulus, $U(r\cdot)/M(r)$ is bounded above and bounded away from zero, and its scaled forcing coefficient $ar^2/M(r)^{q+1}$ tends to zero. Interior Poisson estimates followed by the interior H\"older estimates give the local $C^2$ convergence. In the following three cases $c=0$ and $u=u_0$.

The sector reduction is related to~\cite{Wigley85}, but the singular side integral must be retained here. In the nonlinear range $\beta>\alpha_q$, the angular operator $-\partial_\varphi^2-\alpha_q^2$ has positive first Dirichlet eigenvalue. A small multiple of $\psi$ is a subsolution of~\eqref{eq:negativeangular}. A supersolution is $\psi^\eta+B\psi$, where $\eta=1/(q+1)$ and $B$ is large. Near either endpoint the shifted negative Laplacian of $\psi^\eta$ has leading positive term $\eta(1-\eta)\psi^{\eta-2}(\psi')^2$, which dominates the reaction there. The term $B(\beta^2-\alpha_q^2)\psi$ gives the inequality on the interior. Solve the regularized problems between these barriers and let the regularization tend to zero. The positive first eigenvalue and the decreasing reaction give comparison and uniqueness. For fixed small $\eps>0$, the functions $S_\pm=(1\pm\eps)r^{\alpha_q}\Phi$ are a subsolution and a supersolution in a sufficiently small sector $D_R$, since $|F'|^2\to1$ and $\alpha_q(q+1)=2$. Let $H_+$ be the harmonic lifting of $U$'s outer circular data, and let $H_-$ lift the outer data of $S_-$. Both liftings are zero on the straight sides and at the vertex. Decrease of the reaction gives the comparisons
\begin{equation}\label{eq:sourcelift}
 (S_--H_-)_+\le U\le S_++H_+\quad\text{in }D_R\,.
\end{equation}
For the lower comparison apply the maximum principle only where $S_--H_->U>0$. The outer data are ordered by construction, and all functions tend to zero at the vertex. For the upper comparison the added harmonic term decreases the reaction. On $r<R/2$, the regular harmonic expansion gives $H_\pm\le C_{\eps,R}r^\beta\psi$. The concave angular solution has $\Phi\ge c\psi$, even when its endpoint slopes are infinite. Therefore the harmonic corrections divided by $r^{\alpha_q}\Phi$ tend uniformly to zero. Let $r\downarrow0$ and then $\eps\downarrow0$ in~\eqref{eq:sourcelift}. This proves the nonlinear profile without importing an exact ordinary-boundary coefficient.

In the harmonic range $\beta<\alpha_q$, harmonic comparison and~\eqref{eq:sourceball} give $U\ge c(r^\beta\psi+r^{\alpha_q}\psi^{\alpha_q})$. To verify that the regular harmonic coefficient is finite, integrate the source against the dual mode $r^{-\beta}\psi$. For $q<2$ the angular moment is finite directly. For $q\ge2$ split the angular integral at $\psi=\delta(r)=r^{(\alpha_q-\beta)/(1-\alpha_q)}$. When $q>2$ the resulting radial integrability margin is $(q+1)(\alpha_q-\beta)/(q-1)>0$; when $q=2$ a positive radial power is multiplied by a logarithm. The dual moment is therefore finite. Write $U$ as the positive outer harmonic lifting plus its Dirichlet Green potential on $D_R$. At an observation point $x$ in a fixed interior angular sector, put $r=|x|$. For each fixed pole $z$, the first sine term gives the limit of $r^{-\beta}G_{D_R}(x,z)$ as $r\downarrow0$ along a ray. Off a ball $B_{\rho r}(x)$, with fixed small $\rho>0$, comparison with the full-sector Green kernel gives
\begin{equation*}
 r^{-\beta}G_{D_R}(x,z)\le C|z|^{-\beta}\psi(\arg z)\,.
\end{equation*}
Indeed the power map sends the full sector to the half-plane. With $X=x^\beta$ and $Z=z^\beta$, its Green kernel is bounded by $C\,\operatorname{Im}X\,\operatorname{Im}Z/|X-Z|^2$. Off $B_{\rho r}(x)$ the denominator is bounded below by a multiple of $|z|^{2\beta}$, by scaling and separation in the power chart. The finite dual moment now gives dominated convergence. Inside the ball the source is bounded by $Cr^{-\beta q}$, and the logarithmic half-plane kernel has integral $O(r^2)$ on that ball. The normalized contribution is $O(r^{2-\beta(q+1)})$, which tends to zero. Thus the interior annular limit is $B|z|^\beta\psi$, where $B>0$ because the outer harmonic lifting already contributes a positive first coefficient. Interior estimates give the stated $C^2$ convergence. The near-diagonal rate concerns only the local contribution.

At equality set $U=r^\beta v(t,\varphi)$ and $b=N^{-1}\int_Iv\psi$. In the cylinder equation test with $\xi(t)\psi(\varphi)\chi_\delta(\varphi)$, where $\xi$ is smooth and compactly supported and $\chi_\delta$ cuts off at angular distance $\delta$ from both endpoints, with transition bands between distances $\delta$ and $2\delta$. After moving angular derivatives onto the test, the terms from differentiating $\chi_\delta$ have integrands bounded by $C\delta^{-1}|v|$ on these bands. Their integrals tend to zero by the uniform zero trace of $v$ on each fixed strip. The reaction paired with $\psi$ is dominated by $C\psi^{\beta-1}$ by~\eqref{eq:sourceball} and $\beta(q+1)=2$, so dominated convergence justifies the projection without a finite side derivative. The first equation below holds distributionally and then classically since $J$ is continuous by the same domination. The future integral follows from the growth bound proved just below.
\begin{equation}\label{eq:negsourceb}
 -b''+2\beta b'=J(t),\qquad
 J=N^{-1}\int_I a|F'|^2v^{-q}\psi,\qquad
 b'(t)=\int_t^\infty e^{-2\beta(s-t)}J(s)\dd s\,.
\end{equation}
We establish slow variation before normalizing the solution. At resonance $\alpha_q=\beta$,~\eqref{eq:sourceball} gives $v\ge c\psi^\beta$, and hence
\begin{equation}\label{eq:negativeJbound}
 0<J(t)\le C\int_I\psi^{1-\beta q}\dd\varphi
     =C\int_I\psi^{\beta-1}\dd\varphi<\infty\,.
\end{equation}
The boundedness of $u_0$ gives $b(t)=O(e^{\beta t})$. Integrating $(e^{-2\beta t}b')'=-e^{-2\beta t}J$ shows that a nonzero homogeneous coefficient would produce an $e^{2\beta t}$ term in $b$, incompatible with this bound and positivity. This proves the future integral in~\eqref{eq:negsourceb}, and~\eqref{eq:negativeJbound} gives $0<b'\le C$. Jensen's inequality with weight $\psi$ also gives $J\ge cb^{-q}$. A bounded $b$ would therefore have $b'\ge c>0$, a contradiction. We have proved
\begin{equation}\label{eq:negativeslow}
 b(t)\longrightarrow\infty,\qquad b'(t)/b(t)\longrightarrow0,
 \qquad b(t+s)/b(t)\longrightarrow1
\end{equation}
uniformly for $s$ in a fixed bounded interval.

Now normalize at times $t_j\to\infty$. Set $r_j=e^{-t_j}$ and use $U(r_jz)/(r_j^\beta b(t_j))$ on a fixed sector annulus. Its source is bounded by $Cb(t_j)^{-1}\psi^{\beta-2}$ by~\eqref{eq:sourceball}. Let $g$ solve $-g''=\psi^{\beta-2}$ with zero endpoints. Its one-dimensional Green representation is bounded and continuous up to the endpoints, since the kernel supplies one power of distance and $\int_0^1y^{\beta-1}\dd y<\infty$. On the fixed annulus, $Cb(t_j)^{-1}g(\varphi)$ is a supersolution for the Dirichlet potential, because $-\Delta g(\varphi)=|z|^{-2}\psi^{\beta-2}$. Comparison with compact source truncations therefore gives a potential tending uniformly to zero.

The first angular integrals are bounded by~\eqref{eq:negativeslow}. The positive harmonic liftings are locally bounded on the interior by the mean-value estimate, using these weighted integrals, and up to the sides by harmonic boundary comparison. They give local compactness up to the zero sides. Interior estimates give $C^2$ compactness away from the sides. A translated normalized limit in the $v$ variables solves the homogeneous cylinder equation and has first sine coefficient identically one. Positivity and $|\sin(n\beta\varphi)|\le n\psi$ bound every other sine coefficient for all translated times. Its possible factors are $e^{(1+n)\beta s}$ and $e^{(1-n)\beta s}$. For $n\ge2$ they are unbounded in opposite directions, so their coefficients vanish. The only limit is $\psi$, which proves $v/b\to\psi$ on interior angular intervals.

Apply Lemma~\ref{lem:sourceside} to these normalized annuli, with $\eps=a/b^{q+1}$. Its hypotheses now follow from the interior profile; neither the lemma nor this normalization uses the growth law~\eqref{eq:negsourceresonance}. The endpoint integral is consequently cut off at $y\asymp b^{-(q+1)/(q-1)}$ when $q>1$, and
\begin{equation*}
 J(t)\asymp
 \begin{cases}b^{-q},&q<2,\\b^{-2}\log b,&q=2,\\b^{-2/(q-1)},&q>2.\end{cases}
\end{equation*}
For $q<2$ the finite integral is $\int\psi^{1-q}$; at equality its divergence is logarithmic, and above equality the stated cutoff supplies the power. Each displayed function of $b$ is eventually decreasing. Monotonicity of $b$ bounds the future convolution above by a constant times its value at $t$, while integration over $t<s<t+1$ and~\eqref{eq:negativeslow} give the corresponding lower bound. Thus $b'$ has the same order as the displayed function of $b$. Integration yields~\eqref{eq:negsourceresonance}; for $q=2$ use $t\asymp b^3/\log b$. Already $b\to\infty$ and $b'\le C$ imply $\eps(r)\le C$ and $r|\log\eps(r)|\to0$. None of these assertions is an all-path relative approximation by $\psi$.

\begin{proof}[Proof of Theorem~\ref{thm:jointsource}]
Normalize an annulus by $V(z)=U(rz)/M(r)$. Its equation is $-\Delta V=\eps A_rV^{-q}$ with $A_r=|F'(rz)|^2$. The interior profile supplies the interior bounds required in Lemma~\ref{lem:sourceside}; thus that lemma proves~\eqref{eq:jointvalue}. The independent ball comparison~\eqref{eq:sourceball} gives $\eps A_rV^{-q-1}\le Cy^{-2}$. Interior Harnack and gradient estimates on balls of radius comparable to $y$ therefore give $|\nabla V|\le CV/y$, proving~\eqref{eq:jointgrad}. These estimates also give, on a ball of radius $cy$,
\begin{equation*}
 |f|+y^\eta[f]_{C^{0,\eta}}\le C\eps S_\eps(y)^{-q},\qquad f=\eps A_rV^{-q},
\end{equation*}
for a fixed $0<\eta<1$. The scaled analytic coefficient has the required local bounds by the Cauchy estimate for $F'$.

Use a half-disc contained in the patch where the side bounds hold, and write $V=h+P$ as its harmonic lifting plus its zero-Dirichlet potential. The lifting has bounded Hessian on the smaller patch by odd reflection. In these coordinates $y$ is the distance to the diameter, comparable to the scaled physical distance. The functions $S_\eps$ and $T_\eps$ are comparable at distances differing by fixed factors, uniformly for the present $\eps$ range.

Let $y_1<y_0$ be the height of the chosen half-disc. Suppose $0\le f(z)\le g(z_2)$ for a decreasing envelope $g$, with the local H\"older bound above. The half-plane kernel and the harmonic cap correction give
\begin{equation}\label{eq:weightedGreen}
 |D^2P|\le C\left(g(y)+y^{-2}\int_0^y t g(t)\dd t
          +\int_y^{y_1} g(t)\frac{\dd t}{t}+\int_0^{y_1} t g(t)\dd t\right)\,.
\end{equation}
For poles in a smaller concentric half-disc, the harmonic difference of the two kernels has cap data $O(z_2)$ and vanishes on the diameter. Its second derivatives on a still smaller patch are $O(z_2)$ by odd reflection. For poles outside that concentric half-disc, the Dirichlet kernel itself is harmonic on the observation patch and bounded there by the half-plane kernel, which is $O(z_2)$ by separation. This gives the same derivative bound. For the half-plane part, subtract $f(x)$ at the diagonal and use the local H\"older estimate to obtain $Cg(y)$, including the local term from differentiating the logarithmic potential. Off the diagonal, reflection in the half-plane gives $|D_x^2G(x,z)|\le Cz_2|x-z|^{-3}$. Horizontal integration outside a ball of radius comparable to $y$ bounds the integral of $|x-z|^{-3}$ by $Cy^{-2}$ for $z_2\le y$ and $Cz_2^{-2}$ for $z_2\ge y$. This gives the two middle terms in~\eqref{eq:weightedGreen}. Compact source truncations justify passage to the limit since the first moment $\int_0^{y_1}t g(t)\dd t$ is finite.

For $g=\eps S_\eps^{-q}$, monotonicity gives $g(y)\le2T_\eps(y)$. For $t\ge y$, the three formulas for $S_\eps$ give $g(t)\le Cg(y)(y/t)^\eta$ for some $\eta>0$ depending on $q$. For $q=1$ one can take $\eta=1/2$ by comparing the logarithms. Thus the third term is bounded by $Cg(y)$. Split the last term at $y$ and use monotonicity on $t\ge y$ to bound it by $CT_\eps(y)$. Direct integration gives~\eqref{eq:Texplicit}--\eqref{eq:Tlarge}.

The physical transformation requires one additional term, which we estimate explicitly. Choose an analytic branch $L_0=\log F'$. On a straight side $\operatorname{Im}L_0$ is the turning angle of the image boundary, and its tangential derivative is the curvature times $|F'|$, hence bounded. Subtract the turning angle at the center of a side interval of length comparable to $r$ before multiplying its trace by a cutoff on that interval. The resulting trace is $O(r)$ and uniformly Lipschitz. Its Poisson extension has derivative at distance $ry$ bounded by $C(1+|\log y|)$, by integrating the derivative kernel against trace differences. The remaining harmonic part has zero values on the smaller side interval and amplitude $O(r^{\eta_0})$ by~\eqref{eq:conformal}, taking $\eta_0\le1$. Odd reflection and the harmonic gradient estimate bound its derivative by $O(r^{\eta_0-1})$ on a smaller piece. Therefore
\begin{equation}\label{eq:Fsecond}
 r|F''(rz)|\le C\big(r^{\eta_0}+r(1+|\log y|)\big)=:\mathcal E_r(y)\,.
\end{equation}
The chain rule gives the temporary physical bound
\begin{equation*}
 |D_x^2u|\le C\frac{M(r)}{r^2}
 \left(1+T_\eps(y)+\mathcal E_r(y)S_\eps(y)/y\right)\,.
\end{equation*}
Directly from~\eqref{eq:Ssource}--\eqref{eq:Tlarge},
\begin{equation*}
 \sup_{0<y<y_0}\frac{\mathcal E_r(y)S_\eps(y)/y}{1+T_\eps(y)}
 \le C\big(r^{\eta_0}+r(1+|\log\eps|)\big)\,.
\end{equation*}
For $q>1$, split at $y=\eps^{1/(q-1)}$. In the inner region $(S_\eps/y)/T_\eps\le Cy$. In the outer region $S_\eps/y\le C$ and $|\log y|\le C(1+|\log\eps|)$. For $q<1$ the ratio $(S_\eps/y)/(1+T_\eps)$ is bounded by $C\min\{1,y^q/\eps\}$; splitting at $y=\eps^{1/q}$ bounds its product with $1+|\log y|$ by $C(1+|\log\eps|)$. For $q=1$ write $\ell_y=\log(e/y)$. Where $\eps\ell_y\le1$ and $y\ge\eps$, one has $S_\eps/y\le\sqrt2$ and $\ell_y\le1+|\log\eps|$. Where $\eps\ell_y\le1$ and $y<\eps$, one has $T_\eps\asymp\eps/y$, and $(y/\eps)\ell_y\le C(1+|\log\eps|)$. In the deeper region $\eps\ell_y>1$, the ratio is bounded by $Cy\ell_y$, whose product with $\ell_y$ is uniformly bounded. All these comparisons also hold when $\eps$ stays bounded away from zero. The right side tends to zero for the scales already proved. This proves~\eqref{eq:jointhess}; it permits $F''$ to be unbounded along very close side approaches.
\end{proof}

\begin{remark}
For $c>0$ and $0<q<1$, applying~\eqref{eq:jointvalue} to $u_0$ with its amplitude $M_0(r)$ gives $u_0\le C M_0(r)y$, while $K\asymp r^{-\beta}y$. Hence $(u_c-cK)/K\le C M_0(r)r^\beta\to0$ along all paths. For $q\ge1$, at each sufficiently small fixed $r$ the quotient $S_{\eps(r)}(y)/y$ diverges as $y\downarrow0$. The same two-sided estimates give approaches to $O$ on which $u_c/K$ is unbounded.

The scale $T_\eps$ retains a first moment of the source. For $q\ge2$ it is larger in the outer region than the local source value. The absolute kernel bound~\eqref{eq:weightedGreen} includes this first moment in addition to the local source value. Determining a lossless outer matching formula in that range, and an exact full second-order matching formula in the exponentially thin $q=1$ layer, requires further analysis.
\end{remark}

\section{Sublinear sources and absorption}\label{sec:sublinear}
In this section $a>0$ and $0<p<1$. For the source sign the zero-trace solution is strictly positive. For the absorption sign a nonnegative solution may have an interior zero set; this distinction is necessary when comparing with Definition~\ref{def:natural}.

\subsection{Existence and uniqueness}
The source and absorption signs require different constructions. For the source we combine a zero-trace variational solution with Riesz decomposition and comparison; for absorption we use truncated boundary-value problems, with strict positivity or an interior zero set depending on the size of $c$.

\begin{proof}[Proof of Theorem~\ref{thm:subclass}]
\begin{citationrevision}
For the bounded zero-trace source solution we use the variational construction of~\cite[Section~4]{BO}. Its minimization argument carries over to the present bounded Lipschitz domain: the embedding $H^1_0(\Om)\hookrightarrow L^{p+1}(\Om)$ is compact, $p+1<2$ gives coercivity, and a small positive ray gives a negative value. Elliptic bootstrap and the local boundary barriers, including a sector barrier at $O$, make the resulting positive solution bounded and continuous with zero values on the whole boundary. The singular-trace classification still requires the arguments below.
\end{citationrevision}

For any positive source solution the Riesz decomposition gives $u=cK+a\G[u^p]$, $c\ge0$. The harmonic part is supported at $O$, since it is bounded above by the continuous zero-boundary function $u$ on every ordinary piece. Let $\phi_1$ be the positive first Dirichlet eigenfunction, with eigenvalue $\Lambda_1$. Comparison of $G(x_*,\cdot)$ with $\phi_1$ away from $x_*$ gives $\int u^p\phi_1<\infty$. Also $\int K\phi_1<\infty$, including at the corner. Tonelli's theorem and $\G[\phi_1]=\phi_1/\Lambda_1$ yield $\int u\phi_1<\infty$.

On a fixed sector annulus a nonnegative zero-side subsolution of $-\Delta v\le C(1+v)$ satisfies
\begin{equation}\label{eq:weightedMoser}
 \sup_{\cA}v\le C\left(1+\int_{\cA'}v\psi\right)\,.
\end{equation}
To check the boundary weight, extend by zero across a side; the extra measure in $-\Delta v$ is nonpositive, so the usual local Moser estimate applies. If $S(t)$ is the supremum on a larger annulus and $W$ the weighted integral, split the unweighted integral where $\psi<\eps$ to get $\int v\le\eps^{-1}W+C\eps S(t)$. Choose $\eps$ as a small power of the distance between nested annuli and iterate the resulting inequality $S(s)\le\delta S(t)+C_\delta(t-s)^{-N}(1+W)$. This proves~\eqref{eq:weightedMoser}. Scaling and $\phi_1\circ F\asymp r^\beta\psi$ give the crude bound $u\circ F\le Cr^{-\beta-2}$.

Consequently $u^p\le Cr^{-p(\beta+2)}$. The sector Green formula, or a barrier $r^\sigma f$ with $-f''-\sigma^2f=1$, gives
\begin{equation}\label{eq:potentialimprove}
 \G[u^p]\circ F\le Cr^\sigma\psi,\qquad
 -\beta<\sigma<\min\{\beta,2-p(\beta+2)\}\,.
\end{equation}
The interval is nonempty because $p<1$. As for Lemma~\ref{lem:negativepotential}, compare truncated potentials and then take their monotone limit. This proves the uniform relative trace. If $c=0$, the first improvement supplies a bound $u\le Cr^{-M}$ with $0<M<\beta$ unless it already gives boundedness. Repeat the estimate: when $pM\ge2$ replace this exponent by $pM-1$, which is positive and remains smaller than $\beta$. It decreases by at least one until $pM<2$, when a positive $\sigma$ is available. Thus the zero-trace solution is bounded and $u=a\G[u^p]\in H^1_0$. \citationadded{Once boundedness and $H^1_0$ membership are known, uniqueness is the Brezis--Oswald quotient argument~\cite[Section~2, pp.~56--57]{BO}. Its positive-denominator Picone tests are valid on this Lipschitz domain by approximation in $H^1_0$ and passage to the limit; they do not require a bounded quotient at the corner.}

For $c>0$ put $H_p=\G[K^p]$, $M_p=\sup H_p/K<\infty$. Choose $B$ so large that $B\ge ac^p+aM_p^pB^p$. The functions $cK$ and $cK+BH_p+u_0$ are a subsolution and a supersolution, using $(x+y+z)^p\le x^p+y^p+z^p$. Exhaustion gives existence. For uniqueness at positive trace let $T=\sup u/v$. It is finite by the uniform corner ratio and Hopf's lemma on the ordinary boundary. If $T>1$, then $w=Tv-u\ge0$ and
\begin{equation*}
 -\Delta w\ge a(T-T^p)v^p>0\,.
\end{equation*}
The strong maximum principle, the Hopf boundary lemma, and $w/v\to T-1$ at $O$ give $\inf w/v>0$, contradicting the definition of $T$. The same argument proves order in $c$.

For absorption construct $v_c$ by solving $-\Delta v+av^p=0$ on truncated domains with boundary values $cK$ on the inner arc and zero elsewhere; comparison gives $0\le v_c\le cK$. Its source is dominated by $c^pK^p$, whose Green potential is finite and has vanishing ratio to $K$. Passage to the limit gives~\eqref{eq:subabs}. For uniqueness compare one solution with $v_c+\eps K$, a supersolution for absorption. The uniform trace makes the comparison valid on sufficiently small inner arcs, and monotonicity of $av^p$ gives it on the rest of the truncated domain. Let the arc shrink and then let $\eps\downarrow0$; interchange the two solutions.

For an arbitrary positive natural solution, $u+1$ is positive $(-\Delta+a)$-superharmonic with Riesz source $a(u+1-u^p)\ge0$. The Green kernel of $-\Delta+a$ with pole $x_*$ is comparable to $r^\beta\psi$ near $O$. For the lower bound, in conformal coordinates put $H=r^\beta\psi$ and use a small multiple of $H(1+B r^2)$. It is a subsolution of $-\Delta+a|F'|^2$ for large fixed $B$ and sufficiently small $r$, because $\Delta(r^{\beta+2}\psi)=4(\beta+1)r^\beta\psi$. Hopf's lemma allows the multiple to be chosen below the kernel on the outer arc, and both functions vanish at the vertex and sides. The upper bound follows by comparison with the Dirichlet Green kernel of $-\Delta$. Ordinary boundary comparison and local positivity complete the comparison with $G(x_*,\cdot)$ away from $x_*$. Finiteness of the shifted Green potential at $x_*$ and $u+1-u^p\ge u/2-C_p$ now give $\int_\Om u\phi_1<\infty$. Since $u^p\le u+1$, the ordinary Green potential $\G[u^p]$ is finite. The positive harmonic majorant $u+a\G[u^p]$ has zero ordinary trace by local bounded-source estimates and separated Green comparison. Thus it equals $cK$. Positivity forces $c>0$, and $u\le cK$ yields $\G[u^p]/K\le c^p H_p/K\to0$ uniformly by Lemma~\ref{lem:positivepotential}. Every natural positive solution therefore has the trace required by the comparison argument above.

For large $c$, $cK-ac^pH_p\ge(c-ac^pM_p)K>0$ is a subsolution, proving strict positivity. For small $c$ choose an interior ball $B_R(x_0)$, $r_0<R$, and $\gamma_p=2/(1-p)$. A sufficiently small constant $A>0$ makes $A(|x-x_0|-r_0)_+^{\gamma_p}$ a supersolution of $-\Delta v+av^p=0$ in this ball, since $\gamma_p-2=p\gamma_p$. Its boundary value dominates $cK$ for small $c$, while it is zero in $B_{r_0}$. Comparison yields an actual interior zero set. These solutions are nonnegative, but they are not natural positive solutions. The large-$c$ family supplies the asserted nonuniqueness in the original class.
\end{proof}

\subsection{Sharp second-order estimates for the finite-trace remainder}
Theorem~\ref{thm:remainder} applies to both signs and also to the finite-trace subcritical superlinear solutions considered later. It separates the Green remainder from the stronger Martin singularity of $u$.

\begin{proof}[Proof of Theorem~\ref{thm:remainder}]
In the sector coordinates,
\begin{equation}\label{eq:powerforcing}
 -\Delta(w\circ F)=Q r^{-\beta p}\psi^p(1+o(1))\,.
\end{equation}
We describe the Green expansion rather than differentiate a function-level error. Let $f_n(r)$ be the sine coefficients of the right side. The coefficient $w_n$ satisfies
\begin{equation*}
 -w_n''-r^{-1}w_n'+(n\beta)^2r^{-2}w_n=f_n(r)\,.
\end{equation*}
The condition $w/K\to0$ excludes each negative-power homogeneous term. Its regular Green solution is
\begin{equation*}
 \frac1{2n\beta}\left(r^{-n\beta}\int_0^r s^{n\beta+1}f_n(s)\dd s
       +r^{n\beta}\int_r^R s^{1-n\beta}f_n(s)\dd s\right)
\end{equation*}
plus a regular zero-side outer harmonic correction. At $n=1$ the second integral converges if $\beta<\sigma$, gives the power $r^\sigma$ if $\sigma<\beta$, and gives $r^\beta\log(1/r)$ at equality. Positivity of the Green kernel gives $B>0$. Since $\sigma> -\beta$, the forced angular inverse is positive when $\sigma<\beta$. Projection on the first mode gives~\eqref{eq:Dres}. The bound $|f_n(s)|\le C s^{\sigma-2}$ is uniform in $n$. For $n\beta>|\sigma|+1$ the displayed radial integrals therefore give
\[
 |w_n^{\mathrm{part}}(r)|\le
 \frac{C r^\sigma}{2n\beta}\left(\frac1{n\beta+\sigma}
                        +\frac1{n\beta-\sigma}\right).
\]
These bounds are summable in $n$. The regular harmonic correction is an absolutely convergent geometric series on $r<R/2$, since its coefficients at $R$ are bounded. The finitely many remaining modes are read from their radial integrals. A higher resonance $\sigma=n\beta$, $n\ge2$, gives a term $r^{n\beta}\log(1/r)$ of strictly higher order than $r^\beta$. The series therefore gives the leading profiles uniformly up to a side in absolute value. Scaled zero-Dirichlet estimates on smaller annuli give their weighted first-order versions, using zero side values to divide the error by $\psi$.

To quantify the forcing error, choose a strict exponent $s_0$ between $-\beta$ and $\min\{\beta,\sigma\}$. A sector barrier of degree $s_0$ gives $w/K=O(r^{s_0+\beta})$ by the same compact source truncation and outer-arc comparison used in Lemma~\ref{lem:positivepotential}. Substitution of $u=cK+\epsilon w$ and the expansions~\eqref{eq:conformal}--\eqref{eq:Martin} then replaces $o(1)$ in~\eqref{eq:powerforcing} by $O(r^{\eta_f})$ for some $\eta_f>0$. The same summable coefficient bound applies to this forcing error with $\sigma$ replaced by $\sigma+\eta$, where $0<\eta<\eta_f$ can be chosen to avoid the discrete resonances. The finitely many low modes give a higher-power remainder after the applicable leading terms are removed, and an $O(r^\beta)$ nonlogarithmic remainder at first-mode resonance. The scaled estimates also control first derivatives of these remainders.

Apply Lemma~\ref{lem:physical} to the normalized function and its Laplacian in each row. The forced normalized source converges to $Q|z|^{-\beta p}\psi^p$; in the harmonic and resonant rows it tends to zero. This proves the full $L^\ell$ Hessian limits and the annular bound
\begin{equation}\label{eq:remainderann}
 \norm{D^2w}_{L^\ell(\Om\cap\{r<|x|<2r\})}
 \le C_\ell r^{\min\{\beta,\sigma\}-2+2/\ell}\mathcal L(r),
\end{equation}
where $\mathcal L=\log(1/r)$ in resonance and $1$ otherwise. H\"older's inequality on each annulus and dyadic summation give sufficiency in~\eqref{eq:remainderW}, including $s=1$ by taking $\ell>1$. The function and gradient terms have weaker integrability restrictions. Necessity follows from $\Delta w=-au^p$, which requires $s\beta p<2$, and, in the harmonic-leading case, from the nonzero Hessian of $r^\beta\psi$, which requires $s(2-\beta)<2$. At equality the annular mass stays bounded below; a logarithm increases it.

The only zero leading harmonic Hessian occurs when $\beta=1$ in the harmonic-leading row. That row requires $\beta<\sigma$, hence $0<p<1$ and $\sigma=2-p\in(1,2)$. Subtract $BS_1$ before applying the coefficient formulas. The next forced degree is $2-p<2$, the next harmonic degree. Hence
\begin{equation*}
 (w-BS_1)\circ F=Qr^{2-p}Y_{2-p}+O_\psi^1(r^{2-p+\eta}),\qquad 0<\eta<\min\{p,\eta_f\}\,.
\end{equation*}
The angular inverse exists, although it is not positive. Apply Lemma~\ref{lem:physical} to $w-BS_1$ with $M_r=r^{2-p}$. Since $S_1$ is harmonic, the normalized forcing is unchanged. The joined boundary is $C^{1,1}$ at the flat angle, so local Dirichlet estimates give $S_1\in W^{2,\ell}$ for every finite $\ell$. This proves~\eqref{eq:remainderW} with $S=2/p$. Its harmonic-coordinate Hessian has negligible tail at that finite exponent. The matrix in the forced term has trace $-Q\psi^p$ and is not zero. In the harmonic rows $|\Hmat_\beta[\psi]|=\sqrt2\beta|\beta-1|$, which is nonzero in precisely the asserted cases. Lemma~\ref{lem:tail} now proves all tail statements, including the resonant logarithm.
\end{proof}

\begin{corollary}\label{cor:rawfinite}
Suppose $a>0$, $p>0$, $\beta(p-1)<2$, and~\eqref{eq:remainderhyp} holds. Set $S_0=2/(\beta+2)<1$. For sufficiently small $\delta>0$, the classical Hessian of $u$ satisfies
\begin{equation}\label{eq:rawtail}
 \lim_{L\to\infty}L^{S_0}\abs{\{|D^2u|>L\}\cap B_\delta(O)}
 =\frac\theta2\big(c\kappa\sqrt2\beta(\beta+1)\big)^{S_0}\,.
\end{equation}
On $\Om\cap B_\delta(O)$ the classical Hessian belongs to $L^s$ precisely for $0<s<S_0$. It belongs to $L^{S_0,\infty}$ and fails to belong to $L^{S_0}$. In particular $u\notin W^{2,1}$ near $O$.
\end{corollary}
\begin{proof}
The scaled source is $O(r^{2-\beta(p-1)})$ and tends to zero on complete annuli. Lemma~\ref{lem:physical} gives the Hessian limit of $c\kappa r^{-\beta}\psi$ with uniform finite-exponent bounds. Its polar matrix has the constant norm $c\kappa\sqrt2\beta(\beta+1)$. Lemma~\ref{lem:tail} applies. These assertions concern $u$, whereas~\eqref{eq:remainderW} concerns $w$.
\end{proof}

\subsection{The bounded source solution}
Put $\gamma_p=2/(1-p)>2$. We prove the first three profiles without using a second-order estimate as an input. The critical case requires an annular comparison that is uniform even when the normalized solution is large.

\begin{lemma}\label{lem:subannular}
Let $V>0$ solve $-\Delta V=b(z)V^p$ on an enlarged sector annulus, with zero values on its straight sides, where $0\le b\le C_0$. Suppose that $V$ has a fixed positive lower bound on each fixed interior compact subset. On a smaller complete annulus, if $z_*$ is a fixed interior point, then
\begin{equation}\label{eq:subannular}
 cV(z_*)\psi(\arg z)\le V(z)\le CV(z_*)\psi(\arg z),\qquad
 |\nabla V(z)|\le CV(z_*)\,.
\end{equation}
The constants are independent of an upper bound for $V$.
\end{lemma}
\begin{proof}
Set $c_V=bV^p/(1+V)$. This coefficient is uniformly bounded because $p<1$, and $V+1$ solves $-\Delta(V+1)=c_V(V+1)$. Interior Harnack chains, together with the assumed lower bound, compare $V$ to $V(z_*)$ on interior compact subsets.

Let $D$ be a small fixed half-disc adjoining a straight side. Write $y$ for distance to its diameter. The maximum principle and the Neumann-series inverse for $L=-\Delta-c_V$ hold on $D$ when $\norm{\G_D c_V}_\infty<1/2$. Let $LZ=c_V$ with zero boundary data and let $h_y$ solve $Lh_y=0$ with boundary data $y$. The Neumann series gives a uniform bound on $Z$ and on the nonnegative correction $h_y-y=L^{-1}(c_V y)$. Both corrections have bounded nonnegative forcing for $-\Delta$ and zero boundary data. Comparing with a sufficiently large multiple of the vertical quadratic $y(2R-y)$, where $R$ is the half-disc radius, yields $0\le Z\le Cy$ and $y\le h_y\le Cy$.

The function $H=V-Z$ satisfies $LH=0$ with nonnegative boundary values, hence $H\ge0$. Let $L\sigma=0$ with boundary value one. The Neumann series gives $1\le\sigma\le2$, and both $H/\sigma$ and $h_y/\sigma$ solve $\operatorname{div}(\sigma^2\nabla f)=0$ and vanish on the diameter. If $H\not\equiv0$, apply the boundary Harnack inequality~\cite[Theorem~1.1]{DeSilvaSavin} on a smaller flat patch. The ratio is bounded by a fixed multiple of its value at an interior reference point, since $1\le\sigma\le2$ makes the coefficient uniformly elliptic. The resulting upper bound is immediate if $H=0$. As $H\le V$, interior Harnack and $h_y\asymp y$ give $H\le C V(z_*)y$. The same bound holds for $Z$ since $V(z_*)\ge c>0$. For the lower bound, compare $V$ with the harmonic lifting of its values on a fixed interior portion of the artificial boundary, where interior Harnack gives values at least $cV(z_*)$. A finite covering proves the value estimates in~\eqref{eq:subannular}. Then $bV^p\le CV(z_*)$ on the smaller annuli, and zero-side and interior $W^{2,\ell}$ estimates with $\ell>2$ give the gradient bound.
\end{proof}

\begin{proposition}\label{prop:subzeroprofiles}
Let $U=u_0\circ F$ and $t=\log(1/r)$. Uniformly relative to the displayed positive angular functions, including approaches to the sides,
\begin{equation}\label{eq:subzeroprofiles}
 U(r,\varphi)\sim
 \begin{cases}
 B r^\beta\psi(\varphi),&\beta<\gamma_p,\quad B>0,\\
 C_p r^\beta t^{1/(1-p)}\psi(\varphi),&\beta=\gamma_p,\\
 r^{\gamma_p}\Phi(\varphi),&\beta>\gamma_p,
 \end{cases}
 \qquad
 C_p=\left(\frac{a(1-p)}{2\beta N}\int_I\psi^{p+1}\right)^{1/(1-p)}\,.
\end{equation}
Here $\Phi$ is the unique positive solution of
\begin{equation}\label{eq:subzeroangular}
 -\Phi''-\gamma_p^2\Phi=a\Phi^p,\qquad \Phi>0\text{ in }I,\quad\Phi=0\text{ on }\partial I\,.
\end{equation}
The normalized functions converge in $C^1$ on complete annuli, and their normalized sources converge there uniformly.
\end{proposition}
\begin{proof}
The coordinate equation is $-\Delta U=mU^p$, where $m=a|F'|^2$ is positive and converges to $a$. All comparisons below take place on a small sector with its outer circular boundary included. We use the sublinear weak comparison principle, not a maximum principle for the difference of two solutions. For positive sub- and supersolutions $w,v$ with ordered boundary data, Picone's identity, tested on $\{w>v\}$ with $w-v^2/w$ and $w^2/v-v$, gives a nonnegative Dirichlet expression bounded above by
\begin{equation*}
 \int_{\{w>v\}}m(w^2-v^2)(w^{p-1}-v^{p-1})\le0\,.
\end{equation*}
The strict decrease of $s^{p-1}$ implies $w\le v$. For the regularized tests set $W=w+\eps$ and $V=v+\eps$ and use $(W^2-V^2)_+/W$ and $(W^2-V^2)_+/V$. They have zero Sobolev trace since the boundary data are ordered. Picone's gradient expression is nonnegative for these tests. The positive reaction term on $\{w>v\}$ is bounded by $mWw^p$, and Fatou's lemma applies to the nonnegative term subtracted. Passing to the limit gives the stated inequality without a bounded quotient assumption. The functions compared below belong to $H^1$ and have ordered traces, so the same argument applies at the vertex.

\emph{The harmonic range $\beta<\gamma_p$.} Since $U$ is bounded, its source is bounded. For every $0<s<\min\{\beta,2\}$ the positive angular solution of $-Y_s''-s^2Y_s=1$ gives the supersolution $Cr^sY_s$; its outer values can be chosen above $U$. Thus $U\le Cr^s\psi$. If this bound holds, then $mU^p\le Cr^{ps}$ and the same barrier argument improves it to every exponent $s'<\min\{\beta,2+ps\}$. The iteration $s\mapsto2+ps$ tends to $\gamma_p$. Because $\beta<\gamma_p$, finitely many strict improvements give an $s<\beta$ with $2+ps>\beta$.

The first dual moment of the source is consequently finite: its radial majorant is $\int_0^R r^{1-\beta+ps}\dd r<\infty$. The sector Green formula then gives $U/r^\beta\to B\psi$, with $B>0$ because the positive outer harmonic lifting already has a positive first coefficient. To check convergence up to the sides, the same formula bounds $U/r^\beta$ on complete annuli, and its scaled source is $O(r^{2-\beta+ps})\to0$. Boundary $W^{2,\ell}$ estimates give $C^1$ convergence on smaller complete annuli. The zero endpoint values and the nonzero slopes of $\psi$ turn this into relative convergence.

\emph{The nonlinear range $\beta>\gamma_p$.} \citationadded{Apply~\cite[Theorem~1, pp.~55--56]{BO} on the smooth interval $I$ to $f(s)=\gamma_p^2s+as^p$. Its quotient is strictly decreasing, its limit at zero is infinite, and its limiting linear operator has first eigenvalue $\beta^2-\gamma_p^2>0$. This gives existence and uniqueness in~\eqref{eq:subzeroangular}; the boundary conclusion of~\cite[Lemma~1, pp.~56--57]{BO} gives $\Phi\asymp\psi$.} Small and large multiples of $r^{\gamma_p}\Phi$ are sub- and supersolutions for the variable coefficient $m$. They can be ordered with $U$ on the outer arc, since ordinary boundary regularity and Hopf's lemma give finite positive endpoint quotients there. The comparison principle therefore gives
\begin{equation}\label{eq:subzerobounds}
 c\Phi\le v(t,\cdot)\le C\Phi,\qquad v=r^{-\gamma_p}U\,.
\end{equation}
The equation for $v$ is
\begin{equation*}
 v_{tt}-2\gamma_pv_t+v_{\varphi\varphi}+\gamma_p^2v+mv^p=0\,.
\end{equation*}
Boundary estimates on unit strips give uniform $C^1$ compactness. Its energy
\begin{equation*}
 E(t)=\int_I\left(\frac12v_t^2-\frac12v_\varphi^2+\frac{\gamma_p^2}2v^2+\frac m{p+1}v^{p+1}\right)\dd\varphi
\end{equation*}
satisfies $E'=2\gamma_p\int_Iv_t^2+(p+1)^{-1}\int_Im_tv^{p+1}$. The Cauchy estimate for the H\"older analytic function $F'$ gives $|m_t|\le Ce^{-\eta_0t}\psi^{\eta_0-1}$. In view of~\eqref{eq:subzerobounds}, the last energy term is absolutely integrable. The energy is bounded, hence $\int^\infty\norm{v_t}_2^2\dd t<\infty$. Every translated limit is stationary. The lower bound in~\eqref{eq:subzerobounds} excludes zero, and uniqueness in~\eqref{eq:subzeroangular} identifies every limit with $\Phi$. We obtain convergence in $C^1$; its zero endpoints again give relative convergence.

\emph{The critical range $\beta=\gamma_p$.} A positive outer harmonic lifting gives $U\ge cr^\beta\psi$. Set $v=r^{-\beta}U$, $b(t)=N^{-1}\int_Iv\psi$, and $J(t)=N^{-1}\int_Imv^p\psi$. Lemma~\ref{lem:subannular} applies to each normalized annulus before any upper estimate for $v$ is known. It gives
\begin{equation}\label{eq:subcriticalcompare}
 cb(t)\psi\le v(t,\varphi)\le Cb(t)\psi,\qquad |b'(t)|\le Cb(t),\qquad
 cb(t)^p\le J(t)\le Cb(t)^p\,.
\end{equation}
Projection of the cylinder equation yields
\begin{equation}\label{eq:subcriticalode}
 -b''+2\beta b'=J,\qquad
 b'(t)=\int_0^\infty e^{-2\beta s}J(t+s)\dd s>0\,.
\end{equation}
For this identity no slow variation is assumed: boundedness of $U$ gives $b=O(e^{\beta t})$ and $J=O(e^{\beta pt})$. Integrating $(e^{-2\beta t}b')'=-e^{-2\beta t}J$ excludes a nonzero $e^{2\beta t}$ term, giving the second formula. Its integral over $0<s<1$ and monotonicity of $b$ give $b'\ge cb^p$, so $b\to\infty$.

We also need an upper differential bound. Choose $K$ large enough that $2\beta-C/K>\beta$, where $C$ is the upper constant for $J$ in~\eqref{eq:subcriticalcompare}. At a point where $b'=Kb^p$,
\begin{equation*}
 (b'-Kb^p)'=\big(2\beta K-J/b^p-pK^2b^{p-1}\big)b^p>0
\end{equation*}
for all sufficiently large $t$. Thus an excursion above $Kb^p$ could not cross back below it. On such an excursion $b''=2\beta b'-J\ge(2\beta-C/K)b'$, contradicting $b=O(e^{\beta t})$. Consequently $cb^p\le b'\le Kb^p$ eventually, and
\begin{equation}\label{eq:subcriticalslow}
 b(t)\asymp t^{1/(1-p)},\qquad b'/b\to0,\qquad b(t+s)/b(t)\to1
\end{equation}
for bounded $s$. Integration of the same differential bounds also gives $b(t+s)/b(t)\le C(1+s)^{1/(1-p)}$ for $s\ge0$ and large $t$.

Now normalize $v(t_j+s,\varphi)$ by $b(t_j)$. On every fixed strip the source coefficient $m b(t_j)^{p-1}$ tends to zero, and~\eqref{eq:subcriticalcompare}--\eqref{eq:subcriticalslow} give uniform boundary $W^{2,\ell}$ bounds. A limit $w$ on the whole cylinder satisfies
\begin{equation*}
 w_{ss}-2\beta w_s+w_{\varphi\varphi}+\beta^2w=0,\qquad
 c\psi\le w\le C\psi,\qquad N^{-1}\int_Iw\psi=1\,.
\end{equation*}
The roots for the $n$th sine mode are $(1\pm n)\beta$. Boundedness on the whole cylinder excludes every mode except the constant first mode. Hence $w=\psi$, and uniqueness of this limit gives $v/b\to\psi$ in $C^1$ on unit strips. In particular $J/b^p\to a\Gamma$, where $\Gamma=N^{-1}\int_I\psi^{p+1}$. The polynomial ratio bound in~\eqref{eq:subcriticalslow} permits dominated convergence in~\eqref{eq:subcriticalode}; thus $b'/b^p\to a\Gamma/(2\beta)$. Integration gives $b\sim C_pt^{1/(1-p)}$, with $C_p$ as in~\eqref{eq:subzeroprofiles}. The same $C^1$ convergence and the zero endpoints give the relative profile. In all three cases the equation now gives convergence of the normalized sources, including at the sides.
\end{proof}

The finite Sobolev range requires only boundedness and the sign of the forcing, rather than the full trichotomy.

\begin{lemma}\label{lem:boundedforcing}
Suppose that $w>0$ is bounded near $O$, is continuous and zero there and on the incident sides, and satisfies $0\le-\Delta w\le C$. Then, for some $\delta>0$,
\begin{equation}\label{eq:boundedforcingW}
 w\in W^{2,s}(\Om\cap B_\delta(O))\quad\Longleftrightarrow\quad s<P_\beta,
 \qquad 1\le s<\infty\,.
\end{equation}
\end{lemma}
\begin{proof}
If $\beta<2$, the first dual moment of the bounded transformed source is finite. The sector Green formula gives $w\circ F=B r^\beta\psi+o(r^\beta)$ with $B>0$, since the outer harmonic lifting is positive. Boundary annular estimates give the same convergence in scaled $C^1$. The normalized sources tend uniformly to zero. If $\beta\ne1$, Lemma~\ref{lem:physical} therefore gives a nonzero harmonic Hessian limit and the upper bound $C_\ell r^{\beta-2+2/\ell}$ on each annulus. Dyadic summation gives membership exactly when $s(2-\beta)<2$, including sufficiency for $s=1$ by taking $\ell>1$. At equality the annular $L^s$ mass is bounded below by a positive constant, so the endpoint fails. If $\beta=1$, the joined boundary is $C^{1,1}$ at $O$, and the bounded-source boundary estimate gives every finite $W^{2,s}$ directly.

If $\beta\ge2$, choose any small $\eps>0$ and solve $-Y''-(2-\eps)^2Y=1$ with zero endpoint values. Since $2-\eps<\beta$, this solution is positive and comparable to $\psi$. The barrier $Cr^{2-\eps}Y$ dominates $w$ on a small sector: its negative Laplacian dominates the bounded source, and its outer values can be chosen above those of $w$. Thus $w\circ F\le C_\eps r^{2-\eps}\psi$. Lemma~\ref{lem:physical} gives $\norm{D^2w}_{L^\ell(\{r<|x|<2r\}\cap\Om)}\le C_{\ell,\eps}r^{-\eps+2/\ell}$. For any finite $s\ge1$, take $\ell>s$ and $\eps<2/s$. H\"older's inequality and dyadic summation give $D^2w\in L^s$; the function and gradient terms are less restrictive. This proves~\eqref{eq:boundedforcingW}.
\end{proof}

Applying Lemma~\ref{lem:boundedforcing} to $u_0$ proves
\begin{equation}\label{eq:subzeroW}
 u_0\in W^{2,s}(\Om\cap B_\delta(O))\quad\Longleftrightarrow\quad s<P_\beta,\qquad 1\le s<\infty\,.
\end{equation}
The leading logarithm in Proposition~\ref{prop:subzeroprofiles} occurs only at degree $\gamma_p>2$; this does not exclude logarithms in higher terms of the harmonic case. No $W^{2,\infty}$ conclusion follows from these finite-exponent statements under a bare $C^2$ arc assumption.

\section{Superlinear absorption}\label{sec:superabs}
Let $p>1$ and use $p_c$ from~\eqref{eq:criticalpower}.

\subsection{Existence and prescribed-trace uniqueness}
We first construct the finite-trace family and then transfer the constant-coefficient cone classification to the curved vertex by coefficient squeezing. The same comparison yields the strong singular limit, prescribed-trace uniqueness, and the removability range.

\begin{proof}[Proof of Proposition~\ref{prop:superabs}]
\begin{citationrevision}
Admissibility of $c\delta_O$ for $1<p<p_c$ follows from~\cite[Theorem~5.16]{MV}, or from the Green moment in Lemma~\ref{lem:positivepotential}. With $g(s)=a|s|^{p-1}s$, the existence theorem~\cite[Theorem~3.10]{MV}, the comparison result~\cite[Lemma~3.5]{MV}, and the representation~\cite[Definition~3.3, (3.8)]{MV} give the unique finite-trace family, increasing in $c$, with $u_c+a\G[u_c^p]=cK$. Interior regularity and $0<u_c\le cK$ put these solutions in the natural class. The all-path statement follows from the additional estimate
\begin{equation*}
 0\le1-\frac{u_c}{cK}\le ac^{p-1}\frac{\G[K^p]}K\longrightarrow0\,,
\end{equation*}
using Lemma~\ref{lem:positivepotential}. The bound~\cite[Proposition~5.1]{MV} makes the increasing limit $u_\infty$ locally finite in $\Om$, and interior estimates give a solution. To retain the zero values, center a ball $B_R(P)$ at an ordinary boundary point with $\overline B_R(P)$ avoiding $O$. The usual radial Keller--Osserman supersolution, a sufficiently large multiple of $(R^2-|x-P|^2)^{-2/(p-1)}$, blows up at the sphere. Comparison on $\Om\cap B_R(P)$ bounds every $u_c$ by this function. For each $c$ one can first compare on a concentric ball sufficiently close to the sphere, using boundedness of $u_c$ on the closed patch. The bound is uniform in $c$ on $B_{R/2}(P)$. Harmonic lifting on a smaller patch, using subharmonicity and the common zero values on the true boundary, now gives $u_c\le C d$ uniformly. Thus $u_\infty$ is natural.
\end{citationrevision}

We give the coefficient comparison explicitly, including the critical endpoint. Write $v=u\circ F$, $m_R=\min_{|z|\le R}a|F'|^2$ and $M_R=\max_{|z|\le R}a|F'|^2$. On each punctured sector solve the constant-$m_R$ absorption equation with the boundary values of $v$. If $w_{R,\eps}$ is this solution, monotonicity and the choice $d_R=(m_R/M_R)^{1/(p-1)}$ give
\begin{equation*}
 d_Rw_{R,\eps}\le v\le w_{R,\eps}\,.
\end{equation*}
Keller--Osserman bounds give a locally convergent subsequence as the puncture shrinks. The inequalities $d_Rw_R\le v\le w_R$ persist, and ensure that $w_R$ has continuous zero values on the sides away from the vertex. Restricting to a smaller cap if needed puts this actual constant-coefficient solution in the setting of~\cite[Theorem~5.9]{MV}. The removable alternative would make $u$ extend by zero at $O$, contrary to the maximum principle on $\Om$. The weak and strong alternatives have distinct radial degrees, so comparison for two radii yields the same alternative. In the strong alternative, $d_R\to1$ and $m_R,M_R\to a$ as $R\downarrow0$. The constant-coefficient asymptotics, whose angular profiles scale by $m_R^{-1/(p-1)}$, then give~\eqref{eq:superabsprofiles} uniformly relative to the positive angular profile. The unique positive coefficient-one angular state exists exactly when $d_p>\beta$ by~\cite[Proposition~5.3]{MV}, with $\lambda_S=\beta^2$ and $\lambda_{2,p}=d_p^2$.

In the weak alternative, the comparison gives $c_1K\le u\le c_2K$ near $O$ for positive constants $c_1,c_2$. The upper bound extends to all of $\Om$ by the ordinary boundary estimate for $u$ and $K$. Lemma~\ref{lem:positivepotential} then gives $\G[u^p]/K\to0$ and continuous zero ordinary boundary values for this potential. Consequently $u+a\G[u^p]$ is positive harmonic with zero ordinary values, and hence equals $cK$. This also proves the all-path relative limit. Since $u_\infty\ge u_c$ for every $c>0$, its coefficient in the weak alternative would have to be at least every $c$. Thus $u_\infty$ is strong. Two strong solutions have ratio tending uniformly to one on inner caps. A multiple by $1+\eps$ is a supersolution, so comparison on truncated domains and then $\eps\downarrow0$ give uniqueness.

For $p\ge p_c$, on a local sector let $b$ be the bounded constant-$m_R$ solution with $v$'s outer data and zero values at the vertex and sides. The function $z=(v-b)_+$ is a subsolution of $-\Delta z+m_Rz^p\le0$ with zero outer and side data, since $(b+z)^p-b^p\ge z^p$. Solve the constant-$m_R$ equation above this subsolution on punctured sectors, with its inner values and zero outer and side values. Keller--Osserman bounds and exhaustion produce an actual constant-coefficient solution $W\ge z$. The ordinary zero data persist, including at the cap corners: the radial blow-up comparison gives a uniform local bound up to each zero boundary patch, and a harmonic lifting then gives continuity and zero values. \citationadded{Extend $W$ by zero to the full cone. For compactly supported tests in the open cone only the smooth interior of the cap can contribute a boundary term. There the outward normal derivative is nonpositive because $W\ge0$ and $W=0$ on the cap. The zero extension is therefore a subsolution of the same absorption equation in the full cone. The dominating-solution construction of~\cite[Theorem~4.3(ii) and the proof of Theorem~5.7, Step~4, pp.~960--961]{MV} gives a full-cone solution above it with zero side values. After scaling the constant $m_R$ to one,~\cite[Theorem~5.6]{MV} excludes a nonzero such solution for $p\ge p_c$, including equality. Thus $W=0$.} Hence $v\le b$ is bounded and removable. The maximum principle excludes a positive zero-boundary solution on the bounded original domain. This uses the cone endpoint theorem rather than a strictly supercritical estimate. All comparisons are local at $O$ and require no simple connectivity.
\end{proof}

\subsection{Full second derivatives}
Finite-trace remainders satisfy Theorem~\ref{thm:remainder}. For the strong solution and every finite $1<\ell<\infty$,
\begin{equation}\label{eq:superabsH}
 r^{d_p+2}D_x^2u_\infty(F(r\cdot))
 \longrightarrow D_z^2(|z|^{-d_p}\Psi_p)\quad\text{in }L^\ell(\cA)\,.
\end{equation}
Indeed the normalized functions and sources converge on the complete annulus and vanish on its true sides; Lemma~\ref{lem:physical} applies. With $S_p=2/(d_p+2)<1$, Lemma~\ref{lem:tail} gives
\begin{equation}\label{eq:superabsTail}
 \lim_{L\to\infty}L^{S_p}\abs{\{|D^2u_\infty|>L\}\cap B_\delta(O)}
 =\frac12\int_I|\Hmat_{-d_p}[\Psi_p]|^{S_p}\dd\varphi>0\,.
\end{equation}
Thus its classical Hessian is in weak $L^{S_p}$ and in $L^s$ exactly for $s<S_p$, and the solution is not in $W^{2,1}$ near the vertex. The radial matrix entry $d_p(d_p+1)\Psi_p$ excludes cancellation, including at the flat angle.

\section{Superlinear sources}\label{sec:supersource}
We consider $-\Delta u=au^p$, $a>0$, $p>1$. A bounded zero-boundary positive solution exists for every finite $p$ in two dimensions. This fact alone does not imply existence of an unbounded solution for every parameter and domain.

\subsection{Existence and the limits of uniqueness}
We first determine the positive mass interval. A direct variational
construction gives a bounded solution, and a mountain-pass argument gives
a second solution at the same sufficiently small positive mass as the
strictly stable least solution.

\begin{proof}[Proof of Theorem~\ref{thm:positive-source-mass}]
The smooth-domain result for scalar multiples of a boundary measure is~\cite[Theorem~1.2]{BVV}. Here the subcritical condition comes from the Green moment at the corner.
Every natural solution is positive and superharmonic. Green exhaustion gives
\begin{equation}\label{eq:positive-source-representation}
 u=cK+a\G[u^p],\qquad c\ge0.
\end{equation}
Indeed the harmonic minorant is between $0$ and $u$, so its Martin measure is supported at $O$. If $c>0$, evaluating~\eqref{eq:positive-source-representation} at $x_*$ shows that $\G[K^p](x_*)<\infty$. The radial moment in the proof of Lemma~\ref{lem:positivepotential} is finite exactly when $\beta(p-1)<2$. This proves necessity, including failure at equality.

In the subcritical range put $H_p=\G[K^p]$ and $M_p=\sup_\Om H_p/K$. By Lemma~\ref{lem:positivepotential}, the map $v\mapsto cK+a\G[v^p]$ preserves $[cK,2cK]$ when $a2^pM_pc^{p-1}\le1$. On this interval its Lipschitz constant in $\sup_\Om|v|/K$ is at most $ap(2c)^{p-1}M_p$, which is less than one for all sufficiently small $c$. Iteration from $cK$ gives the asserted small-mass solution and uniqueness in this interval.

More generally, if $v$ is a solution with mass $c'\ge c>0$, then $v$ bounds the same iteration at mass $c$. Its increasing limit is locally bounded and has continuous zero boundary values off $O$ by comparison with $v$. Monotone convergence in the integral equation and interior elliptic regularity give a natural solution of mass $c$. The mass is $c$ because a Green potential has zero greatest nonnegative harmonic minorant. Iteration gives minimality at each realized mass and yields
\[
 U_{c'}-U_c\ge (c'-c)K
\]
whenever both masses are realized.

To bound the admissible masses, fix a ball $B\Subset\Om$. Let $\phi_B$ be its first positive Dirichlet eigenfunction, normalized by $\int_B\phi_B\,\dd x=1$, and let $\Lambda_B$ be its eigenvalue. For $X=\int_Bu\phi_B\,\dd x$, Green comparison on $B$ and Jensen's inequality give
\[
 c\int_B K\phi_B\,\dd x
 \le X-\frac{a}{\Lambda_B}\int_B u^p\phi_B\,\dd x
 \le X-\frac{a}{\Lambda_B}X^p.
\]
The last function has a finite positive maximum for $X\ge0$. Thus the supremum $c^*$ of the admissible masses is finite and positive, and every $0<c<c^*$ is realized.

Fix $c<c^*$. For a solution $v$ at a strictly larger mass, $w=v-U_c>0$ by~\eqref{eq:positive-source-representation} and minimality. Convexity gives $-\Delta w\ge apU_c^{p-1}w$. Testing with $\zeta^2/w$ and completing the gradient square gives~\eqref{eq:positive-source-stability}. For a nonnegative cutoff $\eta$ supported away from $O$, the functions $U_c\eta$ and $U_c\eta^2$ have zero Sobolev trace. They are admissible in the stability inequality and the equation, respectively, by local Dirichlet regularity at the ordinary boundary and density. Combining the two tests yields
\[
 a(p-1)\int_\Om U_c^{p+1}\eta^2\,\dd x
 \le \int_\Om U_c^2|\nabla\eta|^2\,\dd x.
\]
For a smooth nonnegative cutoff $\chi$ on $\overline\Om$ that vanishes near $O$, put $\sigma=(p+1)/(p-1)$ and take $\eta=\chi^\sigma$. Young's inequality gives
\[
 \int_\Om U_c^{p+1}\chi^{2\sigma}\,\dd x
 \le C\int_\Om|\nabla\chi|^{2\sigma}\,\dd x
\]
with $C$ independent of $c$. Consequently $U_c$ is uniformly bounded in $L^{p+1}$ on every region separated from $O$, including ordinary boundary portions.

Set $s=(p+1)/p>1$. Interior and local zero-Dirichlet $W^{2,s}$ estimates bound $U_c$ uniformly on smaller interior and ordinary-boundary neighborhoods, using the preceding $L^{p+1}$ bound and $aU_c^p\in L^s$ there. The local Dirichlet estimates apply since the ordinary boundary is smooth. The planar Sobolev embedding gives local equicontinuity up to this boundary. Thus the increasing limit $U_* = \lim_{c\uparrow c^*}U_c$ is finite in the interior, continuous and zero on the ordinary boundary. Monotone convergence in~\eqref{eq:positive-source-representation} gives $U_*=c^*K+a\G[U_*^p]$. Interior elliptic regularity makes $U_*$ classical. The Green potential has zero greatest nonnegative harmonic minorant, so the coefficient of $K$ in this representation is exactly the mass $c^*$. Every solution of mass $c^*$ bounds each $U_c$ for $c<c^*$ by iteration, so $U_*=U_{c^*}$. Finally, monotone convergence in~\eqref{eq:positive-source-stability} gives the endpoint inequality.
\end{proof}

\begin{proof}[Proof of Proposition~\ref{prop:supersourceexist}]
The compact embedding $H^1_0(\Om)\hookrightarrow L^{p+1}(\Om)$ in two dimensions gives a minimizer of the Dirichlet quotient with the $L^{p+1}$ norm fixed. Rescaling it gives a positive solution. Its right-hand side belongs to $L^2(\Om)$ by the planar Sobolev embedding. The Dirichlet Green representation is bounded and continuous with zero data on the whole boundary, including $O$. Indeed the planar logarithmic Green bound gives a uniform $L^2$-to-$L^\infty$ estimate, and approximation by bounded sources gives continuity up to the boundary by Dirichlet regularity on a Lipschitz domain. Interior elliptic regularity gives the classical equation.

Suppose $1<p<p_c$. Theorem~\ref{thm:positive-source-mass} gives the
least solution $U_c$ with $cK\le U_c\le2cK$ for small $c$. The potential $apU_c^{p-1}$ has quadratic-form norm $O(c^{p-1})$: its corner degree is $-\beta(p-1)>-2$, so H\"older's inequality and the planar Sobolev embeddings apply with some exponent greater than one. Thus
\begin{equation*}
 ap\int U_c^{p-1}\xi^2\le\eta_c\int|\nabla\xi|^2,\qquad \eta_c\to0\,.
\end{equation*}
This proves strict stability.

To obtain a second solution seek $u=U_c+v$ with $v\in H^1_0(\Om)$ and use
\begin{equation*}
 E_c(v)=\frac12\int|\nabla v|^2
 -\frac a{p+1}\int\big((U_c+v_+)^{p+1}-U_c^{p+1}-(p+1)U_c^pv_+\big)\,.
\end{equation*}
The integrand in the nonlinear remainder is bounded by $C(U_c^{p-1}v^2+|v|^{p+1})$. For $s,t\ge0$ the difference of $(U_c+s)^p-U_c^p$ and $(U_c+t)^p-U_c^p$ is bounded in absolute value by $C(U_c^{p-1}+s^{p-1}+t^{p-1})|s-t|$. The form bound and the planar embedding therefore give a $C^1$ functional, with the derivative of the nonlinear remainder given by the positive-part power difference. The same remainder bound makes $E_c$ positive on a small sphere when $c$ is small. Along a fixed nonzero nonnegative smooth ray, the nonlinear energy grows at least as its $(p+1)$-st power, so $E_c$ is negative sufficiently far out. To bound its Palais--Smale sequences choose $2<\tau<p+1$ and put $f_U(s)=(U+s)^p-U^p$, $F_U(s)=\int_0^sf_U(t)\dd t$. The scalar inequality
\begin{equation*}
 F_U(s)-\tau^{-1}sf_U(s)\le C_{p,\tau}U^{p-1}s^2,\qquad U>0,\ s\ge0,
\end{equation*}
follows by scaling $s/U$: the quotient is bounded near zero and has a negative leading coefficient at infinity. Therefore
\begin{equation*}
 E_c(v)-\tau^{-1}E_c'(v)[v]
 \ge\left(\frac12-\frac1\tau-C_{p,\tau}\eta_c\right)\norm{\nabla v}_2^2\,.
\end{equation*}
For small $c$ the coefficient is positive, giving boundedness. Since $U_c^{p-1}\in L^r(\Om)$ for some $r>1$, compactness of $H^1_0$ in $L^{2r/(r-1)}$ gives compactness in $L^2(U_c^{p-1}\dd x)$. The pure-power term is compact by the embedding into $L^{p+1}$. For a bounded Palais--Smale sequence $v_j\weak v$, the power difference estimate and these two compact embeddings make the difference of the nonlinear derivatives paired with $v_j-v$ tend to zero. Subtracting the weakly convergent Dirichlet part from $E_c'(v_j)[v_j-v]\to0$ gives $v_j\to v$ in $H^1_0$. The mountain-pass theorem~\cite[Theorem~2.1, pp.~354--355]{AR} now gives a nonzero critical point. Testing its negative part gives $v\ge0$. Its equation is $-\Delta v=a((U_c+v)^p-U_c^p)\ge0$, so interior regularity and the strong maximum principle give $v>0$.

At an ordinary boundary point, multiply $v$ by a smooth cutoff equal to one near the point and supported away from $O$. The product is in $H^1_0(\Om)$ and its Laplacian is in $L^2$, because $U_c$ is bounded on the support and $v$ belongs to every finite $L^s$. The Dirichlet Green bound used above gives continuous zero data there. Hence $u=U_c+v$ is a natural solution. The nonnegative superharmonic correction $v$ has a Riesz decomposition whose harmonic part is $dK$, $d\ge0$, since $v$ vanishes continuously at every ordinary boundary point. If $d>0$, then $v\circ F\ge c_1e^{\beta t}\psi$ for large $t$ by~\eqref{eq:Martin}. Angular Poincar\'e and conformal invariance would give infinite Dirichlet energy for $v$, a contradiction. The correction has zero atom, and additivity of the Riesz decomposition gives atom $c$ for $u$. Thus $u$ and $U_c$ are distinct solutions with the same prescribed positive trace.
\end{proof}

\subsection{Local classification of an existing solution}
We prove Theorem~\ref{thm:sourcefirst} in this subsection to specify the input to the logarithmic refinement. In particular it is a local alternative, not a global construction of every alternative.

\begin{proof}[Proof of Theorem~\ref{thm:sourcefirst}]
\citationadded{We adapt the doubling argument of~\cite[Proposition~5.1]{BPV}, using the two-dimensional plane Liouville property and~\cite[Corollary~1.2]{DSS} for its bounded half-plane limit, to obtain}
\begin{equation}\label{eq:sourceuniversal}
 u(F(re^{i\varphi}))\le Cr^{-d_p}\psi(\varphi)\,.
\end{equation}
To obtain the estimate without $\psi$, put $U=u\circ F$ in the straight sector. A violating sequence $z_j$, $r_j=|z_j|\to0$, satisfies $r_j U(z_j)^{(p-1)/2}\to\infty$. Apply doubling to $U^{(p-1)/2}$ on the closed annulus $r_j/2\le |z|\le2r_j$, excluding the artificial circular ends and using the zero side data. It gives points $w_j$ with $U(w_j)\ge U(z_j)$ and radii $L_j U(w_j)^{-(p-1)/2}$, where $L_j\to\infty$, on which the scaled functions are uniformly bounded and which miss the circular ends. The scale $U(w_j)^{-(p-1)/2}$ is $o(r_j)$, so at most one straight side can remain at finite scaled distance. On fixed scaled balls the coefficient $a|F'|^2$ tends uniformly to $a$. The rescaled functions have value one at their centers. Flat-side boundary estimates, using the bounded forcing, keep those centers a positive scaled distance from any limiting side. Interior and flat-side compactness give a bounded nonzero nonnegative limit on the whole plane or on a half-plane, with zero data on its boundary in the latter case. The plane case contradicts the Liouville property of positive superharmonic functions in dimension two. The half-plane case contradicts~\cite[Corollary~1.2]{DSS}, whose bounded-on-strips hypothesis holds for the bounded limit. Finally $r^{d_p}U(r\cdot)$ and its Laplacian are bounded on fixed enlarged annuli. The flat-side gradient bound and zero side traces then supply the factor $\psi$ in~\eqref{eq:sourceuniversal}.

Set $t=\log(1/r)$ and $v=r^{d_p}u\circ F$. Then
\begin{equation}\label{eq:sourcecylinder}
 v_{tt}+2d_pv_t+v_{\varphi\varphi}+d_p^2v+m(t,\varphi)v^p=0,\qquad m=a|F'|^2\,.
\end{equation}
The bound~\eqref{eq:sourceuniversal} and flat-side strip estimates give uniform $C^{1,\eta}$ compactness on translated bounded strips. A disk centered at $z=re^{i\varphi}$ of radius $c r\psi(\varphi)$ lies in the straight sector, for fixed sufficiently small $c>0$. The bound on $F'-1$ in~\eqref{eq:conformal} and Cauchy's estimate on this disk give $|m_t|\le Ce^{-\eta_0t}/\psi(\varphi)$. In the energy identity this is multiplied by $v^{p+1}\le C\psi^{p+1}$. The energy
\begin{equation*}
 E(t)=\int_I\left(\frac12v_t^2-\frac12v_\varphi^2+\frac{d_p^2}2v^2+\frac m{p+1}v^{p+1}\right)\dd\varphi
\end{equation*}
satisfies $E'=-2d_p\int v_t^2+(p+1)^{-1}\int m_tv^{p+1}$. Integration by parts on bounded strips is justified by flat-side $W^{2,\ell}$ regularity and the zero side trace of $v_t$. The coefficient derivative is integrable against $v^{p+1}$. Since $E$ is bounded, integration gives $\int_T^\infty\int_I v_t^2<\infty$, and $E$ has a limit. The dissipation on every bounded translated strip tends to zero, so all translated limits are stationary. Their angular equation is~\eqref{eq:sourceangular}. Testing against $\psi$ excludes a positive angular state when $d_p\ge\beta$. \citationadded{When $d_p<\beta$, existence is~\cite[Lemma~9]{HMR}, applied on $I$ with first eigenvalue $\beta^2$ and linear parameter $d_p^2$, followed by the coefficient scaling $a^{-1/(p-1)}$. For uniqueness put $s=\beta\varphi$ and $w(s)=(a/\beta^2)^{1/(p-1)}\Phi_p(s/\beta)$. Then $-w''=(d_p/\beta)^2w+w^p$ on $(0,\pi)$ with zero endpoints, so~\cite[Theorem~3.1]{BPV} applies.} The zero state has energy zero. A positive stationary state has energy $-a(p-1)\int_I\Phi_p^{p+1}/(2(p+1))<0$, by its angular equation. The single limiting energy therefore selects the same state for all translated limits. In the positive case, $C^1$ convergence and the nonzero endpoint slopes give the relative strong profile.

When $v\to0$, let $b=N^{-1}\int v\psi$ and $J=N^{-1}\int mv^p\psi$. On a translated strip about time $\tau$, the function $e^{d_p(t-\tau)}v$ solves $(-\Delta_{t,\varphi}-mv^{p-1})w=0$ with uniformly bounded nonnegative potential. The small-half-disc change of unknown used in Lemma~\ref{lem:subannular} gives the required boundary comparison. In detail, on a fixed small half-disc at a flat side solve the same linear equation for $\sigma$ with boundary data one and for $h$ with fixed nonzero nonnegative cap data supported away from the flat side. The small-domain Green inverse gives $1\le\sigma\le2$ and a uniform bound for $h$. The harmonic lifting gives a lower bound for $h$ by a multiple of the distance to the flat side on the smaller patch, and the flat-side gradient bound gives its upper bound by that distance. Both $w/\sigma$ and $h/\sigma$ solve the same divergence-form equation with coefficient $\sigma^2$. The boundary Harnack inequality~\cite[Theorem~1.1]{DeSilvaSavin} compares them on the smaller flat patch after normalization at an interior point. Interior Harnack chains link these normalizations throughout the section and on neighboring time sections. Flat-side and interior gradient estimates then give
\begin{equation*}
 cb(t)\psi\le v(t,\varphi)\le Cb(t)\psi,\qquad |b'|\le Cb\,.
\end{equation*}
Projection gives $b''+2d_pb'+(d_p^2-\beta^2)b+J=0$, with $0<J\le Cb^p$. Thus $z=-b'/b$ is bounded and $z'=(z-(d_p-\beta))(z-(d_p+\beta))+o(1)$. Away from its two roots the derivative has a fixed sign at large times; boundedness prevents repeated passages in both directions. Hence $z$ tends to one root. Translate and divide by $b(t_j)$. Whole-cylinder positive limits are bounded above and below by constants times the first exponential mode. Every higher sine mode has a different exponential in at least one direction and is excluded. Consequently $v/b\to\psi$ in $C^1$ and
\begin{equation*}
 J/b^p\longrightarrow\Gamma_a:=aN^{-1}\int_I\psi^{p+1}\,.
\end{equation*}
This also gives relative convergence up to the endpoints.

Write $h_-=d_p-\beta$ and $h_+=d_p+\beta$. If the selected root $h$ is positive, $J=O(b^p)$ and $-b'/b\to h$ imply $J(t)\lesssim e^{-p(h-\varepsilon)t}$ for every sufficiently small $\varepsilon>0$. The quantity $e^{h_-t}(b'+h_+b)$ has derivative $-e^{h_-t}J$. At the slow root $h=h_->0$, it tends to a limit, which is nonnegative because $(b'+h_+b)/b\to h_+-h_->0$. If the limit were zero, the tail integral of its derivative and then the equation for $b'+h_+b$ would give $b(t)=O((1+t)e^{-\min\{h_+,p(h_--\varepsilon)\}t})$, contrary to $-b'/b\to h_-$. A positive limit gives $b\sim Ce^{-h_-t}$ with $C>0$ by the same first-order equation. At the fast root the quantity has limit zero, and hence
\begin{equation*}
 b'+h_+b=e^{-h_-t}\int_t^\infty e^{h_-s}J(s)\dd s>0.
\end{equation*}
Thus $e^{h_+t}b$ is increasing and has integrable derivative by the bound on $J$ with $p(h_+-\varepsilon)>h_+$. Its limit is finite and strictly positive. If $d_p<\beta$ the negative root is incompatible with $b\to0$, leaving only the regular branch. If $d_p>\beta$ the two positive roots give the regular and weak singular branches. At $d_p=\beta$, the remaining root zero has $b'/b\to0$. The stable convolution in $b''+2\beta b'=-J$ then gives $b'/b^p\to-\Gamma_a/(2\beta)$, and integration yields~\eqref{eq:criticalfirst}. Dominated convergence in the convolution is justified by $|b'/b|<\eps$ with $p\eps<2\beta$.

In the regular case the full angular $C^1$ rescaling at the fast root gives $|\partial_t(u\circ F)|+|\partial_\varphi(u\circ F)|\lesssim e^{-\beta t}$. Its square is integrable on the strip, giving finite local Dirichlet energy by conformal invariance. In each singular case the angular derivative of the limiting profile is nonzero on some closed interior angular interval. The same $C^1$ rescaling bounds the strip energy there below by a positive multiple of the squared amplitude of $u\circ F$, whose integral in $t$ diverges. Away from $O$, the solution has finite energy by ordinary boundary regularity. For the Martin mass, use the Riesz representation~\eqref{eq:positive-source-representation}. In the regular, critical singular, and supercritical singular alternatives, $u/K\to0$. Since $u\ge cK$, these alternatives have $c=0$.

In the subcritical weak alternative the relative profile bounds $u/K$ near $O$. Away from $O$, the ordinary boundary estimate for the bounded source gives $u\le Cd$, while $K\asymp d$. Hence $u\le CK$ throughout $\Om$. Lemma~\ref{lem:positivepotential} gives
\begin{equation*}
  0\le\frac{u}{K}-c=a\frac{\G[u^p]}{K}
    \le aC^p\frac{\G[K^p]}{K}\longrightarrow0\qquad(x\to O)\,.
\end{equation*}
Thus the weak-profile coefficient is precisely $c\kappa$ with $c>0$. \citationadded{For a genuinely smooth marked boundary point, where $\theta=\pi$, the critical and strong alternatives agree with~\cite[Theorems~1.2--1.3]{BPV}. Those smooth-boundary statements do not replace the curved-angle argument above.}
\end{proof}

\subsection{Stability of the supercritical singular profile}
\begin{proof}[Proof of Proposition~\ref{prop:source-angular-stability}]
For $\lambda=-a$ and $k=p$, pull back $Q_u$ to the conformal chart. With $r=e^{-t}$, the form on a tail is
\[
 \int\left(|\xi_t|^2+|\xi_\varphi|^2
 -ap|F'(re^{i\varphi})|^2 r^2(u\circ F)^{p-1}\xi^2\right)\dd t\dd\varphi.
\]
Since $d_p(p-1)=2$, the potential tends uniformly to $ap\Phi_p^{p-1}$ by
Theorem~\ref{thm:sourcefirst} and the corner chart estimates. For $\nu_1>0$
the limiting angular form controls a positive multiple of
$\int_I|\xi_\varphi|^2$. This follows by combining the positive $L^2$
bound from $\nu_1$ with the derivative bound obtained from boundedness of
the angular potential. Uniform convergence then gives coercivity on a tail.
Away from $O$ the potential is bounded. A smooth partition with squares
summing to one therefore gives
$Q_u[\xi]\ge b\norm{\nabla\xi}_2^2-C\norm{\xi}_2^2$ for some $b>0$.
The compact Dirichlet embedding gives finite Morse index. For $\nu_1<0$
choose a compactly supported smooth angular test with negative angular
form. A time cutoff on a sufficiently long interval makes its time
energy smaller than the angular gain. Disjoint translations of this test
along the tail yield arbitrarily many independent negative directions.

Let $A=\max_I\Phi_p$. The first integral of~\eqref{eq:sourceangular}
and symmetry about the maximum give
\[
 \frac{\theta}{2}=\int_0^1
 \left[d_p^2(1-z^2)+\frac{2a}{p+1}A^{p-1}(1-z^{p+1})\right]^{-1/2}\dd z.
\]
As $p\downarrow p_c$ we have $d_p\uparrow\beta$. If $A^{p-1}$ stayed
bounded below along a sequence, the limit superior of the integral
would be strictly smaller than $\pi/(2\beta)=\theta/2$, by comparison
with the integrable first-term bound. Hence $A^{p-1}\to0$, and Poincar\'e's
inequality gives $\nu_1>0$ for $p$ close enough to $p_c$. Finally, testing
the angular form with $\Phi_p$ and using its equation yields
\[
 \langle\mathcal A_p\Phi_p,\Phi_p\rangle
 =d_p^2\int_I\Phi_p^2-a(p-1)\int_I\Phi_p^{p+1}
 \le \bigl(p d_p^2-(p-1)\beta^2\bigr)\int_I\Phi_p^2.
\]
This proves the sufficient condition for a negative eigenvalue.
\end{proof}

\subsection{Noncritical second derivatives}
For the weak branch, Theorem~\ref{thm:remainder} and Corollary~\ref{cor:rawfinite} apply. For the strong source branch, the analogue of~\eqref{eq:superabsH} is
\begin{equation}\label{eq:sourcestrongH}
 r^{d_p+2}D_x^2u(F(r\cdot))\longrightarrow D_z^2(|z|^{-d_p}\Phi_p)
 \quad\text{in }L^\ell(\cA),\quad1<\ell<\infty\,.
\end{equation}
Its exact tail is~\eqref{eq:superabsTail} with $\Phi_p$ in place of $\Psi_p$. For the regular branch the sharp finite Sobolev range is $s<P_\beta$, as in~\eqref{eq:subzeroW}. These statements follow from normalized source convergence and Lemmas~\ref{lem:physical} and~\ref{lem:tail}. None is obtained by differentiating an unquantified relative error.

\subsection{The second logarithm for every critical singular solution}
We now fix $p=p_c$, put $\mu=\beta/2=1/(p-1)$, and work with the critical alternative~\eqref{eq:criticalfirst}. No relation to a particular constructed family is assumed.

We use the projections $P,Q_0$, the operator $H_Q$, and the constants from~\eqref{eq:critconstants1}--\eqref{eq:critconstants3}. One has $D_*>0$ and $\kappa_*<0$. Indeed $Q_0(\psi^p)\ne0$ and $H_Q$ is positive on this subspace.

\begin{lemma}\label{lem:highinverse}
Let $L=\partial_t^2+2\beta\partial_t+\partial_\varphi^2+\beta^2$ on $[T,\infty)\times I$. Suppose $w=Q_0w$ is polynomially bounded, has zero angular traces, and $Lw=g=Q_0g$, with $\norm{g(t,\cdot)}_\infty\le Ct^{-h}$ for a fixed real $h$. Then
\begin{equation}\label{eq:highestimate}
 w=O_\psi^1(t^{-h}),\qquad
 \norm w_{W^{2,\ell}((t-1/2,t+1/2)\times I)}\le C_\ell t^{-h},\quad1<\ell<\infty\,.
\end{equation}
The same assertion holds with a fixed power of $\log t$ multiplying the weight.
\end{lemma}
\begin{proof}
In the normalized sine basis $e_n=(2/\theta)^{1/2}\sin(n\beta\varphi)$, $n\ge2$, the time Green factor with zero value at $T$ and without an exponentially growing mode is
\begin{equation*}
 G_n(t,s)=-\frac{e^{-\beta(t-s)}}{2n\beta}
 \left(e^{-n\beta|t-s|}-e^{-n\beta(t+s-2T)}\right)\,.
\end{equation*}
The absolute angular kernel is bounded, uniformly in both angular variables, by
\begin{equation*}
 C e^{\beta|t-s|}\sum_{n\ge2}\frac{e^{-n\beta|t-s|}}n
 \le C e^{-\beta|t-s|}\big(1+\log_+(1/|t-s|)\big)\,.
\end{equation*}
This time kernel is integrable. Splitting $s$ into $[t/2,3t/2]$ and its complement shows that convolution preserves every fixed polynomial weight, and likewise a logarithmic factor. The extension of $w(T,\cdot)$ has modes $e^{-(n+1)\beta(t-T)}$ and is exponentially small on smaller strips. The other homogeneous root $e^{(n-1)\beta t}$ is excluded by polynomial boundedness. Hence the representation recovers $w$ and gives its unweighted supremum bound. Local boundary $W^{2,\ell}$ estimates give the second assertion. Taking one $\ell>2$ gives bounded unweighted first derivatives. Since $w$ is zero at both angular endpoints, $|w|\le C\psi\norm{w_\varphi}_\infty$. This proves~\eqref{eq:highestimate} without a weighted hypothesis on $g$ or derivatives of $g$.
\end{proof}

\begin{proof}[Proof of Theorem~\ref{thm:critsecond}]
Set $v=a^\mu r^\beta u\circ F$. By Theorem~\ref{thm:sourcefirst},
\begin{equation}\label{eq:critv}
 Lv+(1+e(t,\varphi))v^p=0,\qquad
 e=|F'|^2-1=O(e^{-\eta_0t}),\qquad v\sim At^{-\mu}\psi\,.
\end{equation}
Write $v=b(t)\psi+\xi$, $\xi=Q_0v$. Lemma~\ref{lem:highinverse}, applied to $-Q_0((1+e)v^p)$, gives $\xi=O_\psi^1(t^{-\mu-1})$. In particular $\xi/(b\psi)=O(1/t)$. Projection gives
\begin{equation}\label{eq:criticalb}
 b''+2\beta b'+J=0,\qquad J=N^{-1}\int_I(1+e)v^p\psi
                 =\Gamma_* b^p+O(t^{-\mu-2})\,.
\end{equation}
The stable convolution first gives $b'=O(t^{-\mu-1})$. It then allows $b(s)^p$ to be replaced by $b(t)^p$ in that convolution, with error $O(t^{-\mu-2})$; only the derivative of $b^p$ is used. Thus
\begin{equation*}
 b'=-\frac{\Gamma_*}{2\beta}b^p+O(t^{-\mu-2}),\qquad
 (b^{-1/\mu})'=\frac{\Gamma_*}{2\beta\mu}+O(1/t)\,.
\end{equation*}
Integration and the first profile imply, with $f=b-At^{-\mu}$,
\begin{equation}\label{eq:firstgain}
 f=O(t^{-\mu-1}\log t),\qquad f'=O(t^{-\mu-2}\log t)\,.
\end{equation}
This is a gain derived from the first profile, rather than an additional remainder hypothesis.

Subtract $t^{-\mu-1}\chi$ from the higher-mode equation. The identity $H_Q\chi=A^pQ_0\psi^p$ cancels its leading forcing. Lemma~\ref{lem:highinverse} gives
\begin{equation*}
 \xi=t^{-\mu-1}\chi+O_\psi^1(t^{-\mu-2}\log t)\,.
\end{equation*}
Taylor expansion in the projected nonlinearity is valid uniformly up to the endpoints. Both arguments are comparable to $At^{-\mu}\psi$; even if $1<p<2$ the possibly singular factor $\psi^{p-2}$ is multiplied by an error squared of order $\psi^2$. Using~\eqref{eq:firstgain} yields
\begin{equation}\label{eq:resonantf}
 f''+2\beta f'+\frac{2\beta(\mu+1)}t f
 =-R_*t^{-\mu-2}+O(t^{-\mu-3}(\log t)^2)\,.
\end{equation}
In particular no derivative of the coefficient error $e$ has been taken. Let $h=f'+(\mu+1)f/t$. Then
\begin{equation*}
 h'+\left(2\beta-\frac{\mu+1}t\right)h
 =-R_*t^{-\mu-2}+O(t^{-\mu-3}(\log t)^2)\,.
\end{equation*}
The additional term $-(\mu+1)(\mu+2)f/t^2$ from this substitution has the displayed error size. The stable integral gives $h=-(R_*/(2\beta))t^{-\mu-2}+O(t^{-\mu-3}(\log t)^2)$. Therefore
\begin{equation*}
 (t^{\mu+1}f)'=-\frac{R_*}{2\beta t}+O(t^{-2}(\log t)^2)\,.
\end{equation*}
The last error is integrable. There is a real constant $Ac_u$ and hence
\begin{equation*}
 f=t^{-\mu-1}\left(-\frac{R_*}{2\beta}\log t+Ac_u\right)
       +O(t^{-\mu-2}(\log t)^2)\,.
\end{equation*}
The same stable formulas control the time derivative. Combining this with the higher-mode estimate and $\kappa_*=-R_*/(2\beta A)$ proves~\eqref{eq:critsecond}, including its weighted first derivatives.
\end{proof}

\subsection{The physical Hessian and the critical endpoint}
Recall $H_0=\Hmat_{-\beta}[\psi]$. In~\eqref{eq:critH} the physical Hessian is represented by $\widehat H_u$ in the stated polar frame. The annular estimates are applied to the physical derivatives before this orthogonal change of basis.

\begin{proof}[Proof of Theorem~\ref{thm:critH}]
Apply~\eqref{eq:matrixidentity} to the two explicit terms of~\eqref{eq:critsecond}. The derivative of $t^{-\mu}$ contributes $+\mu\Dmat_{-\beta}[\psi]/t$ and must be kept. Its trace, together with the angular correction, is
\begin{equation*}
 \tr\big(\Hmat_{-\beta}[Z_*]+\mu\Dmat_{-\beta}[\psi]\big)
 =-A^{p-1}\psi^p\,;
\end{equation*}
This is the first nonzero Laplacian term and has the source sign. Let $U_{\mathrm{exp}}$ be the physical function given by the explicit terms of~\eqref{eq:critsecond}, and set $v_{\mathrm{exp}}=a^\mu r^\beta U_{\mathrm{exp}}\circ F$. The time derivative of the leading term and the angular correction give
\begin{align*}
 L v_{\mathrm{exp}}
 &=t^{-\mu-1}\bigl(-2\beta\mu A\psi-H_Q\chi\bigr)
       +O(t^{-\mu-2}\log t)\\
 &=-A^pt^{-\mu-1}\psi^p+O(t^{-\mu-2}\log t),
\end{align*}
using $H_Q\chi=A^pQ_0(\psi^p)$ and $A^{p-1}\Gamma_*=2\beta\mu$. Equation~\eqref{eq:critv} and the relative expansion of $v$ give the same leading term and error for $Lv$. Thus $E_u=u-U_{\mathrm{exp}}$ and $E_v=a^\mu r^\beta E_u\circ F$ satisfy $|E_v|+|LE_v|=O(t^{-\mu-2}(\log t)^2)$ by subtraction and~\eqref{eq:critsecond}. Conformal covariance now gives on enlarged physical annuli
\begin{equation*}
 \norm{E_u}_\infty+r^2\norm{\Delta_xE_u}_\infty
 \le Cr^{-\beta}t_r^{-\mu-2}(\log t_r)^2\,.
\end{equation*}
The error has zero true side traces. Apply~\eqref{eq:annestimate} directly to it to get its physical Hessian at that rate. For the explicit terms the quantitative part of Lemma~\ref{lem:physical} applies: their normalized profiles have bounded $C^{2,\eta}$ norms and the geometric error is a positive power of $r$, hence smaller than the displayed logarithmic remainder. This proves~\eqref{eq:critH}--\eqref{eq:critHerror} without a pointwise bound for $F''$.

For the exact tail write $b_0=\beta+2$, $D_0=C_*\sqrt2\beta(\beta+1)$, and choose $R$ by $L=D_0R^{-b_0}(\log(1/R))^{-\mu}$. On rescaled bounded annuli the leading norm tends to $|z|^{-b_0}$. The inner-area and outer-Chebyshev argument of Lemma~\ref{lem:tail} applies also to this slowly varying factor; for outer annuli its ratio is bounded by an arbitrarily small power, absorbed by choosing $\ell>S_0$. Since $\log L\sim b_0\log(1/R)$, radial integration gives~\eqref{eq:criticaltail}. The rearrangement therefore satisfies
\begin{equation*}
 (|D^2u|)^*(s)\asymp s^{-1/S_0}(\log(1/s))^{-\mu}\,.
\end{equation*}
The Lorentz integral with finite second index $q_0$ is finite exactly when $\mu q_0>1$. Since $\mu S_0=\beta/(\beta+2)<1$, the strong $L^{S_0}$ endpoint fails.

For two solutions the common angular and logarithmic matrices cancel. If $c_u\ne c_v$, the leading difference is $C_*(c_u-c_v)r^{-b_0}t^{-\mu-1}H_0$. Here $(\mu+1)S_0=1$, and the same calculation gives
\begin{equation*}
 \lim_{L\to\infty}L^{S_0}\log L\,
 \abs{\{|D^2(u-v)|>L\}\cap B_\delta(O)}
 =\frac{\theta b_0}2\big(D_0|c_u-c_v|\big)^{S_0}>0\,.
\end{equation*}
The endpoint integral diverges. If $c_u=c_v$, cover each strip $j<t<j+1$ by finitely many annuli of the type in~\eqref{eq:critHerror}. On them the measure factor $e^{-2j}$ cancels the $S_0$-th power of $e^{(\beta+2)j}$, since $(\beta+2)S_0=2$. H\"older's inequality on each rescaled annulus, from any finite $\ell>1$ to $S_0$, then bounds the $S_0$-power integral on the strip by $Cj^{-(\mu+2)S_0}(\log j)^{2S_0}$. Since $(\mu+2)S_0=1+S_0>1$, these terms are summable. This proves~\eqref{eq:criticaldifference}.
\end{proof}

\subsection{Global realization and separation of the critical branch}\label{subsec:criticalglobal}

We use a first-mode construction of the type used in~\cite[Section~2.1]{DMP}, with the corner chart kept fixed. The argument also controls the difference of two solutions at the order of the constant in~\eqref{eq:critsecond}.

\begin{proof}[Proof of Theorem~\ref{thm:criticalglobal}]
Put $t=\log(1/r)$, $\mu=1/(p-1)=\beta/2$ and $m(t,\varphi)=a|F'(re^{i\varphi})|^2$. For $V=r^\beta u\circ F$, the cylinder equation is
\[
 LV=-mV^p,\qquad L=\partial_t^2+2\beta\partial_t+\partial_\varphi^2+\beta^2.
\]
By~\eqref{eq:conformal}, $m=a+O(e^{-\eta_0t})$ uniformly up to the sides. Write $\Pi g=N^{-1}\int_I g\psi$, $\Pi^\perp g=g-(\Pi g)\psi$, and fix a starting section $T$ in the chart. Set
\begin{align*}
 B(t)&=\Pi(m\psi^p), & M_C(t)&=C+\frac{p-1}{2\beta}\int_T^t B(s)\dd s,\\
 A_C&=M_C^{-\mu}, & R_C&=M_C^{-1},\qquad \lambda_C=(\mu+1)M_C'/M_C.
\end{align*}
Thus $M_C'\asymp1$, $\sup R_C\le C^{-1}$ and $\int_T^\infty R_C^2\lesssim C^{-1}$. In this proof constants in estimates are uniform for large $C$ unless a dependence is indicated.

\emph{The local construction.}
Write $V=y\psi+W$, $\Pi W=0$, and $Z=y+y'/(2\beta)=A_C+z$. Impose $z(T)=0$ and $y(T)=Z(T)$. Recover $y$ from $y'+2\beta y=2\beta Z$ by the retarded formula. In the norms
\[
 \|z\|_c=\sup_{t\ge T}\left(\frac{|z|}{A_C}+\frac{|z'|}{R_CA_C}\right),
 \qquad \|W\|_\perp=\sup_{t\ge T,\,\varphi\in I}\frac{|W|}{R_CA_C\psi},
\]
this gives $|y-Z|\lesssim (1+\|z\|_c)R_CA_C$, with the analogous linear difference bound. Indeed $y-Z$ is the convolution of $-Z'$ with the retarded kernel $e^{-2\beta t}\one_{\{t\ge0\}}$, and the logarithmic derivative of the positive weight $R_CA_C$ is $O(C^{-1})$.

For the orthogonal equation extend the forcing by zero to the left of $T$ and invert each mode $n\ge2$ with kernel
\[
 -\frac{e^{-\beta s-n\beta|s|}}{2n\beta}.
\]
The kernels have weighted $L^1$ norms $O(n^{-2})$ with weight $e^{\beta |s|/2}$. Consequently if the forcing is bounded by a slowly varying weight $H(t)$, with $|(\log H)'|\le\beta/2$, the inverse is bounded by $C H(t)$ in absolute value, by summation in $n$. Local flat-side $W^{2,\ell}$ estimates with $\ell>2$ give the same bound for first derivatives on unit strip sections, including those meeting $T$, and the zero side trace gives a bound $C H(t)\psi$ for the function. The weight may be extended constantly to the left for these estimates. This is the whole-strip version of the mode inverse in Lemma~\ref{lem:highinverse}.

Use this inverse for $LW_{\rm new}=-\Pi^\perp(mV^p)$. On a ball $\|z\|_c\le\varepsilon$, $\|W\|_\perp\le M$ with $M$ fixed and $\varepsilon$ small, $V\asymp A_C\psi$ for large $C$. The forcing bound is $O(R_CA_C)$, and the power mean value formula gives
\[
 \|W_{\rm new}\|_\perp\le C_1,
 \qquad \|\delta W_{\rm new}\|_\perp\le C_1(\|\delta z\|_c+C^{-1}\|\delta W\|_\perp).
\]
For the center solve $z_{\rm new}'+\lambda_C z_{\rm new}=g$, $z_{\rm new}(T)=0$, with
\[
 g=-\frac{\Pi(mV^p)}{2\beta}+\frac{B A_C^p}{2\beta}+\lambda_C z.
\]
Its inverse is $g\mapsto M_C^{-\mu-1}\int_T^t M_C(s)^{\mu+1}g(s)\dd s$ and maps the bound $E R_CA_C$ to $\|z_{\rm new}\|_c\lesssim E$, since $M_C^{\mu+1}R_CA_C=1$ and $t-T\lesssim M_C(t)$. The linear term in $z$ cancels because $pBA_C^{p-1}/(2\beta)=\lambda_C$. To estimate the remainder, expand $V^p$ about $A_C\psi$, dividing both arguments by $A_C\psi$ before Taylor expansion. The relative increment is $z/A_C+O((1+M)R_C)$ and stays in a fixed interval about zero. The quadratic Taylor estimate and the difference estimate for $y-Z$ therefore give
\[
 \|z_{\rm new}\|_c\le C_M(\varepsilon^2+C^{-1}),\qquad
 \|\delta z_{\rm new}\|_c\le C_M\big((\varepsilon+C^{-1})\|\delta z\|_c+C^{-1}\|\delta W\|_\perp\big).
\]
Choose $M$ large, $\varepsilon$ small and then $C$ large. The two maps preserve the ball and contract in the sum of the center norm and a sufficiently small fixed multiple of the orthogonal norm. The fixed point $V_C$ solves the cylinder equation. Its compensated moment $Z_C$ obeys
\[
 V_C=Z_C\psi(1+O(R_C)),\qquad
 Z_C'=-\Pi(mV_C^p)/(2\beta),\qquad Z_C(T)=C^{-\mu}.
\]
Integration gives
\begin{equation}\label{eq:criticalglobalmoment}
 Z_C^{1-p}=M_C+O(E_C),\qquad E_C(t)=\int_T^t R_C(s)\dd s .
\end{equation}
In particular $V_C\sim A_C\psi$ uniformly in relative form.

\emph{Correction on the given domain.}
Cut off $e^{\beta t}V_C$ by a nonnegative smooth time cutoff supported in the chart and equal to one on a fixed deep tail, and call the physical function $U_C$. Its defect $f_C=-\Delta U_C-aU_C^p$ is supported on the fixed transition annulus and satisfies $\|f_C\|_{H^{-1}(\Om)}\lesssim C^{-\mu}$. This follows from the cylinder equation, the derivative bounds above, and the bounded chart factors on that annulus. In conformal strip coordinates the weight in $\int U_C^{p-1}\xi^2$ is $r^2|F'|^2 U_C^{p-1}\lesssim C^{-1}$, since $\beta(p-1)=2$. Angular Poincar\'e and conformal invariance thus give
\[
 \int_{\Om} U_C^{p-1}\xi^2\lesssim C^{-1}\int_{\Om}|\nabla\xi|^2 .
\]
On the Dirichlet space $H^1_0(\Om)$ solve
\[
 w_C=(-\Delta_D)^{-1}\left(a\big((U_C+w_C)_+^p-U_C^p\big)-f_C\right).
\]
The preceding form bound and the embedding into $L^{p+1}$ bound the Lipschitz constant of the right side on a ball of radius $A C^{-\mu}$ by $O(C^{-1})$, for fixed $A$. The defect bound gives invariance for large $A$, hence a fixed point. Off $O$, a cutoff of $u_C=U_C+w_C$ is in $H^1_0$ and has Laplacian in $L^2$ by planar Sobolev embedding. Its Dirichlet Green potential is continuous with zero boundary values, by the uniform $L^2$ bound for the planar logarithmic Green kernel and approximation by bounded forcings. Thus $u_C$ has the ordinary zero boundary data.

On the tail write $w_C$ also for its pullback. It has finite strip Dirichlet energy and zero side traces, so its $H^1$ norm on unit strips is uniformly bounded by angular Poincar\'e. There the equation gives $|\Delta_{t,\varphi}w_C|\lesssim R_C|w_C|+e^{-2t}|w_C|^p$. Local Sobolev embedding and flat-side Poisson estimates give a uniform $C^1$ bound on smaller strips, hence $w_C=O(\psi)$. This is negligible relative to $e^{\beta t}A_C\psi$. The solution is positive on the deep tail, and superharmonic comparison on the remaining domain and the strong maximum principle give positivity throughout. Its classical interior regularity follows from the equation. The relative profile agrees with~\eqref{eq:criticalfirst} because
\[
 M_C(t)/t\longrightarrow b_\infty=\frac{a}{\beta^2N}\int_I\psi^{p+1},
 \qquad b_\infty^{-\mu}=C_*.
\]
Theorem~\ref{thm:sourcefirst} gives zero Martin atom. The same strip form bound for $u_C^{p-1}$ tends to zero as the tail is moved inward, giving stability near $O$. On the rest of the domain the linearized potential is bounded. A cutoff partition of the Dirichlet form thus gives a lower bound $c\|\nabla\xi\|_2^2-C'\|\xi\|_2^2$ for some $c>0$, which implies finite Morse index by Rellich compactness.

\emph{Separation at the phase scale.}
Take $C<D\le(1+\eta)C$, with $\eta>0$ to be fixed small, and put $R=R_C$. The difference of the local equations is $L(V_C-V_D)=-b(V_C-V_D)$, where
\begin{align*}
 b(t,\varphi)&=p m(t,\varphi)\int_0^1(sV_C+(1-s)V_D)^{p-1}\dd s,\\
 \Lambda&=\frac{\Pi(b\psi)}{2\beta},\qquad
 j(t)=\exp\left(-\int_T^t\Lambda(s)\dd s\right).
\end{align*}
Here $0\le b\lesssim R$. Using~\eqref{eq:criticalglobalmoment} and $M_D=M_C+D-C$, one obtains
\[
 |\Lambda-\lambda_C|\lesssim R\left(R+\frac{E_C}{M_C}+\frac{D-C}{M_C}\right).
\]
Indeed the relative deviations of $V_C,V_D$ from $A_C\psi$ have the three orders in parentheses, and the mean value bound for the power of the relative profiles gives the displayed estimate. This difference is bounded in absolute value by $C_2(C^{-1}+\eta)R$ and has $L^1(T,\infty)$ norm $O(C^{-1}+(D-C)/C)$. For the middle term use $M_C'\asymp1$ and Fubini to get $\int_T^\infty R^2 E_C\lesssim\int_T^\infty R^2$. It follows that
\[
 j_\infty:=\lim_{t\to\infty}j(t)M_C(t)^{\mu+1}\asymp C^{\mu+1},
\]
where the limit is positive and the comparison is uniform in the indicated pairs.

We give the comparison with the actual difference. Put $\delta_0=C^{-\mu}-D^{-\mu}>0$ and seek its compensated moment as $Z=\delta_0j+z$. Use the preceding norms with $A_C$ replaced by $j$, retaining $R$ in the derivative and orthogonal norms. Recover $y$ from $y'+2\beta y=2\beta Z$, $y(T)=Z(T)$, and use the same zero-extended inverse for $LW_{\rm new}=-\Pi^\perp(b(y\psi+W))$. Since $|j'/j|\lesssim R$ and $|R'/R|\lesssim R$, the retarded error is $|y-Z|\lesssim Rj(|\delta_0|+\|z\|_c)$. The orthogonal output has norm at most $C_3(|\delta_0|+\|z\|_c+C^{-1}\|W\|_\perp)$. The center output is
\[
 z_{\rm new}(t)=j(t)\int_T^t
 \frac{-\Lambda(y-Z)-\Pi(bW)/(2\beta)}{j(s)}\dd s .
\]
The integrand is bounded by $C_3 R^2(|\delta_0|+\|z\|_c+\|W\|_\perp)$, so the center norm of the output is bounded by $C_4 C^{-1}$ times the same sum, including the derivative bound obtained from its equation. These linear estimates give an affine contraction in a weighted sum of the norms. Its solution has both norms $O(\delta_0)$ and $Z/j\to\delta_0+O(\delta_0/C)>0$.

This solution is the actual difference. Both have the same initial compensated moment, the same retarded condition, and the same zero-extended orthogonal inverse. The actual difference lies in the original $A_C$-weighted spaces because $A_D\le A_C$ and $R_D\le R_C$. The solution constructed with weight $j$ also lies in those spaces, by $j=O(M_C^{-\mu-1})$ and $|j'|\lesssim R_Cj$ for the fixed pair. In those spaces its center forcing for $\partial_t+\lambda_C$ is $(\lambda_C-\Lambda)Z-\Lambda(y-Z)-\Pi(bW)/(2\beta)$, with center norm estimate $O((C^{-1}+\eta)\|Z\|_c+C^{-1}\|W\|_\perp)$. The orthogonal norm estimate is $O(\|Z\|_c+C^{-1}\|W\|_\perp)$. Taking $\eta$ small and $C_0$ large gives uniqueness with zero initial moment. Since both $y-Z$ and $W/\psi$ are $O(\delta_0 Rj)$, we conclude
\[
 \frac{V_C-V_D}{M_C^{-\mu-1}\psi}\longrightarrow \ell_{C,D}>0,
 \qquad \ell_{C,D}\asymp C^{\mu+1}(C^{-\mu}-D^{-\mu})\asymp D-C .
\]
The estimates on $w_C,w_D$ show that replacing the local profiles by $r^\beta u_C\circ F,r^\beta u_D\circ F$ preserves this limit. Subtract~\eqref{eq:critsecond} for the two solutions. The universal terms cancel, so comparison of the $t^{-\mu-1}\psi$ coefficients gives $C_*(c_{u_C}-c_{u_D})=b_\infty^{-\mu-1}\ell_{C,D}$. The constants are uniform for nearby parameters, and a finite chain of such pairs proves~\eqref{eq:criticalglobalphase} for every $C<D$. The Hessian assertion follows from~\eqref{eq:criticaldifference}.
\end{proof}

\begin{remark}\label{rem:noexhaustion}
The constant $c_u$ is an asymptotic invariant. Neither~\eqref{eq:critsecond} nor~\eqref{eq:criticaldifference} proves that two solutions with the same invariant are identical, or that every critical singular solution belongs to any particular globally constructed family. \citationadded{For smooth bounded domains,~\cite[Theorem~1.1]{DMP} constructs singular solutions at and slightly above the boundary-critical exponent. For a domain contained in a cone and agreeing with it near the vertex,~\cite[Theorem~1]{HMR} gives two singular solutions for $p\in(p_*,p_*+\varepsilon)$, where $p_*$ is its cone threshold. The cone and higher-dimensional wedge constructions in~\cite[Theorems~1.1--1.3]{Gkikas} have their stated geometry assumptions and strict supercritical exponent intervals.} The local classification and the refinements here do not require a global family to exhaust the natural class.
\end{remark}

\section{The existence range}\label{sec:range}
The preceding results now provide the proof of Theorem~\ref{thm:existencetable}. Here bounded solutions at $O$ are allowed, as in Definition~\ref{def:natural}.

\begin{proof}[Proof of Theorem~\ref{thm:existencetable}]
For negative powers and sublinear powers the positive constructions were given in Sections~\ref{sec:negabs}--\ref{sec:sublinear}. The superlinear absorption range is Proposition~\ref{prop:superabs}; the superlinear source has a bounded variational solution for every finite power by Proposition~\ref{prop:supersourceexist}. Proposition~\ref{prop:elementary} covers zero and first powers.

It remains to exclude $k\le-1$ for absorption. Put $k=-p$, $p\ge1$,
and flatten a small ordinary boundary patch by a conformal chart with
bounded and boundedly invertible derivative. After rescaling, the
pullback $v$ is positive on the upper half-disc $D$ of radius one,
vanishes on the diameter, and satisfies $\Delta v=A(z)v^{-p}$ with
$A\ge a_0>0$. Let $h$ be the harmonic lifting of its boundary data.
Subharmonic comparison gives $v\le h$, and near the middle of the
diameter $h(z)\le C\operatorname{Im}z$. For compact truncations
$f_j\uparrow A v^{-p}$ the same comparison gives $v+\G_D[f_j]\le h$, where $\G_D$ is the Dirichlet Green operator on $D$.

Write $x=is$ and $z=\xi+iy$. For $2s<y<y_0$ and $|\xi|<y$, the
half-plane Green kernel is bounded below by $c s/y$, since it equals
$(4\pi)^{-1}\log(1+4sy/|x-z|^2)$. The harmonic correction for the
semicircular boundary is at most $Csy$. Indeed, for $z$ near zero,
the half-plane kernel on that semicircle is bounded by
$C\operatorname{Im}(x)\operatorname{Im}(z)$, and this bound extends
by harmonic comparison. Taking $y_0$ small therefore gives
$G_D(is,z)\ge c' s/y$ in the indicated region. Monotone convergence
in the truncated comparison and integration over that region yield
\[
 Cs\ge h(is)\ge c' s\int_{2s}^{y_0} y^{-p}\dd y .
\]
The integral tends to infinity as $s\downarrow0$ for every $p\ge1$.
Thus no positive solution exists for $k\le-1$, independently of
the corner angle.
\end{proof}

\section*{Funding}
W.~Zou acknowledges support from the National Key R\&D Program of China under grant 2023YFA1010001 and from the National Natural Science Foundation of China under grant 12171265.

\section*{Conflict of interest}
The authors declare that they have no known competing financial interests or personal relationships that could have appeared to influence the work reported in this paper.

\section*{Data availability}
No data were used for the research described in this article.

\section*{Declaration on the use of artificial intelligence}
During research discussions and preparation of the manuscript, generative artificial intelligence was used to assist with literature searches and reference checks, discussion and elaboration of mathematical arguments, language revision, and \LaTeX{} typesetting. The authors assume full responsibility for the mathematical validity, accuracy, and integrity of the proofs presented in the manuscript.

\end{document}